\documentclass[secnum]{rims27}
\usepackage[matrix,arrow,curve]{xy}
\usepackage{euscript, amssymb, amscd}
\makeatletter 
\@addtoreset{equation}{subsection} 
\@addtoreset{footnote}{section}
\makeatother

\renewcommand{\theequation}{\thesubsection.\arabic{equation}}
\newcommand{\type}[1]{$\mathrm{(#1)}$}

\newcommand{\operatorname}[1]{\mathop{\rm #1}\nolimits}
\newcommand{\ord}{\operatorname{ord}}
\newcommand{\wt}{\operatorname{wt}}
\newcommand{\rk}{\operatorname{rk}}
\newcommand{\ow}{\operatorname{ow}}
\newcommand{\lcm}{\operatorname{lcm}}
\newcommand{\muu}{\mbox{\boldmath $\mu$}}

\newcommand{\Coker}{\operatorname{Coker}}
\newcommand{\Clsc}{\operatorname{Cl^{\operatorname{sc}}}}
\newcommand{\gr}{\operatorname{gr}}
\newcommand{\red}{\operatorname{red}}
\newcommand{\len}{\operatorname{len}}
\newcommand{\ii}{\operatorname{i}}
\newcommand{\mod}{\operatorname{mod}}
\newcommand{\size}{\operatorname{siz}}
\newcommand{\Spec}{\operatorname{Spec}}

\newcommand{\qldeg}{ql \deg}

\newcommand{\Supp}{\operatorname{Supp}}
\newcommand{\Sing}{\operatorname{Sing}}
\newcommand{\Pic}{\operatorname{Pic}}
\newcommand{\mt}[1]{\operatorname{#1}}
\newcommand{\down}[1]{\left\lfloor #1\right\rfloor}

\newcommand{\HHom}{\operatorname{\EuScript{H}\it{om}}}

\newcommand{\MMM}{{\EuScript{M}}}
\newcommand{\LLL}{{\EuScript{L}}}
\newcommand{\EEE}{{\EuScript{E}}}
\newcommand{\OOO}{\mathcal{O}}
\newcommand{\FFF}{\mathcal{F}}
\newcommand{\CC}{\mathbb{C}}
\newcommand{\ZZ}{\mathbb{Z}}
\newcommand{\PP}{\mathbb{P}}
\newcommand{\QQ}{\mathbb{Q}}

\newcommand{\qq}{\mathbin{\sim_{\scriptscriptstyle{\QQ}}}}
\newcommand{\toplus}{\mathbin{\tilde\oplus}}
\newcommand{\totimes}{\mathbin{\tilde\otimes}}
\newcommand{\xref}[1]{{\rm \ref{#1}}}
\newcommand{\eqref}[1]{{\rm (\ref{#1})}}
\newcommand{\comp}{\mathbin{\scriptstyle{\circ}}}

\newcommand{\comment}[1]{}

\newcounter{THN}[section]
\renewcommand{\theTHN}
{(\arabic{section}.\arabic{subsection})}
\newcounter{THNO}[section]
\renewcommand{\theTHNO}
{(\arabic{section}.\arabic{subsection}.\arabic{equation})}

\newenvironment{mparag}[1]{
\setcounter{THN}{\value{subsection}}
\refstepcounter{subsection}\refstepcounter{THN}
\par\medskip\noindent\begingroup \rm
{\bf\theTHN\ #1\ }}{\par\smallskip\endgroup}

\newenvironment{mtparag}[1]{
\setcounter{THN}{\value{subsection}}
\refstepcounter{subsection}\refstepcounter{THN}
\par\medskip\noindent\begingroup \it
{\bf\theTHN\ #1\ }}{\par\smallskip\endgroup}

\newenvironment{parag}[1]{
\setcounter{THNO}{\value{equation}}
\refstepcounter{equation}\refstepcounter{THNO}
\par\medskip\noindent\begingroup \rm
{\bf\theTHNO\ #1\ }}{\par\smallskip\endgroup}

\newenvironment{tparag}[1]{
\setcounter{THNO}{\value{equation}}
\refstepcounter{equation}\refstepcounter{THNO}
\par\medskip\noindent\begingroup \it
{\bf\theTHNO\ #1\ }}{\par\smallskip\endgroup}

\newenvironment{example}{\begin{parag}{Example.}}{\end{parag}}
\newenvironment{remark}{\begin{parag}{Remark.}}{\end{parag}}
\newenvironment{pusto}{\begin{parag}{}}{\end{parag}}
\newenvironment{proposition}{\begin{tparag}{Proposition.}}{\end{tparag}}
\newenvironment{lemma}{\begin{tparag}{Lemma.}}{\end{tparag}}
\newenvironment{claim}{\begin{tparag}{Claim.}}{\end{tparag}}
\newenvironment{theorem}{\begin{tparag}{Theorem.}}{\end{tparag}}
\newenvironment{emptytheorem}{\begin{tparag}{}}{\end{tparag}}
\newenvironment{corollary}{\begin{tparag}{Corollary.}}{\end{tparag}}

\newenvironment{theoremm}{\begin{mtparag}{Theorem.}}{\end{mtparag}}
\newenvironment{lemmam}{\begin{mtparag}{Lemma.}}{\end{mtparag}}
\newenvironment{propositionm}{\begin{mtparag}{Proposition.}}{\end{mtparag}}
\newenvironment{definitionm}{\begin{mparag}{Definition.}}{\end{mparag}}

\title{On $\QQ$-conic bundles}
\author{Shigefumi 
\textsc{Mori}\footnote{RIMS, 
Kyoto University, Oiwake-cho, Kitashirakawa, Sakyo-ku, Kyoto
606-8502, Japan.
\newline
e-mail: \texttt{mori@kurims.kyoto-u.ac.jp}
} \hspace{2pt} and Yuri 
\textsc{Prokhorov}\footnote{Department 
of Algebra, Faculty of Mathematics, Moscow State
University, Moscow 117234, Russia.
\newline
e-mail: \texttt{prokhoro@mech.math.msu.su}
}} 
\AuthorHead{Shigefumi Mori and Yuri Prokhorov}

\VolumeNo{4x}  
\YearNo{200x}  
\PagesNo{000--000}

\begin{document}
\maketitle
\begin{abstract}
A $\mathbb Q$-conic bundle is a proper morphism from a threefold
with only terminal singularities to a normal surface such that
fibers are connected and the anti-canonical divisor is
relatively ample. We study the structure of $\mathbb Q$-conic
bundles near their singular fibers.
One corollary to our main results is that the base surface of
every $\mathbb Q$-conic bundle has only Du Val singularities of
type A (a positive solution of a conjecture by Iskovskikh).
We obtain the complete classification of $\mathbb Q$-conic
bundles under the additional assumption that the singular fiber
is irreducible and the base surface is singular.
\end{abstract} 
\tableofcontents

\section{Introduction}
In this paper we study the local structure of extremal contractions
from threefolds to surfaces. Such contractions naturally appear in
the birational classification of three-dimensional algebraic
varieties of negative Kodaira dimension. More precisely, according
to the minimal model program every algebraic projective threefold
$V$ with $\kappa(V)=-\infty$ is birationally equivalent to a
$\QQ$-factorial terminal threefold $X$ having a $K_X$-negative
extremal contraction to a lower dimensional variety $Z$. There are
three cases:

a) $Z$ is a point and then $X$ is a $\QQ$-Fano variety with
$\rho(X)=1$,

b) $Z$ is a smooth curve and then $X/Z$ is a del Pezzo fibration,

c) $Z$ is a normal surface and then $X/Z$ is a rational curve
fibration.

We study the last case.
\begin{definitionm}
By a \textit{$\QQ$-conic bundle} we mean a projective morphism
$f\colon X \to Z$ from an (algebraic or analytic) threefold to a
surface that satisfies the following properties:
\begin{enumerate}
\item
$X$ is normal and has only terminal singularities,
\item
$f_*\OOO_X=\OOO_Z$,
\item
all fibers are one-dimensional,
\item
$-K_X$ is $f$-ample.
\end{enumerate}
For $f\colon X\to Z$ as above and for a point $o\in Z$, we call the
\textit{analytic} 
germ $(X, f^{-1}(o)_{\red})$ a \textit{$\QQ$-conic
bundle germ}.
\end{definitionm}

The easiest example of $\QQ$-conic bundles is a \textit{standard
Gorenstein conic bundle}: $Z$ is smooth and $X$ is embedded in the
projectivization $\PP_Z(\EEE)$ of a rank $3$ vector bundle so that
the fibers $X_z$, $z\in Z$ are conics in $\PP_Z(\EEE)_z$. More
complicated examples can be constructed as quotients:

\begin{parag}{Example-Definition {(toroidal example)}.}
\label{ex-toric}
Consider the following action of $\muu_m$ on $\PP^1_x\times
\CC^2_{u,v}$:
\[
(x;u,v) \longmapsto(\varepsilon x; \varepsilon^a u,
\varepsilon^{b} v),
\]
where $\varepsilon$ is a primitive $m$-th root of unity and $\gcd
(m,a)=\gcd(m,b)=1$. Let $X:=\PP^1\times\CC^2/\muu_m$,
$Z:=\CC^2/\muu_m$ and let $f\colon X\to Z$ be the natural
projection. Since $\muu_m$ acts freely in codimension one, $-K_X$ is
$f$-ample. Two fixed points on $\PP^1\times \CC^2$ gives two cyclic
quotient singularities of types $\frac1m(1,a,b)$ and
$\frac1m(-1,a,b)$ on $X$. These points are terminal if and only if
$a+b\equiv 0 \mod m$. In this case, $f$ is a $\QQ$-conic bundle and
the base surface $Z$ has a Du Val singularity of type $A_{m-1}$. We
say that a $\QQ$-conic bundle germ is \textit{toroidal} if it is
biholomorphic to $f\colon (X,f^{-1}(0)_{\red}) \to (Z,0)$ above
(with $a+b\equiv 0 \mod m$).
\end{parag}

Our first main result is a complete classification of $\QQ$-conic
bundle germs with irreducible central fiber under the assumption
that the base surface is singular:

\begin{theoremm}
\label{th-main-1}
Let $f\colon (X,C)\to (Z,o)$ be a $\QQ$-conic bundle germ, where $C$
is irreducible. Assume that $(Z,o)$ is singular. Then one of the
following holds:

\textbf{Cases where $X$ is locally primitive.}

\begin{emptytheorem}
\label{item=main--th-pr-toric}
$(X,C)$ is toroidal.
\end{emptytheorem}
\begin{emptytheorem}
\label{item=main--th-pr-ex3}
$(X,C)$ is biholomorphic to the quotient of the smooth $\QQ$-conic
bundle
\[
X'= \{ y_1^2+uy_2^2+vy_3^2=0\}\subset \PP^2_{y_1,y_2,y_3} \times\CC^2_{u,v}
\longrightarrow \CC^2_{u,v}.
\]
by $\muu_{m}$-action:
\[
(y_1,y_2,y_3,u,v)\longmapsto (\varepsilon^{a} y_1,\varepsilon^{-1}y_2,y_3,\varepsilon
u,\varepsilon^{-1} v).
\]
Here $m=2a+1$ is odd and $\varepsilon$ is a primitive $m$-th root of
unity. The singular locus of $X$ consists of two cyclic quotient
singularities of types $\frac 1m(a,-1,1)$ and $\frac 1m(a+1,1,-1)$.
The base surface $\CC^2/\muu_m$ has a singularity of type $A_{m-1}$.
\end{emptytheorem}

\textbf{Cases where $X$ is not locally primitive.} 
Let $P\in X$ be the imprimitive
point and let $m$, $s$ and $\bar m$ be its index, splitting degree
and subindex, respectively. In this case, $P$ is the only
non-Gorenstein point and $X$ has at most one Gorenstein singular
point. We refer to \xref{prop-imp-types} for the definition of types
\type{IA^\vee} -- \type{IE^\vee}.

\begin{emptytheorem}
\label{item-main-th-impr-barm=2-s=4}
$(X,C)$ is of type \type{IE^\vee} at $P$, $s=4$, $\bar m=2$, $(Z,o)$
is Du Val of type $A_3$, $X$ has a cyclic quotient singularity $P$
of type $\frac18(5,1,3)$ and has no other singular points.
Furthermore, $(X,C)$ is the quotient of the index-two $\QQ$-conic
bundle germ given by the following two equations in
$\PP(1,1,1,2)_{y_1,\dots,y_4}\times \CC^2_{u,v}$
\[
\label{eq-imp-exc-8-eq-a}
\begin{array}{lll}
y_1^2-y_2^2&=&u \psi_1(y_1,\dots,y_4;u,v)+v\psi_2(y_1,\dots,y_4;u,v),
\\[5pt]
y_1y_2-y_3^2&=&u \psi_3(y_1,\dots,y_4;u,v)+v\psi_4(y_1,\dots,y_4;u,v)
\end{array}
\]
by $\muu_{4}$-action:
\[
y_1\mapsto -\ii y_1,\quad y_2\mapsto \ii y_2,\quad y_3\mapsto - y_3,\quad y_4\mapsto
\ii y_4,\quad u\mapsto \ii u,\quad v\mapsto -\ii v,
\]
\textup(see Example \xref{example-imp-4}\textup).
\end{emptytheorem}

\begin{emptytheorem}
\label{item-main-th-impr-barm=1}
$(X,C)$ is of type \type{ID^\vee} at $P$, $s=2$, $\bar m=1$, $(Z,o)$
is Du Val of type $A_1$, $(X,C)$ is a quotient of a Gorenstein conic
bundle given by the following equation in $\PP^2_{y_1,y_2,y_3}\times
\CC^2_{u,v}$
\[
y_1^2+y_2^2+\psi(u,v)y_3^2=0, \qquad \psi(u,v)\in\CC\{u^2,\, v^2,\, uv\},
\]
by $\muu_{2}$-action:
\[
u\mapsto -u,\quad v\mapsto -v,\quad y_1\mapsto -y_1,\quad y_2\mapsto y_2,\quad y_3\mapsto y_3.
\]
Here $\psi(u,v)$ has no multiple factors. In this case, $(X,P)$ is
the only singular point and it is of type $cA/2$ or $cAx/2$.
\end{emptytheorem}

\begin{emptytheorem}
\label{item-main-th-impr-barm=2-s=2-cycl}
$(X,C)$ is of type \type{IA^\vee} at $P$, $\bar m=2$, $s=2$, $(Z,o)$
is Du Val of type $A_1$, $(X,P)$ is a cyclic quotient singularity
of type $\frac{1}{4}(1,1,3)$,
and $(X,C)$ is the quotient of the index-two $\QQ$-conic bundle germ
given by the following two equations in
$\PP(1,1,1,2)_{y_1,\dots,y_4}\times \CC^2_{u,v}$
\[
\begin{array}{lll}
y_1^2-y_2^2&=&u \psi_1(y_1,\dots,y_4;u,v)+v\psi_2(y_1,\dots,y_4;u,v),
\\[5pt]
y_3^2&=&u \psi_3(y_1,\dots,y_4;u,v)+v\psi_4(y_1,\dots,y_4;u,v)
\end{array}
\]
by $\muu_{2}$-action:
\[
y_1\mapsto y_1,\quad y_2\mapsto - y_2,\quad y_3\mapsto y_3,\quad y_4\mapsto -
y_4,\quad u\mapsto - u,\quad v\mapsto - v
\]
\textup(see Example \xref{example-imp-IA}\textup).
\end{emptytheorem}

\begin{emptytheorem}
\label{item-main-th-impr-barm=2-s=2-cAx/4}
$(X,C)$ is of type \type{II^\vee} at $P$, $\bar m=2$, $s=2$, $(Z,o)$
is Du Val of type $A_1$, $(X,P)$ is a singularity of type $cAx/4$,
and $(X,C)$ is the quotient of the same form as in
\xref{item-main-th-impr-barm=2-s=2-cycl} \textup(see Example
\xref{example-imp-II}\textup).
\end{emptytheorem}
\par\medskip\noindent
All the cases \xref{item=main--th-pr-toric} --
\xref{item-main-th-impr-barm=2-s=2-cAx/4} occur.
\end{theoremm}

By running MMP over the base $Z$ we immediately 
obtain the following fact which was
conjectured by Iskovskikh:
\begin{theorem}
\label{theo-main-DuVal}
Let $f\colon X \to Z$ be a $\QQ$-conic bundle \textup(possibly with
reducible fibers\textup). Then $Z$ has only Du Val singularities of
type $A$.
\end{theorem}

We also note that the singularity $(Z,o)$ is unbounded only in
locally primitive cases \ref{item=main--th-pr-toric} and
\ref{item=main--th-pr-ex3}. In all other cases $(Z,o)$ is either of
type $A_1$ or $A_3$. Theorem \ref{theo-main-DuVal} has important
applications to rationality problem of conic bundles
\cite{Iskovskikh-1996-conic-re}.

Note that the condition that $X$ has only terminal singularities is
essential in Theorem \ref{theo-main-DuVal}: put in Example
\ref{ex-toric} $a=1$ and $b=-2$ ($m$ is odd). We get an extremal
contraction having two singular points which are canonical
Gorenstein of type $\frac1m(1,1,-2)$ and terminal of type
$\frac1m(-1,1,-2)$. The base surface has a singularity of type
$\frac 1m (1,-2)$ which is not Du Val.

\begin{corollary}
\label{cor-main-sing-s}
If in notation of \xref{th-main-1} the base $(Z,o)$ is not of type
$A_1$, then $X$ has only cyclic quotient singularities.
\end{corollary}

In the case of smooth base our results are not so strong:
\begin{theoremm}
Let $(X,C\simeq \PP^1)$ be a $\QQ$-conic bundle germ over a smooth
base $(Z,o)$. Then $(X,C)$ is locally primitive and the
configuration of singular points is one of the following
\textup(notation \type{IA} -- \type{III} are explained in
\xref{prop-prim-types-def}\textup)\textup:
\begin{emptytheorem}
\label{item=main--th2-Gor}
$\emptyset$,\ \type{III},\ \type{III}$+$\type{III} \quad \textup($X$
is Gorenstein\textup).
\end{emptytheorem}
\begin{emptytheorem}
\label{item=main--th2-IA}
\type{IA}, \ \type{IA}$+$\type{III},\
\type{IA}$+$\type{III}$+$\type{III}.
\end{emptytheorem}
\begin{emptytheorem}
\label{item=main--th2-II}
\type{IIA}, \ \type{IIA}$+$\type{III}.
\end{emptytheorem}
\begin{emptytheorem}
\label{item=main--th2-ICIB}
\type{IC},\ \type{IIB}.
\end{emptytheorem}
\begin{emptytheorem}
\label{item=main--th2-IAIA}
\type{IA}$+$\type{IA} of indices $2$ and odd $m\ge 3$.
\end{emptytheorem}
\begin{emptytheorem}
\label{item=main--th2-IAIAIII}
\type{IA}$+$\type{IA}$+$\type{III} of indices $2$, odd $m\ge 3$ and
$1$.
\end{emptytheorem}
\end{theoremm}

In contrast with Theorem \ref{th-main-1} we can say nothing about
the existence of $\QQ$-conic bundles as in \ref{item=main--th2-II}
-- \ref{item=main--th2-IAIAIII}. There are examples of index-two
$\QQ$-conic bundles as in \ref{item=main--th2-IA} (see \cite[\S
3]{Prokhorov-1997_e} and \ref{th-index=2}). One can also easily
construct examples of Gorenstein standard conic bundles of type
\ref{item=main--th2-Gor}.

\begin{tparag}{Proposition
(Reid's conjecture about general elephant.}
\label{propo-main-ge}
Let $(X,C\simeq \PP^1)$ be a $\QQ$-conic bundle germ. Then, except
possibly for cases \xref{item=main--th2-ICIB},
\xref{item=main--th2-IAIA}, and \xref{item=main--th2-IAIAIII}, a
general member of $|-K_X|$ has only Du Val singularities. In these
exceptional cases a general member of $D\in |-2K_X|$ does not contain
$C$ and the log divisor $K_X+\frac 12 D$ is log terminal.
\end{tparag}
Proposition \ref{propo-main-ge} follows from
Remarks
\ref{rem-imp-ge}, \ref{rem-prim-ge} and Theorem \ref{th-semist-ge}.

\begin{mparag}{Comments on the approach.}
Contractions similar to $\QQ$-conic bundles were considered in
\cite{Mori-1988}. In fact, \cite{Mori-1988} deals with birational
contractions of threefolds $f\colon X\to (Z,o)$ such that $X$ has
only terminal singularities, $-K_X$ is $f$-ample, and
$C:=f^{-1}(o)_{\red}$ is a curve. In this case, we have vanishings
$R^1f_*\OOO_X=0$ and $R^1f_*\omega_X=0$. (cf. \ref{th-vanish}).
Though the former vanishing was used all over the places in
\cite{Mori-1988}, the latter vanishing $R^1f_*\omega_X=0$ was used
only occasionally. It is easy to find the places where the
corollaries of $R^1f_*\omega_X=0$ were used. In this paper we follow
the arguments of \cite{Mori-1988} paying special attention to those
corollaries of $R^1f_*\omega_X=0$ and furthermore give comments to
modify the arguments when the corollaries are used.

Though fewer vanishing conditions are available, we have new tools 
Lemma \ref{lemma-KC} and Theorem \ref{theo-sect-1-gen-fiber-m}
for $\QQ$-conic bundles. These results together with
\cite{Mori-1988} form the basis of our approach.
\end{mparag}

\section{Preliminaries}
\setcounter{subsection}{1}
\begin{mparag}{}
Let $f\colon (X,C)\to (Z,o)$ be a $\QQ$-conic bundle germ. The
following is an immediate consequence of the Kawamata-Viehweg
vanishing theorem.
\end{mparag}

\begin{theoremm}
\label{th-vanish}
$R^if_*\OOO_X=0$ for $i>0$.
\end{theoremm}

\begin{tparag}{Corollary
(cf. {\cite[Remark 1.2.1, Cor. 1.3]{Mori-1988}}).}
\label{cor-C-pa=0}
\begin{enumerate}
\item
If $J$ be an ideal such that $\Supp \OOO_X/J\subset C$, then
$H^1(\OOO_X/J)=0$.
\item
$p_a(C)=0$ and $C$ is a union of smooth rational curves.
\item
$\Pic X\simeq H^2(C,\ZZ)\simeq \ZZ^\rho$, where $\rho$ is the number
of irreducible components of $C$.
\end{enumerate}
\end{tparag}

\begin{remark}
\label{rem-prel-extr-nbd}
If $C$ is reducible, then $\rho(X/Z)>1$ and for every closed curve
$C'\subsetneq C$ the germ $(X,C')$ is an extremal neighborhood
\textup(isolated or divisorial\textup). These were classified in
\cite{Mori-1988} and \cite{Kollar-Mori-1992} under the condition
that $C'$ is irreducible.
\end{remark}

\begin{remark}
In general, we do not assume that $X$ is $\QQ$-factorial (i.e., a
Weil divisor on $X$ is not necessarily $\QQ$-Cartier). In fact, the
following are equivalent
\begin{enumerate}
\item
$X$ is $\QQ$-factorial and $\rho (X/Z)=1$,
\item
the preimage of an arbitrary irreducible curve $\Gamma\subset Z$ is
also irreducible.
\end{enumerate}
Indeed, the implication (i) $\Rightarrow$ (ii) is obvious. To show
(ii) $\Rightarrow$ (i), consider a $\QQ$-factorialization $Y\to X$
\cite{Kawamata-1988-crep} and run the MMP over $(Z,o)$. If (i) does
not hold, $\rho (Y/Z)>1$. On the last step of the MMP we get a
divisorial contraction $Y_{n-1}\to Y_n$ over $(Z,o)$. Let $E$ be the
corresponding exceptional divisor and let $\Gamma$ be its image on
$Z$. Then $\Gamma$ is an irreducible curve and $f^{-1}(\Gamma)$ has
two components.
\end{remark}

\begin{mparag}{}
\label{base-change}
We need the following easy construction which is to be used
throughout the paper. First, we claim that $(Z,o)$ is a quotient
singularity. Indeed, the general hyperplane section $H\subset X$ is
smooth and the restriction $f_H\colon H\to Z$ is a finite morphism.
Thus $(Z, o)$ is a log terminal singularity \cite[Prop.
5.20]{KM-1998}. Therefore, $(Z,o)$ is a quotient of a smooth germ
$(Z',o')$ by a finite group $G$ which acts freely outside of $o'$
\cite[Th. 9.6]{Kawamata-1988-crep}. Consider the base change
\begin{equation}
\label{eq-diag--base-change}
\begin{CD}
(X',C') @>{g}>> (X,C)
\\
@VV{f'}V @VV{f}V
\\
(Z',o')@>{h}>>(Z,o)
\end{CD}
\end{equation}
where $X'$ is the normalization of $X\times_Z Z'$ and
$C':=f'^{-1}(C)_{\red}$. The group $G$ naturally acts on $X'$ so
that $X=X'/G$. Since $X$ has only isolated singularities, $g$ is
\'etale in codimension $2$. Moreover, $K_{X'}=g^*K_X$ and
singularities of $X'$ are terminal. In particular, $f'\colon
(X',C')\to (Z',o')$ is a $\QQ$-conic bundle germ.
\end{mparag}

\begin{tparag}{Corollary ({\cite{Cutkosky-1988}}).}
\label{cor-Cut}
Let $f\colon X\to Z$ be a $\QQ$-conic bundle. If $X$ is Gorenstein
\textup(and terminal\textup), then $Z$ is smooth and there is a
vector bundle $\EEE$ of rank $3$ on $Z$ and an embeddings
$X\hookrightarrow \PP(\EEE)$ such that every scheme fiber $X_z$,
$z\in Z$ is a conic in $\PP(\EEE)_z$.
\end{tparag}

\begin{proof}[Sketch of the proof]
The question is local, so we assume that $f\colon (X,C)\to (Z,o)$ is
a $\QQ$-conic bundle germ. If $(Z,o)$ is smooth, the assertion can
be proved in the standard way: $f$ is flat because $X$ is
Cohen-Macaulay and we can put $\EEE=f_*\OOO_X(-K_X)$ (see, e.g.,
\cite{Cutkosky-1988}). Assume that $(Z,o)$ is singular. Consider the
base change \eqref{eq-diag--base-change}. Then $(Z',o')$ is smooth
and $G\neq\{1\}$. Since $X$ is Gorenstein terminal, the action of
$G$ on $X'$ and $C'$ is free. On the other hand, $X'$ is also
Gorenstein. By the above arguments $f'$ is a standard Gorenstein
conic bundle. In particular, the central fiber
$C':=f^{-1}(o')_{\red}$ is a conic. If $C'$ is reducible, then the
singular point $P'\in C'$ is $G$-invariant, a contradiction. Hence,
$C'\simeq \PP^1$. This contradicts the fact that the action of
$G$ on $C'$ is free.
\end{proof}

\begin{mparag}{Definition ({\cite[(0.4.16), (1.7)]{Mori-1988}}).}
Let $(X,P)$ be a terminal $3$-dimensional singularity of index $m$
and let $C\subset X$ be a smooth curve passing through $P$. We say
that $C$ is (locally) \textit{primitive} at $P$ if the natural map
\[
\varrho \colon \ZZ\simeq \pi_1(C\setminus \{P\}) \to
\pi_1(X\setminus \{P\})\simeq \ZZ/m\ZZ
\]
is surjective and \textit{imprimitive} 
at $P$
otherwise. The order $s$ of
$\Coker \varrho$ is called the \textit{splitting degree} and the
number $\bar m=m/s$ is called the \textit{subindex} of $P\in C$.
\end{mparag}
It is easy to see that the splitting degree coincides 
with the number of
irreducible components of the preimage $C^\sharp$ of $C$ under the
index-one cover $X^\sharp\to X$ near $P$. If $P$ is primitive, we
put $s=1$ and $\bar m=m$.

\begin{mparag}{}
>From now on we assume that $f\colon (X,C)\to (Z,o)$ is a $\QQ$-conic
bundle germ with $C\simeq\PP^1$. There are two cases (cf.
\cite[(1.12)]{Mori-1988}):
\end{mparag}

\begin{parag}{Case: $C'$ is irreducible.}
\label{case-prim-1}
\end{parag}

\begin{parag}{Case:  $C'=\cup_{i=1}^s C_i'$, where $s>1$ and $C'_i\simeq \PP^1$.}
\label{case-imprim-1}
In this case, $G$ acts on $\{C'_i\}$ transitively.
\end{parag}

\begin{claim}
\label{claim-imprimitive}
In the case \xref{case-imprim-1},
all the irreducible components
$C_i'$ pass through one point $P'$ and do not intersect each other
elsewhere.
\end{claim}
\begin{proof}
Since $G$ acts on $\{C'_i\}$ transitively and $p_a(C')=0$, each
component $C_i'$ meets the closure of $C'-C_i'$ at one point.
\end{proof}

\begin{tparag}{Proposition (cf. {\cite[1.11-1.13]{Mori-1988}}).}
\label{prop-first-prop-imp}
Notation as in \xref{base-change} and \xref{claim-imprimitive}.
\begin{enumerate}
\item
In the case \xref{case-imprim-1},
$C$ is imprimitive at $g(P')$.
Conversely, if $C$ is imprimitive at some point $P$, then $f$ is
such as in \xref{case-imprim-1}
and $P=g(P')$. Moreover, the splitting degree
$s$ coincides with the number of irreducible components of $C'$.
\item
$(X,C)$ has at most one imprimitive point.
\end{enumerate}
\end{tparag}

\begin{mtparag}{Proposition ({\cite[Lemma 1.10]{Prokhorov-1997_e}}).}
\label{prop-cyclic-quo}
$(Z,o)$ is a cyclic quotient singularity.
\end{mtparag}
\begin{proof}
It is sufficient to show that $G$ is a cyclic group. In the 
case where $X$ is locally
primitive, $G$ effectively acts on $C'\simeq\PP^1$ and on the
tangent space $T_{Z',o'}\simeq \CC^2$. This gives us two embeddings:
$G\subset PGL(2,\CC)$ and $G\subset GL(2,\CC)$. Assume that $G$ is
not cyclic. By the classification of finite subgroups in
$PGL(2,\CC)$ $G$ is either $\mathfrak A_5$, $\mathfrak S_4$,
$\mathfrak A_4$, or the dihedral group $\mathfrak D_n$ of order $2n$
(see, e.g., \cite{Springer-1977}). In all cases there are at least
two different elements of order two in $G$. But then at least one of
them is a reflection in $G\subset GL(2,\CC)$, a contradiction.

In the case where $f$ is not locally 
primitive, by Claim \ref{claim-imprimitive},
$G$ has a fixed point $P'\in X'$. Let $P=g(P')$ and let $U\ni P$ be
a small neighborhood. There is a surjection $\pi_1(U\setminus
\{P\})\twoheadrightarrow G$. On the other hand, $\pi_1(U\setminus
\{P\})$ is cyclic \cite[Lemma 5.1]{Kawamata-1988-crep}.
\end{proof}

Thus we may assume that $G=\muu_d$ and $Z\simeq \CC^2/\muu_d$, where
the action of $\muu_d$ on $\CC^2\simeq Z'$ is free outside of $0$.
We call this $d$ the \textit{topological index} of $f\colon (X,C)\to
(Z,o)$.

Let $\Clsc X$ be the subgroup of the divisor class group $\mt{Cl} X$
consisting of Weil divisor classes which are $\QQ$-Cartier.
\begin{corollary}
\label{cor-cyclic-quo}
$\pi_1(X\setminus \Sing X)\simeq \muu_d$ and $\Clsc X\simeq
\ZZ\oplus \ZZ_d$, where $\ZZ_d=\ZZ/d\ZZ$ and $d$ is the topological
index of $f$.
\end{corollary}

\begin{corollary}
\label{cor-cyclic-quo-new}
In the above notation, let $P_1$, \dots, $P_l$ be all the
non-Gorenstein points and let $m_1,\dots,m_l$ be their indices
\textup(the case $l=1$ is not excluded\textup). Assume that
$P_2,\dots, P_l$ are primitive. Let $s_1$ and $\bar m_1$ be the
splitting degree and the subindex of $P_1$.
\begin{enumerate}
\item
For each prime $p$ the number of the $m_i$'s divisible by $p$ is
$\le 2$.
\item
There is a $\QQ$-Cartier Weil divisor $D$ on $X$ such that $D\cdot
C=d/m_1\cdots m_l$. Moreover, $D$ generates $\Clsc
X/\mt{Torsion}=\Clsc X/{\equiv}$.
\item
$\prod\limits_{i=1}^l m_i=d\cdot \lcm (\bar m_1,m_2,\dots, m_l)$.
\end{enumerate}
\end{corollary}

\begin{proof}
Let $H$ be an ample generator of $\Pic X$ so that $H\cdot C=1$.
Clearly, the following sequence
\begin{equation}
\label{eq-ex-seq-Pic_Cl}
0 \longrightarrow \Pic X \longrightarrow \Clsc X
\stackrel{\varsigma}{\longrightarrow} \oplus_i\Clsc (X,P_i)
\longrightarrow 0
\end{equation}
is exact. Here $\Clsc (X,P_i)\simeq \ZZ_{m_i}$ by \cite[Lemma
5.1]{Kawamata-1988-crep}. Then (i) immediately follows by
\ref{cor-cyclic-quo}.

Let us prove (ii). We have $\Clsc X/\ZZ_d\simeq \ZZ$ by
\ref{cor-cyclic-quo} and the order of $(\Clsc X/\ZZ_d)/\Pic X$ is
$\frac 1d\prod {m_i}$ by \eqref{eq-ex-seq-Pic_Cl}. Let $D$ be an
ample Weil divisor generating $\Clsc X/\ZZ_d$. Since $H\cdot C=1$,
We have $\frac 1d \prod {m_i} D\cdot C=H\cdot C=1$. This proves
(ii).

Finally, by \cite[1.9, 1.7]{Mori-1988} the $\ZZ$-module $\Clsc X$ is
generated by $H$ and some ample Weil divisors $D_1,\dots, D_l$ with
relations $m_iD_i-n_iH\sim 0$, where $n_i=m_i D_i\cdot C$, $\gcd
(m_1,n_1)=s_1$, and $\gcd (m_i,n_i)=1$ for $i=2,\dots, l$. Now (iii)
can be proved by considering the $p$-prime component of $\Clsc X$
for each prime $p$.
\end{proof}

\begin{corollary}
\label{cor-pr-top-1}
In the locally primitive case, $\muu_d$ has exactly two fixed points
on $(X',C'\simeq \PP^1)$. Therefore, there are two points on $(X,C)$
whose indices are divisible by $d$. Conversely, if there are two
primitive points on $(X,C)$ whose indices divisible by $r$, then $r$
divides $d$.
\end{corollary}

\begin{corollary}
\label{cor-imp-top-1}
In the case $(X,C)$ is imprimitive at $P=g(P')$,
the
splitting degree $s$ $(>1)$ divides $d$, and let $r:=d/s$. 
Put $X^\flat:=X'/\muu_r$, $Z^\flat:=Z'/\muu_r$,
and $C^\flat:=C'/\muu_r$. We have the following decomposition:
\begin{equation}
\label{eq-cor-imp-top-diag}
\begin{CD}
(X',C') @>{g''}>> (X^\flat,C^\flat) @>{g^\flat}>> (X,C)
\\
@V{f'}VV @V{f^\flat}VV @V{f}VV
\\
(Z',o')@>{h'}>> (Z^\flat,o^\flat) @>{h^{\flat}}>> (Z,o)
\end{CD}
\end{equation}
and the following hold:
\begin{enumerate}
\item 
The group $\muu_{r}$ does not permute components of $C'$, so
$C^\flat$ has exactly $s$ irreducible components $C^\flat_i$ 
passing through one
point $P^\flat=g''(P')$. The group $\muu_s=\muu_d/\muu_r$ naturally
acts on $X^\flat$ so that $X=X^\flat/\muu_{r}$. 
\item 
If $d>s$, then
$\muu_{r}$ has two fixed points on each component $C_i'\subset C'$,
$P'$ and $Q_i'\neq P'$.
\item 
$(X^{\flat},C^{\flat}_i)$ is 
a locally primitive
extremal neighborhood 
\xref{rem-prel-extr-nbd},
$X^{\flat} \to X$ is \'etale outside $P^{\flat}$
and 
$C^{\flat}_i \to C$ is an isomorphism.
\end{enumerate}
\end{corollary}
The base change $g^\flat$ as in \eqref{eq-cor-imp-top-diag} is called
the \textit{splitting cover} \cite[1.12.1]{Mori-1988}.

\begin{proof}
Let $G\subset \muu_d$ be the stabilizer of some component
$C_i'\subset C'$. Then $G=\muu_{r}$ and $X^\flat:=X'/G$ satisfies
the desired properties.
\end{proof}

\begin{lemmam}
\label{lemma-KC}
Let $(X,C)$ be a $\QQ$-conic bundle germ with $C\simeq \PP^1$. Let
$d$ be the topological index of $(X,C)$ and let $m_1,\dots,m_r$ be
indices of all the non-Gorenstein points. Assume that $f$ is not
toroidal. Then
\[
-K_X\cdot C=d/m_1\cdots m_r.
\]
\end{lemmam}

\begin{proof}
Take $D$ as in (ii) of Corollary \ref{cor-cyclic-quo-new}. Then
$-K_X\equiv r D$ for some $r\in \ZZ_{>0}$. We claim that $r=1$.
Indeed, for the general fiber $L$ we have $2=-K_X\cdot L=rD\cdot L$.
Since $D\cdot L$ is an integer, $r=1$ or $2$. If $r=2$, then $D\cdot
L=1$, i.e., $D$ is $f$-ample with $\deg=1$ on the general fiber.
Apply construction \eqref{eq-diag--base-change}. Then $D':=f'^{*} D$
satisfies the same property: it is $f'$-ample with $\deg=1$ on the
general fiber. Since $X'\setminus f'^{-1}(o')\to Z'\setminus \{o'\}$
is a standard conic bundle (see \ref{cor-Cut}), this implies that
all the fibers over $Z'\setminus \{o'\}$ are smooth rational curves.
In particular, the morphism $f'$ is smooth outside of $C'$. We claim
that $f'$ is smooth everywhere. Denote $\FFF:=\OOO_{X'}(D')$. Then
locally near a singular point $P'\in X'$, $\FFF$ is a direct summand
of $\pi_*\FFF^\sharp$, where $\pi\colon (X^\sharp,P^\sharp) \to
(X',P')$ is the index-one cover and $\FFF^\sharp$ is the lifting of
$\FFF$. Since $\FFF^\sharp$ is Cohen-Macaulay and $Z'$ is smooth,
$\FFF$ is flat over $Z'$. Then by the base change theorem
(\cite[Lect. 7, (iii), p. 51]{Mumford-Lectures-on-curves})
$f'_*\FFF$ is locally free. Put $\hat{X}:=\PP(f'_*\FFF)$ with
natural projection $\hat{f}\colon \hat{X}\to Z'$. We have a
bimeromorphic map $\hat{X}\dashrightarrow X'$ over $Z'$ that indices
an isomorphism $(\hat{X}\setminus\hat{C})\simeq (X'\setminus C')$,
where $\hat{C}=\hat{f}^{-1}(o)$. Since $f'$, $\hat{f}$ are
projective and $\rho(X'/Z')=\rho(\hat{X}/Z')=1$, we have
$\hat{X}\simeq X'$. But then $X'$ is smooth and so is the morphism
$f'$. This proves our claim.

Thus we may assume $X'\simeq Z'\times \PP^1$. Recall that $X$ is the
quotient of $X'$ by $\muu_d$. By \cite[\S 2]{Prokhorov-1997_e} the
action of $\muu_d$ is as in \ref{ex-toric} and $f$ is toroidal, a
contradiction to our assumption. Therefore, $r=1$ and $-K_X\cdot
C=D\cdot C$. This proves our equality.
\end{proof}

\begin{corollary}
\label{cor-eqiv-imp}
If $X$ has a unique non-Gorenstein point which is imprimitive of
splitting degree $s$ and subindex $\bar m$, then $2\bar m\equiv 0
\mod s$.
\end{corollary}

\begin{proof}
Let $f'^{-1}(o')$ be the scheme fiber. Then $f'^{-1}(o')\equiv r C'$
for some $r\in \ZZ_{>0}$. Thus $2=-K_{X'}\cdot f'^{-1}(o')=
-rK_{X'}\cdot C'=-rsK_X\cdot C=rs/\bar m$. This proves our
statement.
\end{proof}
The following fact will be used freely.

\begin{mtparag}{Proposition (\cite[Th. 2.4]{Prokhorov-1997_e}).}
\label{prop-Gor-quot}
In notation of \eqref{eq-diag--base-change} assume that $X'$ is
Gorenstein \textup(we do not assume that $C$ is irreducible\textup).
Assume further that $d>1$.
Then $(X,C)$ is in one of the cases \xref{item=main--th-pr-toric},
\xref{item=main--th-pr-ex3}, \xref{item-main-th-impr-barm=1}.
\end{mtparag}

\begin{proof}[Sketch of the proof]
By \ref{cor-Cut} there is a $\muu_d$-equivariant embedding
$X'\hookrightarrow \PP^2\times Z'$ over $Z'$. Then one can choose a
suitable coordinate system in $\PP^2$ and $Z'\simeq \CC^2$.
\end{proof}

\section{Numerical invariants $i_P$, $w_P$ and $w_P^*$}
For convenience of the reader we recall some basic notation of
\cite{Mori-1988}.

\begin{mparag}{}
\label{not-iPwP}
Let $X$ be an analytic threefold with terminal singularities and let
$C\subset X$ be a reduced curve. Let $I_C\subset \OOO_X$ be the
ideal sheaf of $C$ and let $I_C^{(n)}$ be the symbolic $n$th power
of $I_C$, that is, the saturation of $I_C^n$ in $\OOO_X$. Put
$\gr_C^n\OOO:=I^{(n)}_C/I^{(n+1)}_C$. Further, let $F^n\omega_X$ be
the saturation of $I^n_C\omega_X$ in $\omega_X$ and let
$\gr_C^n\omega:=F^n\omega_X/F^{n+1}\omega_X$. Let $m$ be the index
of $K_X$. There are natural homomorphisms
\[
\begin{array}{lllll}
\alpha_1 &\colon& \bigwedge^2 \gr_C^1\OOO &\longrightarrow&
\HHom_{\OOO_C}(\Omega_C^1,\gr_C^0\omega),
\\[7pt]
\alpha_n &\colon& S^n \gr_C^1\OOO &\longrightarrow& \gr_C^n\OOO,
\quad n\ge 2,
\\[7pt]
\beta_0&\colon& (\gr_C^0\omega)^{\otimes m} &\longrightarrow&
(\omega_X^{\otimes m})^{**}\otimes \OOO_C,
\\[7pt]
\beta_n&\colon& \gr_C^0\omega\otimes S^n \gr_C^1\OOO
&\longrightarrow& \gr_C^n\omega, \quad n\ge 1,
\end{array}
\]
where $M^*$ for an $\OOO_X$-module $M$ denotes its dual,
$\HHom{\OOO_X} (M, \OOO_X)$. Denote
\[
i_P(n):=\len_P\Coker \alpha_n,\qquad w_P(0):=\len_P\Coker \beta_0/m,
\]
\[
w_P(n):=\len_P\Coker \beta_n, \quad w_P^*(n):=
{{n+1}\choose {2}}
i_P(1)-w_P(n), \quad n\ge 1.
\]
\end{mparag}
Assume that $C\simeq \PP^1$. Then we have by \cite[2.3.1]{Mori-1988}
\begin{equation}
\label{eq-grw-w}
-\deg \gr_C^0\omega=-K_X\cdot C+\sum_P w_P(0)
\end{equation}
\begin{equation}
\label{eq-grO-iP1}
2+\deg \gr_C^0\omega-\deg \gr_C^1\OOO=\sum_P i_P(1).
\end{equation}
\begin{equation}
\label{eq-grO-iPn}
\deg \gr_C^n\OOO=\frac12 n(n+1) \deg \gr_C^1\OOO +\sum_P i_P(n),
\qquad n\ge 2,
\end{equation}
and therefore the following corollaries to $\rk \gr_C^i\OOO=i+1$ and
$R^1f_*\OOO_X=0$:
\begin{equation}
\label{eq-grO-iP1-1}
\sum_{i=1}^n (\deg \gr_C^i\OOO + i+1)\ge 0,\qquad n\ge 1,
\end{equation}
\begin{equation}
\label{eq-grO-iP1-2}
4\ge -\deg \gr_C^0\omega+\sum_P i_P(1)= -K_X\cdot C+\sum_P
w_P(0)+\sum_P i_P(1).
\end{equation}

\begin{remark}
In the case of extremal neighborhoods by the Grauert-Riemenshneider
vanishing one has $\gr_C^0\omega=\OOO_{C}(-1)$ (see
\cite[2.3]{Mori-1988}). This is no longer true for $\QQ$-conic
bundles: in Example \xref{ex-toric} easy computations show 
$\deg \gr_C^0\omega= -2$ (see \eqref {eq-grw-w}). Similarly, in
\ref{item-main-th-impr-barm=1} we also have $\deg \gr_C^0\omega=
-2$.
We will show below that these two examples are the only
exceptions (see Corollaries \ref{cor-prop-grw=2-2-points-prim} and
\ref{cor-impr-grw-21cases}).
\end{remark}

\begin{lemma}
\label{lemma-grw-ex-esq}
If $\gr_C^0\omega=\OOO(-1)$, then
\begin{equation}
\label{eq-grw-1}
\deg \gr_C^n\omega =\frac12 (n+1)(n-2)-\sum_P w_P^*(n),\quad n\ge1.
\end{equation}
If furthermore $H^1(\omega_X/F^{n+1}\omega_X)=0$, then
\begin{equation}
\label{eq-grw-1-2}
\sum_{i=1}^n \left(\deg \gr_C^i\omega +i+1\right)\ge 0, \qquad n\ge
1.
\end{equation}
\end{lemma}
\begin{proof}
Follows by the exact sequences
\[
0\longrightarrow \gr_C^i\omega \longrightarrow
\omega_X/F^{i+1}\omega_X \longrightarrow\omega_X/F^i\omega_X
\longrightarrow 0
\]
(see \cite[2.3]{Mori-1988}).
\end{proof}

\begin{tparag}{Lemma ({\cite[2.15]{Mori-1988}}).}
If $(X,P)$ is singular, then $i_P(1)\ge 1$.
\end{tparag}

\begin{proof}
The proof of \cite[2.15]{Mori-1988} applies because it uses only
local computations near $P$ that are not based on $R^1f_*\omega_X$.
\end{proof}

\begin{corollary}
A $\QQ$-conic bundle germ $(X,C\simeq\PP^1)$ has at most three
singular points.
\end{corollary}

\section{Sheaves $\gr_C^n\omega$}
\begin{lemmam}
\label{lemma-omega-main}
Let $f\colon X \to Z$ be a $\QQ$-conic bundle. Assume that the base
surface $Z$ is smooth. Then there is a canonical isomorphism
$R^1f_*\omega_X\simeq \omega_Z$.
\end{lemmam}

\begin{proof}
Let $g\colon W\to X$ be a resolution. By {\cite[Prop.
7.6]{Kollar-1986-I}} we have $R^1(f \comp g)_* \omega_W=\omega_Z$.
Since $X$ has only terminal singularities, $g_*\omega_W=\omega_X$
and by the Grauert-Riemenshneider vanishing, $R^ig_* \omega_W=0$ for
$i>0$. Then the Leray spectral sequence gives us $R^1f_*
\omega_X=R^1(f \comp g)_* \omega_W= \omega_Z$.
\end{proof}

For convenience of the reader we recall basic definitions
\cite[8.8]{Mori-1988}.
\begin{mparag}{}
\label{pre-def-totimes}
Let $(X,P)$ be three-dimensional terminal singularity of index $m$
and let $\pi\colon (X^\sharp,P^\sharp)\to (X,P)$ be the index-one
cover. Let $\LLL$ be a coherent sheaf on $X$ without submodules of
finite length $>0$. An \textit{$\ell$-structure} of $\LLL$ at $P$ is
a coherent sheaf $\LLL^\sharp$ on $X^\sharp$ without submodules of
finite length $>0$ with $\muu_m$-action endowed with an isomorphism
$(\LLL^\sharp)^{\muu_m}\simeq \LLL$. An \textit{$\ell$-basis of
$\LLL$ at $P$} is a collection of $\muu_m$-semi-invariants
$s_1^\sharp,\dots,s_r^\sharp\in \LLL^\sharp$ generating
$\LLL^\sharp$ as an $\OOO_{X^\sharp}$-module at $P^\sharp$. Let $Y$
be a closed subscheme of $X$ such that $P\notin \mt{Ass} \OOO_Y$ and
let $Y^\sharp\subset X^\sharp$ be the canonical lifting. Note that
$\LLL$ is an $\OOO_Y$-module if and only if $\LLL^\sharp$ is an
$\OOO_{Y^\sharp}$-module. We say that $\LLL$ is \textit{$\ell$-free}
$\OOO_Y$-module at $P$ if $\LLL^\sharp$ is a free
$\OOO_{Y^\sharp}$-module at $P^\sharp$. If $\LLL$ is $\ell$-free
$\OOO_Y$-module at $P$, then an $\ell$-basis of $\LLL$ at $P$ is
said to be \textit{$\ell$-free} if it is a free
$\OOO_{Y^\sharp}$-basis.

Let $\LLL$ and $\MMM$ be $\OOO_Y$-modules at $P$ with
$\ell$-structures $\LLL\subset \LLL^\sharp$ and $\MMM\subset
\MMM^\sharp$. Define the following operations:
\begin{itemize}
\item
$\LLL\toplus \MMM\subset (\LLL\oplus \MMM)^\sharp$ is an
$\OOO_Y$-module at $P$ with $\ell$-structure
\[
(\LLL\toplus \MMM)^\sharp =\LLL^\sharp\oplus\MMM^\sharp.
\]
\item
$\LLL\totimes \MMM\subset (\LLL\otimes \MMM)^\sharp$ is an
$\OOO_Y$-module at $P$ with $\ell$-structure
\[
(\LLL\totimes \MMM)^\sharp
=(\LLL^\sharp\otimes_{\OOO_{X^\sharp}}\MMM^\sharp) /
\mt{Sat}_{\LLL^\sharp\otimes\MMM^\sharp}(0),
\]
where $\mt{Sat}_{\mathcal F_1}\mathcal F_2$ is the saturation of
$\mathcal F_2$ in $\mathcal F_1$.
\end{itemize}
These operations satisfy standard properties (see
\cite[8.8.4]{Mori-1988}). If $X$ is an analytic threefold with
terminal singularities and $Y$ is a closed subscheme of $X$, then
the above local definitions of $\toplus$ and $\totimes$ patch with
corresponding operations on $X\setminus \Sing X$. Therefore, they
give well-defined operations of global $\OOO_Y$-modules.
\end{mparag}

\begin{mparag}{}
\label{pre-not-cb}
Let $f\colon (X,C)\to (Z,o)$ be a $\QQ$-conic bundle germ \textup(we
do not assume that $C$ is irreducible\textup).
\end{mparag}

\begin{theoremm}{}
\label{theo-sect-1-gen-fiber-m}
Assume that $(Z,o)$ is smooth. Let $J\subset \OOO_X$ be an ideal
such that $\Supp \OOO_X/J \subset C$ and $\OOO_X/J$ has no embedded
components. Assume that $H^1(\omega_X\totimes \OOO_X/J)\neq 0$. Then
$\Spec_X \OOO_X/J\supset f^{-1}(o)$, where $f^{-1}(o)$ is the scheme
fiber \textup(in other words, $J\subset \mathfrak{m}_{Z,o}\OOO_X$,
where $\mathfrak{m}_{Z,o}\subset \OOO_Z$ is the maximal ideal of
$o$\textup).
\end{theoremm}

\begin{proof}
First we assume that $\Spec_X \OOO_X/J\subsetneq f^{-1}(o)$. Denote
$\Gamma:= f^{-1}(o)$ and $V:=\Spec \OOO_X/J$. Then
$\omega_{\Gamma}\simeq \omega_X\otimes \OOO_{\Gamma}$, and so
$\omega_X\totimes \OOO_V \simeq \omega_\Gamma \totimes \OOO_V$ in
this case.

Let $\mathcal{I}_V$ be the ideal sheaf of $V$ in $\Gamma$
where we note $\mathcal{I}_V \not=0$ by $V\subsetneq f^{-1}(o)$, and
let
$\mathcal{I}_D$ be an associated prime of $\mathcal{I}_V$
(i.e. $\mathcal{I}_D \in \operatorname{Ass} (\mathcal{I}_V)$), and let
$D\subset C$ be the corresponding irreducible component. By the Serre
duality, we have
\[
\omega_{D}=\HHom_{\OOO_\Gamma}(\OOO_{D},
\omega_\Gamma)=\HHom_{\OOO_\Gamma}(\OOO_{D}, \OOO_\Gamma)\totimes
\omega_\Gamma.
\]
Hence
$\HHom_{\OOO_\Gamma}(\OOO_D,\OOO_\Gamma)$ is a torsion-free 
$\OOO_D$-module of
rank $1$.
We also see
$\HHom_{\OOO_\Gamma}(\OOO_D,\mathcal I_V)\neq 0$
by $\mathcal{I}_D \in \operatorname{Ass} (\mathcal{I}_V)$.  
Thus the cokernel of the inclusion
\[
0\neq \HHom_{\OOO_\Gamma}(\OOO_D,\mathcal I_V)\hookrightarrow
\HHom_{\OOO_\Gamma}(\OOO_D,\OOO_\Gamma).
\]
is of finite length and is a submodule of $\OOO_V=\OOO_\Gamma/\mathcal
I_V$. Since $\OOO_V=\OOO_\Gamma/\mathcal I_V$ has
no embedded primes, we have
$\HHom_{\OOO_\Gamma}(\OOO_D,\mathcal I_V)=
\HHom_{\OOO_\Gamma}(\OOO_D,\OOO_\Gamma)$ and
\[
\mathcal I_{V}=\HHom_{\OOO_\Gamma}(\OOO_\Gamma,\mathcal I_V)
\supset
\HHom_{\OOO_\Gamma}(\OOO_D,\OOO_\Gamma).
\]

Considering the trace map one can see that $\CC\simeq
H^1(\omega_{D})\to H^1(\omega_{\Gamma})$ is an injection (and
moreover $H^1(\omega_{D})\simeq H^1(\omega_{\Gamma}) \simeq \CC$).
Since $\omega_\Gamma$ is $\ell$-invertible, the composition map
\[
\upsilon\colon H^1(\omega_D) \stackrel{{\sim}\phantom{{\sim}}}{\to}
H^1(\omega_\Gamma) \to H^1(\omega_\Gamma \totimes \bigl(\OOO_\Gamma/
\HHom_{\OOO_\Gamma}(\OOO_{D}, \OOO_\Gamma)\bigr))
\]
is zero.

On the other hand, $\omega_\Gamma \to \omega_X \totimes \OOO_V$
has the following decomposition
\[
\omega_\Gamma \to \omega_\Gamma \totimes \bigl(\OOO_\Gamma/
\HHom_{\OOO_\Gamma}(\OOO_{D}, \OOO_\Gamma) \bigr) \to \omega_\Gamma
\totimes \OOO_\Gamma/\mathcal I_{V}\simeq \omega_X \totimes \OOO_V,
\]
and the induced surjective map
\[
H^1(\omega_\Gamma)\longrightarrow H^1(\omega_\Gamma\totimes
\OOO_V)\neq 0
\]
factors through $\upsilon$ which is zero, a contradiction. 

This
proves that $\Spec_X \OOO_X/J=f^{-1}(o)$ if $\Spec_X \OOO_X/J\subset
f^{-1}(o)$.

Now we treat the general case. 
By Nakayama's lemma 
$H^1(\omega_X\totimes \OOO_X/J)\otimes_{\OOO_Z} \OOO_Z/\mathfrak{m}_{Z,o} \neq 0$.
Since $H^1$ is right exact for $\OOO_X$-sheaves, we see that
$$H^1((\omega_X\totimes \OOO_X/J)\otimes \OOO_X/\mathfrak{m}_{Z,o}\OOO_X) \simeq H^1(\omega_X\totimes \OOO_X/J)\otimes_{\OOO_Z} 
\OOO_Z/\mathfrak{m}_{Z,o}\not=0.$$

Let us consider the homomorphism
$$
(\omega_X\totimes \OOO_X/J)\otimes \OOO_X/\mathfrak{m}_{Z,o}\OOO_X
\to
\omega_X\totimes \OOO_X/J^s,
$$
where $J^s$ is the 
saturation of $J+\mathfrak{m}_{Z,o}\OOO_X$ in $\OOO_X$.
It is surjective and its 
kernel 
is supported at 
a finite number of 
points.
Thus
$$
H^1\bigl(\omega_X\totimes \OOO_X/J^s\bigr)\simeq
H^1\bigl((\omega_X\totimes \OOO_X/J)\otimes \OOO_X/\mathfrak{m}_{Z,o}\OOO_X\bigr)
\neq 0
$$
and $J^s\supset
\mathfrak{m}_{Z,o}\OOO_X$. By the special case treated above we have
$J+\mathfrak{m}_{Z,o}\OOO_X \subset J^s= \mathfrak{m}_{Z,o}\OOO_X$,
i.e., $J\subset \mathfrak{m}_{Z,o}\OOO_X$.

\end{proof}

\begin{corollary}
\label{cor-gr-omega-=0}
Assume that $(Z,o)$ is smooth. If $H^1(\gr_C^0\omega)\neq 0$, then
$C=f^{-1}(o)$.
\end{corollary}
\begin{proof}
Apply Theorem \xref{theo-sect-1-gen-fiber-m} with $J=I_C$.
\end{proof}

\begin{tparag}{Lemma ({\cite[Prop. 4.2]{Kollar-1999-R}}).}
\label{lemma-int-non-Gor}
If $X$ is not Gorenstein, then $X$ has index $>1$ at all singular
points of $C$.
\end{tparag}
\begin{proof}
If $C$ has at least three irreducible components, the assertion
follows by Remark \ref{rem-prel-extr-nbd} and \cite[Cor.
1.15]{Mori-1988}. Thus we assume that $C=C_1\cup C_2$ and $X$ is
Gorenstein at $P\in C_1\cap C_2$. First we consider the case when
$(Z,o)$ is smooth. By our assumption $\gr_C^0\omega=\omega_X\otimes
\OOO_C$ is invertible at $P$. Consider the injection $\varphi\colon
\gr_C^0\omega\hookrightarrow \gr_{C_{1}}^0\omega \oplus
\gr_{C_{2}}^0\omega$. Recall that $(X,C_i)$ is an extremal
neighborhood by Remark \ref{rem-prel-extr-nbd}. Then by \cite[Prop.
1.14]{Mori-1988} $\gr_{C_{i}}^0\omega=\OOO_{C_i}(-1)$, so
$H^0(\Coker \varphi)= H^1(\gr_C^0\omega)$. On the other hand,
$\Coker \varphi$ is a sheaf of finite length supported at $P$. Since
$\gr_C^0\omega$ is invertible, $\Coker \varphi$ is non-trivial. So,
$H^1(\gr_C^0\omega)\neq 0$ and by Corollary \ref{cor-gr-omega-=0}
$C_1\cup C_2=f^{-1}(o)$. Thus $X$ is smooth outside of $\Sing C$.
Since $P$ is the only singular point of $C$ by Corollary
\ref{cor-C-pa=0}, we are done.

Now we assume that $(Z,o)$ is singular. Consider the base change
\eqref{eq-diag--base-change}. Since $X$ is Gorenstein terminal at
$P$, so is $X'$ at all the points $P_i'\in g^{-1}(P)$. Moreover, $g$
is \'etale over $P$. Hence, the central curve $C'$ is singular at
$P_i'$. By the above, $X'$ is Gorenstein and by Corollary
\ref{cor-Cut} $f'\colon X'\to Z'$ is a standard Gorenstein conic
bundle. In particular, $C'$ is a plane conic. Since the set
$g^{-1}(P)$ is contained in the singular locus of $C'$, it consists
of one point, a contradiction.
\end{proof}

\begin{tparag}{Corollary (cf. {\cite[Prop. 1.14]{Mori-1988}}).}
\label{cor-prop-grw=2-2-points-prim}
Assume that $C$ is irreducible. If $\gr_C^0\omega \not\simeq
\OOO_C(-1)$, then in notation of \eqref{eq-diag--base-change} we
have $f'^{-1}(o')=C'$. If furthermore $(X,C)$ is locally primitive,
then it is toroidal.
\end{tparag}

\begin{proof}
Let $m$ be the index of $X$. Since there is an injection
$(\gr_C^0\omega)^{\otimes m}\hookrightarrow \OOO_C(mK_X)$, $\deg
\gr_C^0\omega<0$. Since $\gr_C^0\omega \not\simeq \OOO_C(-1)$,
$H^1(\gr_C^0\omega)\neq 0$. In notation of
\eqref{eq-diag--base-change} we have $H^1(\gr_{C'}^0\omega)\neq 0$
(because $H^1(\gr_{C}^0\omega)=H^1(\gr_{C'}^0\omega)^{\muu_d}$).

By
Corollary \ref{cor-gr-omega-=0} $C'=f'^{-1}(o')$. If $f$ is locally
primitive, $C'$ is irreducible (see \ref{case-prim-1}). So $C'\simeq
\PP^1$ and $X'$ is smooth. Up to analytic isomorphism we may assume
that $X'\simeq Z'\times \PP^1$. Then in some coordinate system the
action of $\muu_d$ on $X'$ is as in \ref{ex-toric} (see \cite[\S
2]{Prokhorov-1997_e}), so $f$ is toroidal.
\end{proof}

\begin{remark}
\label{remark-sect-1-gen-fiber}
In notation of Theorem \xref{theo-sect-1-gen-fiber-m} assume that
the map $H^0(I_C)\to H^0(I_C/J)$ is zero. Then $\Spec_X
\OOO_X/J\subset f^{-1}(o)$. Therefore, the nonvanishing
$H^1(\omega_X\totimes \OOO_X/J)\neq 0$ implies $\Spec_X \OOO_X/J=
f^{-1}(o)$.
\end{remark}

\begin{corollary}
\label{cor-lemma-sect-1-gen-omega-OOO-ne0}
Notation as in \xref{pre-not-cb}. Assume that $(Z,o)$ is smooth. If
the map $H^0(I_C)\to H^0(\gr_C^1\OOO)$ is zero, then
$H^1(\gr_C^1\omega)= 0$.
\end{corollary}

\begin{proof}
Assume that $H^1(\gr_C^1\omega)\neq 0$. In notation of Theorem
\xref{theo-sect-1-gen-fiber-m}, put $J=I_{C}^{(2)}$ and $V:=\Spec_X
\OOO_X/I_{C}^{(2)}$. From the exact sequence
\[
\begin{CD}
0 @>>> \gr_{C}^1\OOO \totimes \omega_X
 @>>> \OOO_X/ I_{C}^{(2)} \totimes \omega_X
@>>> \OOO_C \totimes \omega_X @>>> 0
\\
@. @| @| @|
\\
@. \gr_C^1\omega @.\OOO_V\totimes \omega_X @. \gr_C^0\omega @.
\end{CD}
\]
and $\deg \gr_{C_i}^0\omega<0$ for each $i$ we get
$H^1(\OOO_V\totimes \omega_X)\neq 0$. Then by Theorem
\xref{theo-sect-1-gen-fiber-m} and Remark
\xref{remark-sect-1-gen-fiber} $V=f^{-1}(o)$. Let $P\in C$ be a
general point. Then in a suitable coordinate system $(x,y,z)$ near
$P$ we may assume that $C$ is the $z$-axis. So, $I_C=(x,y)$ and
$I_C^{(2)}=(x^2,xy,y^2)$. But then $V$ is not a local complete
intersection near $P$, a contradiction.
\end{proof}

\begin{corollary}
\label{cor-wPstar}
Assume that $(Z,o)$ is smooth and $C$ is irreducible. If $\sum_P
i_P(1)\ge 3$, then $\sum w_P^*(1)\le 1$.
\end{corollary}
\begin{proof}
By \eqref{eq-grO-iP1-2} $\gr_C^0\omega=\OOO(-1)$. Further, by
\eqref{eq-grO-iP1} $\deg \gr_{C}^1\OOO\le -2$. Hence,
$H^0(\gr_{C}^1\OOO)=0$ (cf. \cite[Remark 2.3.4]{Mori-1988}). Now the
desired inequality follows by Corollary
\ref{cor-lemma-sect-1-gen-omega-OOO-ne0} and \eqref{eq-grw-1}.
\end{proof}

\section{Preliminary classification of singular points}
\begin{mparag}{Notation.}
\label{not-nazalo-loc}
Let $f\colon (X,C\simeq\PP^1)\to (Z,o)$ be a $\QQ$-conic bundle
germ. Let $P\in C$ be a point of index $m\ge 1$. Let $s$ and $\bar
m$ be the splitting degree and subindex, respectively. Thus $m=s\bar
m$. Consider the canonical $\muu_m$-cover $\pi\colon
(X^\sharp,P^\sharp)\to (X,P)$ and let $C^\sharp:=\pi^{-1}(C)$. Take
normalized $\ell$-coordinates $(x_1,\dots,x_4)$ and $t$ and let
$\phi$ be an $\ell$-equation of $X\supset C\ni P$ (see
\cite[2.6]{Mori-1988}). Put $a_i=\ord x_i$.

Note that $a_i<\infty$ and $\wt x_i\equiv a_i\mod \bar m$. If $m=1$,
then $X=X^\sharp$. In this case, $P$ is said to be \textit{of type
\type{III}}.
\end{mparag}

\begin{mparag}{Primitive point.}
Consider the case when $P$ is primitive and $m>1$. Then $s=1$ and
$\bar m=m$. In this case, the classification coincides with that in
\cite{Mori-1988} as shown next:
\end{mparag}

\begin{tparag}{Proposition (cf. {\cite[Prop. 4.2]{Mori-1988}}).}
\label{prop-prim-types-def}
Let $P$ and $m$ be as above. Modulo permutations of $x_i$'s, the
semigroup $\ord C^\sharp$ is generated by $a_1$ and $a_2$. Moreover,
exactly one of the following holds:
\begin{enumerate}
\item [\type{IA}]
$a_1+a_3\equiv 0\mod m$, $a_4=m$, $m\in \ZZ_{>0} a_1+\ZZ_{>0} a_2$,
where we may still permute $x_1$ and $x_3$ if $a_2=1$,
\item [\type{IB}]
$a_1+a_3\equiv 0\mod m$,\ $a_2=m$,\ $2\le a_1$,
\item [\type{IC}]
$a_1+a_2=a_3=m$,\ $a_4\not\equiv a_1,\, a_2 \mod m$,\ $2\le a_1<a_2$,
\item [\type{IIA}]
$m=4$, $P$ is of type $cAx/4$, and $\ord x=(1,1,3,2)$,
\item [\type{IIB}]
$m=4$, $P$ is of type $cAx/4$, and $\ord x=(3,2,5,5)$.
\end{enumerate}
\end{tparag}
\begin{proof}
If $X$ has an imprimitive point ($\neq Q$),
then 
$P$ is as classified in \cite[Prop.4.2]{Mori-1988}
by \ref {cor-imp-top-1}, (iii). 
So we can assume that $X$ has
no imprimitive points.
If $(X,C)$ is toroidal, then at both singular points $\ord
x=(1,a,m-1,m)$. So, these points are of type \type{IA}. Taking
Corollary \ref{cor-prop-grw=2-2-points-prim} into account, we may assume
that $\gr_C^0\omega=\OOO_{C}(-1)$. By \eqref{eq-grw-w} and
\eqref{eq-grO-iP1-2} we have $w_P(0)<1$ and $i_P(1)\le 3$. We claim
that $C^\sharp$ is planar, i.e., $\ord C^\sharp$ is generated by two
elements. Indeed, in the contrary case by \cite[Lemma
3.4]{Mori-1988} we have $i_P(1)=3$. Hence $P$ is the only singular
point (see \eqref{eq-grO-iP1-2}) and $(Z,o)$ is smooth (Corollary
\ref{cor-pr-top-1}). By Corollary \ref{cor-wPstar} $w_P^*(1)\le 1$.
In this case, arguments of \cite[3.5]{Mori-1988} work. This shows
that $C^\sharp$ is planar. Now we can apply \cite[Proof of
4.2]{Mori-1988} and obtain the above classification.
\end{proof}

\begin{mparag}{Imprimitive point.}
Now assume that $P$ is imprimitive. Then in diagram
\eqref{eq-diag--base-change} the central fiber $C'$ has exactly $s$
$(>1)$ irreducible components. Note that the classification is
different from that in \cite{Mori-1988} only in the case $C'=
f'^{-1}(o')$.
\end{mparag}

\begin{tparag}{Proposition (cf. {\cite[Prop. 4.2]{Mori-1988}}).}
\label{prop-imp-types}
Let $P$, $C^\sharp$, and $s$ be as above. Modulo permutations of
$x_i$'s and changes of $\ell$-characters, the semigroup $\ow
C^\sharp$ is generated by $\ow x_1$ and $\ow x_2$ except for the
case \type{IE^{\vee}} below. Moreover, exactly one of the following
holds:
\begin{enumerate}
\item[\type{IA^{\vee}}]
$\bar m>1$, $\wt x_1+\wt x_3\equiv 0\mod m$, $\ow x_4=(\bar m,0)$,
$\ow C^\sharp$ is generated by $\ow x_1$ and $\ow x_2$, and
$w_P(0)\ge 1/2$.

\item[\type{IC^{\vee}}]
$s=2$, $\bar m$ is an even integer $\ge 4$, and
\[
\begin{array}{cccccc}
&x_1&x_2&x_3&x_4&
\\
\wt&1&-1&0&\bar m+1&\mod m
\\
\ord&1&\bar m-1&\bar m&\bar m+1&
\end{array}
\]

\item[\type{II^{\vee}}]
$\bar m=s=2$, $P$ is of type $cAx/4$, and
\[
\begin{array}{cccccc}
&x_1&x_2&x_3&x_4&
\\
\wt&1&3&3&2&\mod 4
\\
\ord&1&1&1&2&
\end{array}
\]

\item[\type{ID^{\vee}}]
$\bar m=1$, $s=2$, $P$ is of type $cA/2$ or $cAx/2$, and
\[
\begin{array}{cccccc}
&x_1&x_2&x_3&x_4&
\\
\wt&1&1&1&0&\mod 2
\\
\ord&1&1&1&1&
\end{array}
\]

\item[\type{IE^{\vee}}]
$\bar m=2$, $s=4$, $P$ is of type $cA/8$, and
\[
\begin{array}{cccccc}
&x_1&x_2&x_3&x_4&
\\
\wt&5&1&3&0&\mod 8
\\
\ord&1&1&1&2&
\end{array}
\]
\end{enumerate}
Moreover, we are in the case \type{ID^{\vee}} or \type{IE^{\vee}} if
only if $C'=f'^{-1}(o')$. In this case, $P$ is the only
non-Gorenstein point.
\end{tparag}

\begin{remark}
It is easy to show that in cases \type{IC^{\vee}} and
\type{IE^{\vee}} the point $(X,P)$ is a cyclic quotient singularity
(cf. \cite[Lemma 4.4]{Mori-1988}).
\end{remark}

\begin{proof}
First assume that $C'\neq f'^{-1}(o')$. By Corollary
\ref{cor-gr-omega-=0} we have $H^1(\gr_{C'}^0\omega)=0$. Therefore,
\[
H^1(\gr_{C}^0\omega)=H^1(\gr_{C'}^0\omega)^{\muu_d}=0.
\]
This implies $\gr_{C}^0\omega=\OOO_C(-1)$ and $w_P(0)<1$. In
particular, $\bar m>1$ (see \cite[Cor. 2.10]{Mori-1988}). Let
$g^\flat \colon (X^\flat,C^\flat)\to (X,C)$ be the splitting cover.
Consider the exact sequence
\begin{equation}
\label{eq-ex-seq-imp-gr}
0\longrightarrow \gr_{C^\flat}^0 \omega
\stackrel{\varphi}{\longrightarrow} \bigoplus_{i=1}^s
\gr_{C^\flat_i}^0\omega \longrightarrow \Coker \varphi
\longrightarrow 0.
\end{equation}
Note that $\gr_{C^\flat_i}^0\omega=\OOO(-1)$ (see
\cite[2.3.2]{Mori-1988}). Hence,
\[
H^0(\Coker \varphi)= H^1(\gr_{C^\flat}^0\omega)=
H^1(\gr_{C'}^0\omega)^{\muu_{d/s}}=0.
\]
Since the support of $\Coker \varphi$ is zero-dimensional, $\varphi$
is an isomorphism. Therefore, the classification
\cite[4.2]{Mori-1988} holds for $(X,P)$ in this case (see
\cite[3.6-3.8]{Mori-1988}).

Now we consider the case where $C'=f'^{-1}(o')$. If $\bar m=1$, then
by Lemma \ref{lemma-int-non-Gor} 
the splitting cover $X^{\flat}$
is Gorenstein and  
$-K_{X^{\flat}}\cdot C_i^{\flat}$
is an
integer for any component 
$C_i^{\flat}\subset C^{\flat}$. 
Hence, 
$2=-K_{X^{\flat}}\cdot C^{\flat}=-sK_{X^{\flat}}\cdot C_i^{\flat}$. 
This implies $s=2$. We get the case
\type{ID^\vee}. Furthermore by Proposition \ref{prop-Gor-quot}
we are in the case \ref{item-main-th-impr-barm=1}
and 
hence $P$ is the only
non-Gorenstein point.
>From now on we assume that $\bar m>1$.

\begin{parag}{}
\label{parag-new-imp-2points}
We claim that $P$ is the only non-Gorenstein point. 
Indeed, assume first that there are at least two non-Gorenstein
points 
other than $P$ on 
$C$. 
Then on the splitting cover $X^{\flat}$ any 
irreducible component 
$C^{\flat}_i$ 
of $C^{\flat}$ contains at least three
non-Gorenstein points by 
$\bar{m} > 1$ \ref{cor-imp-top-1}, (iii).
Since the extremal neighborhood 
$(X^{\flat}, C^{\flat}_i)$
\ref{cor-imp-top-1}, (iii) 
can have at most two non-Gorenstein
points 
\cite[Thm. 6.2]{Mori-1988}, this is impossible.
So we
assume that $(X,C)$ contains exactly two non-Gorenstein points,
$P$ and $Q$. Let $n$ be the index of $Q$.
Clearly,
$-K_{X'}\cdot C'=2$. On the other hand, by Lemma \ref{lemma-KC}
$-K_{X'}\cdot C'=-sK_X\cdot C=sd/m n$. Let $r=\gcd (\bar m,n)$. Then
$d=rs$ and $s=2 n \bar m /r$. Since $\bar m>1$, we have $s=2s_1$,
where $s_1=n \bar m /r>1$. Consider the quotient $X''$ of
$X^\flat$ by $\muu_{s_1}\subset \muu_s$. We get extremal
neighborhoods $(X'',C''_i)$ with two non-Gorenstein points:
imprimitive of index $\bar ms_1$ and primitive of index $n$. By
\cite[Th. 6.7, 9.4]{Mori-1988} this is impossible. Thus the claim is
proved. In particular, $X^\flat=X'$.
\end{parag}

As above we have $-K_{X'}\cdot C'=2=s/\bar m$. Hence, $s=2\bar m$.
In particular, $m=2\bar m^2\neq 4$ and $P$ is not of type $cAx/4$.
Up to permutation of $x_i$'s we may assume that $\wt x_4\equiv \wt
x_1 x_3\equiv \wt \phi \equiv 0\mod m$. Since $-K_X\cdot C=1/\bar m$
and $P$ is the only non-Gorenstein point, $\ord (x_1\cdots x_4/\phi)
\equiv -\bar m K_X\cdot C \equiv 1 \mod \bar m$ (see \cite[Corollary
2.10]{Mori-1988}). So, $a_2\equiv 1\mod \bar m$.

Consider the map $\varphi\colon \gr_{C'}^0 \omega \to \oplus
\gr_{C'_i}^0\omega$ (see \eqref{eq-ex-seq-imp-gr}) and the induced map
\begin{equation}
\label{eq-imp-Phi-gr}
\Phi\colon \gr_{C'}^0 \omega= (\OOO_{C^\sharp}\bar
\omega)^{\muu_{\bar m}} \to \bigoplus_{i=1}^s \gr_{C'_i}^0 \omega
\otimes \CC(P'),
\end{equation}
where $\bar \omega$ is a semi-invariant generator of
$\omega_{X^\sharp}$ at $P^\sharp$. For example we can take
\[
\bar \omega= \frac{dx_2 \wedge dx_3 \wedge d x_4} {\partial \phi/
\partial x_1}.
\]
Since $\bar m w_{P'_{(i)}}(0)=\bar m-1$, we see that $\Phi(\nu \bar
\omega)=0$ for $\muu_{\bar m}$-semi-invariants $\nu$ with $\ord \nu
\ge \bar m$.

\begin{lemma}
If $C'=f'^{-1}(o')$, then
\[
H^1(\gr_{C'}^0\omega)=H^0(\Coker \varphi)=
H^0(\Coker \Phi)= \CC.
\]
\end{lemma}
\begin{proof}
Since $H^j(\gr_{C'_i}^0\omega)=0$, we have
\[ 
\begin{array}{l}
H^0(\Coker \varphi)\simeq H^1(\gr_{C'}^0\omega)\simeq
H^1(\omega_{X'}\otimes \OOO_{f'^{-1}(o')})
\\
\qquad
\simeq (R^1f'_*\omega_{X'})\otimes \CC(o') \simeq \omega_{Z'}\otimes
\CC(o')\simeq \CC.
\end{array}
\]
(We used the base change theorem and Lemma \ref{lemma-omega-main}.)
\end{proof}

There are two cases.
\begin{pusto}
\textbf{Case $a_2\ge \bar m$.} Clearly, $a_1+a_3\ge \bar m$. Thus we
may assume that $a_1\le a_3 \ge \bar m/2$. In this case, $\Phi$
factors through
\[
\Bigl(\OOO_{C^\sharp, P^\sharp}/ (x_1^{\bar m}, x_1x_3, x_3^2, x_2,
x_4)\cdot \bar \omega\Bigr)^{\muu_{\bar m}} \simeq \CC(P')
x_1^\lambda \cdot\bar \omega \oplus \bigl(\CC(P')x_3\cdot\bar
\omega\bigr)^{\muu_{\bar m}}
\]
for a unique $0<\lambda<\bar m$ such that $\lambda a_1+\wt \bar
\omega\equiv 0\mod \bar m$. Since $\dim \Coker \Phi\le 1$, by
\eqref{eq-imp-Phi-gr} we have $2\bar m=s\le 2+1=3$, a contradiction.
\end{pusto}

\begin{pusto}
\textbf{Case $a_2=1$.} As above, $\Phi$ factors through
\[
R:=\bigl(\OOO_{C^\sharp, P^\sharp}/ (x_1^{\bar m}, x_1x_3, x_3^2,
x_2^{\bar m}, x_4)\cdot \bar \omega\bigr)^{\muu_{\bar m}}.
\]
This $R$ is generated by the images of
\[
x_1^ix_2^{\bar m-1-a_1i}, \quad x_3^jx_2^{\bar m-1-a_3j},\qquad 0\le
i\le (\bar m-1)/a_1,\ 1\le j\le (\bar m-1)/a_3.
\]
Therefore,
\[
2\bar m=s\le \dim R+1\le \frac{\bar m-1}{a_1}+1+ \frac {\bar
m-1}{a_3}+1\le \frac{\bar m-1}{1}+ \frac {\bar m-1}{1}+2=2\bar m.
\]
This immediately implies $a_1=a_3=1$. Since $a_1+a_3\equiv 0\mod
\bar m$, we have $\bar m=2$, $s=4$, and $m=8$. Changing
$\ell$-characters \cite[2.5]{Mori-1988} and permuting $x_1$ and
$x_3$, we may assume that $\wt x_2\equiv 1\mod 8$ and $\wt x_1\equiv
1$ or $5 \mod 8$. If $\wt x_1\equiv 1 \mod m$, then $\ow C^\sharp$
is generated by $\ow x_1$ and $\ow x_3$. In particular, $x_1/x_2$ is
constant on $C^\sharp$. This means that $R$ is generated by $x_1$
and $x_3$. Hence, by \eqref{eq-imp-Phi-gr} $4=s\le \dim R+1\le 3$, a
contradiction. Therefore, we have the case \type{IE^\vee}.
\end{pusto}
\end{proof}

\section{Deformations of $\QQ$-conic bundles}
We recall the following
\begin{mtparag}{Proposition ({\cite[1b.8.2]{Mori-1988}}).}
\label{prop-cor-def-f}
Let $(X, C)$ a the $\QQ$-conic bundle germ and let $P\in C$. Then
every deformation of germs $(X,P)\supset (C,P)$ can be extended to a
deformation of $(X, C)$ so that the deformation is trivial outside
some small neighborhood of $P$.
\end{mtparag}

\begin{proof}[Proof \textup(cf. {\rm \cite[11.4.2]{Kollar-Mori-1992}}\textup)]
Let $P_i\in X$ be singular points. Consider the natural morphism
\[
\Psi\colon \mt{Def} X \longrightarrow \prod\mt{Def} (X,P_i).
\]
It is sufficient to show that $\Psi$ is smooth (in particular,
surjective). The obstruction to globalize a deformation in
$\prod\mt{Def} (X,P_i)$ lies in $R^2 f_*T_X$. Since $f$ has only
one-dimensional fibers, $R^2 f_*T_X=0$.
\end{proof}

\begin{propositionm}
Let $f\colon (X,C\simeq \PP^1)\to (Z,o)$ be a $\QQ$-conic bundle
germ. Let $(X_t,C_t)$, $t\in\mathfrak{T}\ni 0$ be an one-parameter
deformation as in Proposition \xref{prop-cor-def-f} and let
$\mathfrak{X}\to \mathfrak{T}$ be the corresponding family so that
$\mathfrak{X}_0=X$. There exists a contraction $\mathfrak f\colon
\mathfrak{X}\to \mathfrak{Z}$ over $\mathfrak{T}\ni 0$ such that
$\mathfrak{Z}_0=Z$ and for all $t\in \mathfrak{T}\ni 0$, $\mathfrak
f_t\colon \mathfrak{X}_t\to \mathfrak{Z}_t$ is a $\QQ$-conic bundle
germ.
\end{propositionm}

\begin{proof}
Consider the base change \eqref{eq-diag--base-change}. Let
$\Gamma'=f'^{-1}(o')$ be the scheme fiber (so that
$\Gamma_{\red}'=C'$) and let
$g^{-1}(P_i)=\{P'_{i,1},\dots,P_{i,s_i}\}$. We claim that arbitrary
deformation in $\mt{Def}(X,P_i)$ determines a $\muu_d$-equivariant
deformation in $\prod_j \mt{Def} (X',P'_{i,j})$. Indeed, the total
space $\mathfrak Q$ of a deformation of a terminal singularity
$(X,P)$ is $\QQ$-Gorenstein (see \cite[\S 6]{Stevens-1988}) and
index-one cover of $\mathfrak Q$ is the total deformation space of
the index-one cover $(X^\sharp,P^\sharp)$ of $(X,P)$. Therefore
every deformation of a terminal singularity of index $m$ is induced
by some $\muu_m$-equivariant deformation of its index-one cover.
This proves our claim. This implies that a deformation in
$\mt{Def}(X,P_i)$ determines a deformation of $X'$ which must be
$\muu_d$-equivariant. Therefore, the cover $X'\to X$ induces a cover
$\mathfrak{X}'\to \mathfrak{X}$ so that
$\mathfrak{X}=\mathfrak{X}'/\muu_d$.

Since $\Gamma'$ is a complete intersection in $\mathfrak{X}'$, the
conormal sheaf $\mathcal{N}_{\Gamma'/\mathfrak{X}'}^*$ is locally
free. We have the exact sequence
\begin{equation}
\label{eq-seq-N}
\begin{CD}
0@>>> \mathcal{N}^*_{X'/\mathfrak{X}'}|_{\Gamma'} @>>>
\mathcal{N}^*_{\Gamma'/\mathfrak{X}'} @>>>
\mathcal{N}^*_{\Gamma'/X'} @>>> 0
\\
@.@|@.@|
\\
@.\OOO_{\Gamma'}@.@.\OOO_{\Gamma'}^{\oplus 2}
\end{CD}
\end{equation}
Since $\mt{Ext}^1(\OOO_{\Gamma'}^{\oplus 2},\OOO_{\Gamma'})=
H^1(\OOO_{\Gamma'}^{\oplus 2})=0$, the sequence splits.

Therefore the germ $\mathfrak{D}$ of the Douady space of
$\mathfrak{X}'$ at $[\Gamma']$ is smooth, where $[\Gamma']$ is the
point representing $\Gamma'$. Let $\mathfrak{U}\to \mathfrak{D}$ be
the corresponding universal family. There is a natural embedding
$\mathfrak{U}\subset \mathfrak{X}'\times \mathfrak{D}$ such that
$\mathfrak{U}\to \mathfrak{D}$ is induced by the projection
$\mathfrak{X}'\times \mathfrak{D}\to \mathfrak{D}$. Thus we have the
following diagram:
\[
\xymatrix{\mathfrak{X}'\times\mathfrak{D}\ar[d]^{\mt{pr}_2}
&\supset&\mathfrak{U}\ar[r]^(0.35){\alpha}\ar[dll]^{\beta}
&\mathfrak{X}'\supset\Gamma'
\\
\mathfrak{D}&&&}
\]
The natural embedding $\Gamma'=\Gamma'\times [\Gamma']\subset
\mathfrak{U}$ induces an isomorphism $\alpha|_{\Gamma'}\colon
\Gamma'\to \Gamma'$. Further, $\Gamma'=\Gamma'\times
[\Gamma']\subset \mathfrak{U}$ is a fiber of $\beta$, so
$\mathcal{N}^*_{\Gamma'/\mathfrak{U}}$ is locally free and
isomorphic to $\OOO_{\Gamma'}^{\oplus 3}$. Hence,
\[
d\alpha\colon \mathcal{N}^*_{\Gamma'/\mathfrak{X}'}\to
\mathcal{N}^*_{\Gamma'/\mathfrak{U}}.
\]
is an isomorphism. Shrinking $\mathfrak{U}\supset \Gamma'$ and
$\mathfrak{X}'\supset \Gamma'$ we may assume that $\alpha$ is an
isomorphism. This induces a $\muu_d$-equivariant contraction
morphism $\Phi=\alpha^{-1}\beta\colon \mathfrak{X}'\to \mathfrak{D}$
such that $\Phi(\Gamma')$ is a point. Put $\mathfrak Z:=\mathfrak
D/\muu_d$. Since the morphism $\mathfrak p\colon \mathfrak X\to
\mathfrak T$ maps $X$ to $0$, by shrinking $\mathfrak{X}$ we may
assume that $\mathfrak p$ is constant on fibers of $\mathfrak f$.
Then $\mathfrak p$ defines $\mathfrak Z\to \mathfrak T$. We obtain
the following diagram
\[
\xymatrix{\mathfrak X'\ar[d]^{\mathfrak f'}\ar[rr]^{\mathfrak r}&
&\mathfrak X\ar[d]^{\mathfrak f}\ar@/_1.4pc/[ddl]_(0.4){\mathfrak p}
\\
\mathfrak D\ar@{->}'[r][rr]^(-0.4){\mathfrak q}\ar[rd]&&\mathfrak
Z\ar[ld]
\\
&\mathfrak T&}
\]
We have
\[
\mathfrak f_*\OOO_{\mathfrak X}= \mathfrak f_*(\mathfrak
r_*\OOO_{\mathfrak X'})^{\muu_d}= (\mathfrak f_*\mathfrak
r_*\OOO_{\mathfrak X'})^{\muu_d}= (\mathfrak q_*\mathfrak
f_*'\OOO_{\mathfrak X'})^{\muu_d}= (\mathfrak q_*\OOO_{\mathfrak
D})^{\muu_d}=\OOO_{\mathfrak Z}.
\]
Therefore $\mathfrak f$ has connected fibers. Clearly, $\mathfrak X$
is $\QQ$-Gorenstein and $-K_{\mathfrak X}$ is $\mathfrak f$-ample.
\end{proof}

\begin{remark}
In general, it is not true that $f_{t}^{-1}(o_t)_{\red}=C_{t}$. It
is possible that $f_{t}^{-1}(o_t)$ is reducible and $C_t$ is one of
its components. In this case, $(X_{t},C_{t})$ is an extremal
neighborhood by Remark \xref{rem-prel-extr-nbd}.
\end{remark}

\section{The case where $X$ is not locally primitive }
\label{sect-imprimitive}
In this section we consider the 
case where $X$ is not locally primitive.
We
classify configurations of singular points and prove Theorem
\ref{th-main-1} except for one case when a $\QQ$-conic bundle germ
has two non-Gorenstein points. Some weaker results were obtained in
\cite{Prokhorov-1996b}.

\begin{mparag}{Notation.}
\label{not-imprim}
Let $f\colon (X,C)\to (Z,o)$ be a $\QQ$-conic bundle germ, where $C$
is irreducible. Assume that $(X,C)$ contains an imprimitive point
$P$. Let $m$, $\bar m$ and $s$ be the index, the subindex and the
splitting degree of $P$, respectively.
\end{mparag}

\begin{parag}{}
\label{zam-imp-Gor-fac}
First we note that if $\bar m=1$, then by Lemma
\ref{lemma-int-non-Gor} $X'$ is Gorenstein and we have the case
\ref{item-main-th-impr-barm=1} by \ref{prop-Gor-quot}. From now on
we assume that $\bar m>1$ (in particular, $P$ is not of type \type{ID^\vee}).
\end{parag}

\begin{tparag}{Lemma (cf. {\cite[Th. 6.1 (ii)]{Mori-1988}}).}
\label{lemma-IA-dual-IA-2}
In notation of \xref{not-imprim}, assume that $X$ has a singular
point $Q\neq P$. Then $\Sing X=\{P,\, Q\}$ and $Q$ is of type
\type{IA} or \type{III}. If $Q$ is of type \type{IA}, then
$\size_Q=1$.
\end{tparag}

Recall \cite[4.5]{Mori-1988} that a point $P\in X$ is said to be \emph{ordinary}  
iff $(X,P)$
is either an ordinary double point or a cyclic quotient
singularity.

\begin{proof}
Assume that $(X,C)$ has two more singular points $Q$ and $R$.
Then by Proposition \ref{prop-imp-types},
$P$ is of type 
\type{IA^\vee}, \type{IC^\vee}, or \type{II^\vee}.
By Proposition \ref{prop-first-prop-imp} both $Q$ and $R$ are primitive. 
Replace $(X,C)$ with $L$-deformation
\cite[Prop.-Def. 4.7]{Mori-1988} 
so that $P$, $Q$, $R$ are
ordinary
(cf. \cite[Rem. 4.5.1]{Mori-1988}). 
If this new $(X,C)$ is an extremal neighborhood, the
assertion follows by \cite[Th. 6.1 (ii)]{Mori-1988}. Thus we may
assume that $(X,C)$ is a $\QQ$-conic bundle germ. Consider the cover
$g\colon (X',C') \to (X,C)$ from \eqref{eq-diag--base-change}. By
Lemma \ref{lemma-int-non-Gor} we may assume that $P'\in X'$ is not
Gorenstein (otherwise $(X',C')$ is a standard Gorenstein conic
bundle germ and then $(X,C)$ is as in
\ref{item-main-th-impr-barm=1}, see Proposition
\ref{prop-Gor-quot}). 
Let $C'_i \subset C'$ be any irreducible component.
By \ref{cor-imp-top-1} (iii), $(X', C'_i)$ is an extremal neighborhood
with at least three singular points.
Then by \cite[(2.3.2)]{Mori-1988} $\deg \gr_{C'_i}^1\OOO\le
-2$. Hence, $H^0(\gr_{C'_i}^1\OOO)=0$ \cite[Remark
2.3.4]{Mori-1988}. This implies that $H^0(\gr_{C'}^1\OOO)\subset
\bigoplus\limits_i H^0(\gr_{C'_i}^1\OOO)=0$. Therefore,
$H^0(I_{C'}^{(2)})= H^0(I_{C'})$ and $H^0(\OOO_{X'}/
I_{C'}^{(2)})=\CC$. By Corollary
\ref{cor-lemma-sect-1-gen-omega-OOO-ne0} we have
$H^1(\gr_{C'}^1\omega)=0$. Therefore,
$H^1(\gr_{C}^1\omega)=H^1(\gr_{C'}^1\omega)^{\muu_d}=0$. 
With this extra condition, the proof of \cite[Th. 6.1 (i)]{Mori-1988}
(resp. \cite[Th. 6.1 (ii)]{Mori-1988}) works
if $P$ is of type \type{IC^{\vee}} (resp. \type{IA^{\vee}}).
Thus $\Sing X=\{P,\, Q\}$.

Now assume that $Q$ is not of type \type{III}. Consider the
splitting cover $g^\flat$ from \eqref{eq-cor-imp-top-diag}. Each
$(X^\flat,C_i^\flat)$ is an extremal neighborhood having two
non-Gorenstein points: $P^\flat$ and $Q^\flat_i$. By \cite[Th. 6.7,
Th. 9.4]{Mori-1988} $Q^\flat_i$ is of type \type{IA} with $\size=1$
and so is $Q$.
\end{proof}

\begin{propositionm}{}
\label{prop-imp-fiber}
If $C'= f'^{-1}(o')$ \textup(and $\bar m>1$\textup) or, equivalently
if we have a point of type
\type{IE^\vee}, then $f$ is as in
\xref{item-main-th-impr-barm=2-s=4}.
\end{propositionm}
\begin{proof}
We note that $C'= f'^{-1}(o')$ \textup(and $\bar m>1$\textup) if and
only if we are in the case \type{IE^\vee} by Proposition
\ref{prop-imp-types}. In some (non-normalized) coordinate system,
$C^\sharp\subset \CC^{4}_{y_1,y_2,y_3,y_4}$ is a complete
intersection given by
\begin{equation}
\label{eq-imp-exc-8-eq}
y_1^2-y_2^2=y_1y_2-y_3^2=y_4=0,
\end{equation}
where $\wt (y)\equiv (5,1,3,0) \mod 8$. Thus we may fix an embedding
$C^\sharp \subset \CC^3_{y_1,y_2,y_3}$ and $X^\sharp\subset
\CC^{4}_{y_1,\dots,y_4}$. Let $(u,v)$ be $\muu_8$-semi-invariant
coordinates in $Z'=\CC^2$. Since $C'=f'^{-1}(o')$, we may regard
$u,\, v$ as $\muu_8$-semi-invariant generators of the ideal of
$C^\sharp$ in $X^\sharp$. Therefore the ideal of $C^\sharp$ in
$\CC^{4}_{y_1,\dots,y_4}$ has two systems of semi-invariant
generators:
\[
y_1^2-y_2^2,\ y_1y_2-y_3^2, \ y_4\quad \mbox{and} \quad u,\, v,\,
\phi.
\]
Up to permutation of $u$ and $v$ we may assume that
\[
\begin{array}{l}
\wt u \equiv \wt y_1^2\equiv 2, \quad \wt v \equiv \wt y_1y_2 \equiv
-2 \mod 8,
\\[7pt]
\phi= (\mbox{unit}) y_4 +(y_1^2-y_2^2)\phi_1 + (y_1y_2-y_3^2) \phi_2
\end{array}
\]
because $\wt y_1^2,\ \wt y_3^2 \not \equiv \wt y_4 \mod 8$. In
particular, $X^\sharp$ is smooth and $(X,P)$ is a cyclic quotient of
type $\frac 18(5,1,3)$. Hence $(X',P')$ is a singularity of type
$\frac12(1,1,1)$ and coordinates $y_1$, $y_2$, $y_3$ can be regarded
as sections of $|-K_{X'}|$ on $X'$. Note that the linear system
$|-K_{X'}|$ has a unique base point $P'$ and $|-2K_{X'}|$ is base
point free. Let $z$ be a section of $|-2K_{X'}|$. Then
$y_1,y_2,y_3,z$ define a map $\vartheta\colon X'\dashrightarrow
\PP\times \CC^2$, where $\PP:=\mt{Proj}
\CC[y_1,y_2,y_3,y_4]=\PP(1,1,1,2)$. Since this $\vartheta$ is regular
on each component of $C'$ and on the tangent space to $C'$ at $P'$,
it is an embedding. Therefore, $X'$ can be naturally embedded into
$\PP\times \CC^2$ and by \eqref{eq-imp-exc-8-eq} the defining
equations are of the form
\[
\begin{array}{lll}
y_1^2-y_2^2&=&u \psi_1+v\psi_2,
\\[7pt]
y_1y_2-y_3^2&=&u \psi_3+v\psi_4,
\end{array}
\]
where $\psi_i=\psi_i(y_1,y_2,y_3,z, u,v)$. 
This proves our proposition.
\end{proof}

\begin{corollary}
\label{cor-impr-grw-21cases}
In the notation of \xref{not-imprim}  the following are equivalent\textup:\
\textup{(i)}\ \nolinebreak $\gr_C^0\omega \not\simeq \OOO(-1)$,\
\textup{(ii)}\ we are in the case \xref{item-main-th-impr-barm=1}, and\
\textup{(iii)}\ $P$ is  a point of type \type{ID^{\vee}}.
\end{corollary}

\begin{proof}
Assume that
$\gr_C^0\omega \not\simeq \OOO(-1)$. 
Then 
by \ref{cor-prop-grw=2-2-points-prim} we see that $f'^{-1}(o')=C'$. 
If $\bar m>1$, then we are in the case
\xref{item-main-th-impr-barm=2-s=4} by Proposition
\ref{prop-imp-fiber}, in which case we have $w_P(0)=1/2$ and hence
$\gr_C^0\omega= \OOO(-1)$ by \eqref {eq-grw-w} and Lemma
\ref{lemma-KC}. Hence $\bar m=1$, then we are in the case \ref{item-main-th-impr-barm=1}
as explained in \ref{zam-imp-Gor-fac}. Thus (i) implies (ii) and (iii).

If $P$ is  of type \type{ID^{\vee}}, then $\bar m=1$ and again by 
\ref{zam-imp-Gor-fac} we are in the case \xref{item-main-th-impr-barm=1}.
Then by \cite[(2.10)]{Mori-1988} we have $w_P(0)\ge 1$ and so
$\deg \gr_C^0\omega<-1$ (see \eqref{eq-grw-w}).
\end{proof}

\begin{propositionm}{}
\label{prop-IC-dual}
In notation of \xref{not-imprim} $(X,C)$ has no type \type{IC^\vee}
points.
\end{propositionm}

\begin{proof}
Assume that $P\in (X,C)$ is a type \type{IC^\vee} point. Recall that
$s=2$ and $\bar m$ is even $\ge 4$ in this case. Then $(X,C)$ has at
most one more (primitive) singular point. Applying $L$-deformation
we may assume one of the following:

\begin{parag}{}
\label{imp-IC-Case-a}
$P$ is the only singular point of $X$, or
\end{parag}

\begin{parag}{}
\label{imp-IC-Case-b}
$(X,C)$ has one more ordinary singular point $Q$ of index $n>1$.
\end{parag}

By \cite[Th. 6.1 (i)]{Mori-1988} this new $(X,C)$ is a conic bundle
germ. Following the proof of \cite[(i) Th. 6.1]{Mori-1988} we get
$H^1(\OOO_{X}/I_{C}^{(2)}\totimes \omega_{X})\neq 0$. Hence,
$H^1(\OOO_{X'}/I_{C'}^{(2)}\totimes \omega_{X'})\neq 0$.

\begin{parag}{}
\label{imp-IC-fiber}
Let $V:= \Spec \OOO_{X'}/I_{C'}^{(2)}$. By Theorem
\xref{theo-sect-1-gen-fiber-m} we have $f'^{-1}(o')\subset V$.
Moreover, as in the proof of \ref{cor-lemma-sect-1-gen-omega-OOO-ne0}
one can see that $V$ is not a local complete intersection at the
general point. Therefore, $f'^{-1}(o')\neq V$. Since $V\equiv 3C'$
(as a cycle), we have
\[
2=-K_{X'}\cdot f'^{-1}(o') < -K_{X'} \cdot V=-3K_{X'} \cdot C'.
\]
Taking account of $-K_{X'}\cdot C'=d^2/mn$ (see \ref{lemma-KC}) we
get $2nm<3d^2$, where we put $n=1$ in the case \ref{imp-IC-Case-a}.
Recall that $s=2$. Write $m=s\bar m=2\bar m$, $n=rn'$ and $\bar
m=r\bar m'$, where $r=\gcd(\bar m, n)$. Then $d=2r$ (see Corollary
\ref{cor-cyclic-quo-new} (iii)) and $n'\bar m'<3$. Note that $\bar
m'$ is the index of $P'$. If $\bar m'=1$, then $X'$ is Gorenstein by
\ref{lemma-int-non-Gor} and $(X,P)$ cannot be of type \type{IC^\vee}
by Proposition \ref{prop-Gor-quot}. So, $\bar m'=2$, $n'=1$, and
$\bar m=2n$. In particular, $n>1$ and the case \ref{imp-IC-Case-a}
is impossible.
\end{parag}

Note that $(X,P)$ is a cyclic quotient singularity by \cite[Lemma
4.4]{Mori-1988}. Thus we may assume that $X^\sharp
=\CC^3_{x_1,x_2,x_4}$ and $C^\sharp$ is given by the equations
$x_4=x_2^2-x_1^{m-2}=0$. Thus $C'$ near $P'$ is isomorphic to
$\{x_4=x_2^2-x_1^{m-2}=0\}/\muu_2(1,1,1)$. Putting $w_1=x_1^2$,
$w_2=x_2^2$, $w_3=x_1x_2$ we get that near $P'$ the curve $C'\subset
\CC^3_{w_1,w_2,w_3}$ can be given by two equations $w_2=w_1^{\bar
m-1}$ and $w_1w_2=w_3^2$. Eliminating $w_2$ we obtain
$C':=\{w_3^2=w_1^{2n}\}$. It is easy to see that $C'$ has an
ordinary double point at the origin only if $n=1$. This contradicts
$p_a(C')=0$.
\end{proof}

The following lemma was proved in \cite[\S 3]{Prokhorov-1996b}.
However it was implicitly assumed in the proof that $X$ is
$\QQ$-factorial. Below is a corrected version.

\begin{lemma}
\label{lemma-imp-s=2}
In notation of \xref{not-imprim} assume that $X$ is $\QQ$-factorial.
Then $s=2^k$. If furthermore $X$ has two non-Gorenstein points, then
$s=2$.
\end{lemma}

\begin{proof}
Write $d=sr$ and $s=2^kq$, where $q$ is odd.
We will derive a contradiction
assuming $q >1$.
Consider the quotient
$X^{\dagger}/Z^{\dagger}$ of $X'/Z'$ from
\eqref{eq-diag--base-change} by $\muu_{2^kr}\subset \muu_d$:
\[
\begin{CD}
X^{\dagger}@>{g^{\dagger}}>>X
\\
@V{f^{\dagger}}VV @V{f}VV
\\
Z^{\dagger}@>{h^{\dagger}}>> Z
\end{CD}
\]
where $h^{\dagger}\colon Z^{\dagger}\to Z$ is a $\muu_{q}$-cover.
Then $C^{\dagger}:={g^{\dagger}}^{-1}(C)$ has $q$ irreducible
components because $X^{\dagger}$ is a $\muu_{2^k}$-quotient of the
splitting cover \ref{cor-imp-top-1} and therefore, $\rho(X^{\dagger}/Z^{\dagger})=q$ by
Corollary \ref{cor-C-pa=0}. There is a curve $V\subset Z^{\dagger}$
such that $f^{\dagger -1}(V)$ has exactly two components, say $E_1$
and $E_{1'}$. For a general point $z\in V$ the preimage
${f^{\dagger}}^{-1}(z)$ is a reducible conic, so
${f^{\dagger}}^{-1}(z)=\ell_1+\ell_2$. Consider the orbit $\{E_1,
E_2, \dots, E_t\}$ of $E_1$ under the action of $\muu_{q}$.
Obviously, every $f^{\dagger}(E_i)$ is a curve on $Z^{\dagger}$.
Further, $\sum_{i=1}^t E_i\qq g^{\dagger *} M$, where $M$ is a Weil
$\QQ$-Cartier divisor on $X$. On the other hand, $\rho(X/Z)=1$ and
$M$ is $f$-vertical. Hence, $M\qq 0$ and $\sum E_i\qq 0$. We can
choose components $\ell_1,\ell_2\subset {f^{\dagger}}^{-1}(z)$,
$z\in V$ so that $\ell_1\cdot E_1<0$ and $\ell_1\cdot E_{1'}>0$.
This gives us $\ell_1\cdot E_i>0$ for some $E_i\in\{E_1, E_2, \dots,
E_t\}$. Then $E_{1'}=E_i$, i.e., there exists $\sigma\in\muu_{q}$
such that $\sigma(E_1)=E_{1'}$. From the symmetry we get that the
orbit $\{E_1, E_2, \dots, E_t\}$ may be divided into pairs of
divisors $E_j$, $E_{j'}$ such that
${f^{\dagger}}(E_j)={f^{\dagger}}(E_{j'})$ is a curve. Thus both $t$
and $q$ are even. This proves the first statement.

Now assume that $X$ has two non-Gorenstein points and let $s=2^k$,
$k\ge 2$. Consider the quotient $(X'',C'')$ of $(X',C')$ by
$\muu_{2^{k-1}}$. Then the central fiber $C''$ is reducible and
every germ $(X'',C_i'')$ is an extremal neighborhood having two
non-Gorenstein points: imprimitive and primitive. By the
classification \cite[Th. 6.7, 9.3]{Mori-1988} this is impossible.
\end{proof}

\begin{mtparag}{Proposition (cf. {\cite[Th. 6.1 (iii)]{Mori-1988}}).}
\label{prop-imp-IA-size}
Notation as in \xref{not-imprim}. Assume that $(X,C)$ has one more
non-Gorenstein point $Q$. Then $P$ is of type \type{IA^\vee},
$\size_P=1$, and $w_P(0)\ge 2/3$.
\end{mtparag}

\begin{corollary}
\label{cor-imp-IA-size}
In the above notation we have $w_Q(0)<1/3$. In particular, the index
of $Q$ is $\ge 4$.
\end{corollary}

\begin{proof}
\ref{cor-imp-IA-size} immediately follows from
\ref{prop-imp-IA-size} and \ref{cor-impr-grw-21cases}. 
We assume that $\size_P\ge 2$ or $w_P(0)<
2/3$ and we will derive a contradiction. Let $n$ be the index of
$Q$. By Lemma \ref{lemma-IA-dual-IA-2} $\Sing X=\{P,\, Q\}$ and $Q$
is of type \type{IA} and by 
Propositions \ref{prop-imp-fiber}, \ref{prop-IC-dual}
and 
Corollary \ref{cor-impr-grw-21cases}, $P$ is of
type \type{IA^\vee} or \type{II^\vee}. Replacing $(X,C)$ with
$L$-deformation, we may assume that $X$ has only ordinary points (in
particular, $P$ is of type \type{IA^\vee}). If this new $(X,C)$ is
an extremal neighborhood, the assertion follows by \cite[Th. 6.1
(iii)]{Mori-1988}. Thus we may assume that $(X,C)$ is again a
$\QQ$-conic bundle germ.

If $H^1(\OOO_{X}/I_{C}^{(2)}\totimes \omega_{X})= 0$, then following
the proof of \cite[(iii) Th. 6.1]{Mori-1988} we derive a
contradiction. Hence, in notation of \eqref{eq-diag--base-change} we
have $H^1(\OOO_{X'}/I_{C'}^{(2)}\totimes \omega_{X'})\neq 0$. Let
$V:= \Spec \OOO_{X'}/I_{C'}^{(2)}$. As in \ref{imp-IC-fiber} by
Theorem \xref{theo-sect-1-gen-fiber-m} $f'^{-1}(o')\subset V$, and
so
\[
2=-K_{X'}\cdot f'^{-1}(o') < -K_{X'} \cdot V =-3K_{X'} \cdot C'.
\]
Taking account of $-K_{X'}\cdot C'=d^2/mn$ (see \ref{lemma-KC}) we
obtain
\[
2mn< 3d^2.
\]
Since $X$ has only ordinary points of index $>1$, $X$ is
$\QQ$-factorial. By Lemma \xref{lemma-imp-s=2} $s=2$. Write $m=s\bar
m=2\bar m$, $n=rn'$ and $\bar m=r\bar m'$, where $r=\gcd(\bar m,
n)$. By Corollary \ref{cor-cyclic-quo-new} (iii) we have $d=sr=2r$.
Then $n'\bar m'<3$. Since $\bar m'$ is the index of $P'$, we may
assume that $\bar m'>1$ (otherwise by Lemma \ref{lemma-int-non-Gor}
and Proposition \ref{prop-Gor-quot} we have the case
\ref{item-main-th-impr-barm=1}). Therefore, $\bar m'=2$, $n'=1$, and
$\bar m=2n$.

We may assume that $X^\sharp =\CC^3_{x_1,x_2,x_3}$ at $P^\sharp$,
the curve $C^\sharp$ is given by the equations
$x_3=x_2^{a_1s}-x_1^{a_2s}=0$, and $(X',P')=
\CC^3_{x_1,x_2,x_3}/\muu_2(1,1,1)$. Putting $w_1=x_1^2$,
$w_2=x_2^2$, $w_3=x_1x_2$ we get that near $P'$ the curve $C'\subset
\{x_3=0\}/\muu_2=\CC^2_{x_1,x_2}(1,1)/\muu_2$ can be given by two
equations: $w_2^{a_1}=w_1^{a_2}$ and $w_1w_2=w_3^2$. We claim that
$a_1=a_2=1$. Indeed, assume for example that $a_1>1$. Since
$p_a(C')=0$, $C'$ has an ordinary double point at the origin. Hence,
$a_2=1$. Eliminating $w_1$ we get the following equation for $C'$:
$w_2^{a_1+1}=w_3^2$. Again the origin is an ordinary double point
only if $a_1=1$, a contradiction. Thus, $a_1=a_2=1$ and
$w_P(0)=1-1/\bar m$ by \cite[Th. 4.9]{Mori-1988}. This gives
\[
w_Q(0)=1-w_P(0)+K_X\cdot C=1/\bar m-1/\bar m=0,
\]
a contradiction. Hence we have $\size_P=1$ and $w_P(0)\ge 2/3$. If
$P$ is of type \type{II^\vee}, then $w_P(0)=1/2$ (see \cite[Th.
4.9]{Mori-1988}). So, $P$ is of type \type{IA^\vee}.
\end{proof}

\begin{propositionm}
\label{prop-75}
In notation of \xref{not-imprim} assume that $P$ is of type
\type{IA^\vee}. Then $P$ is of index $4$, splitting degree $2$ and
subindex $2$. Moreover, $(X,P)$ is a cyclic quotient and
\[
\begin{array}{cccccc}
&x_1&x_2&x_3&x_4&
\\
\wt&1&-1&-1&0&\mod 4
\\
\ord&1&1&1&2&
\end{array}
\]
\end{propositionm}
\begin{proof}
By \ref{cor-2nonGor-IAIA} below $P$ is the only non-Gorenstein point
on $X$.
\footnote{In \ref{not-2-points-m} -- \ref{cor-2nonGor-IAIA},
no results in  \ref{prop-75} 
--  \ref{prop-931} 
are used when 
$P$ is imprimitive.
Thus the 
back
reference \ref{cor-2nonGor-IAIA} 
here does
not cause any trouble.}
Since $-K_X\cdot C=1/\bar m$, $w_P(0)=1-1/\bar m$. 
Hence we  have
$a_2=1$ by \cite[Th. 4.9.(i)]{Mori-1988}.
The general member $F\in
|-K_{(X,P)}|$ has only Du Val singularity of type $A_{mk-1}$, $k\in
\ZZ_{>0}$. It is easy to see that $F^{\sharp}$ is given by $x_2=0$,
so $F\cdot C=1/\bar m$. Hence, $K_X+F$ is a numerically trivial
Cartier divisor. Since $\Pic X\simeq H^2(X,\ZZ)\simeq \ZZ$,
$K_X+F\sim 0$. Thus, the general member $F\in |-K_X|$ does not
contain $C$ and has only Du Val singularity of type $A_{mk-1}$.
Consider the double cover $f|_F\colon (F,P)\to (Z,o)$. Diagram
\eqref{eq-diag--base-change} induces the following
\[
\begin{CD}
(F',P') @>{g_{F'}}>> (F,P)
\\
@V{f'_{F'}}VV @V{f_F}VV
\\
(Z',o') @>{h}>> (Z,o)
\end{CD}
\]
where $F'\in |-K_{X'}|$, $P'=g^{-1}(P)$, $F'\cap C'=\{P'\}$, and
$g_{F'}$ is \'etale outside of $P'$. Since $Z'\to Z$ is of degree
$s$, $s$ divides $mk$ and $(F',P')$ is of type $A_{n-1}$, where
$n=mk/s=\bar mk$. We see $n=2$ because otherwise we have a
contradiction by Lemma \ref{lemma-vspom-cover-DuVal} below. Thus
$n=2$ and $\bar m=2$ (recall that $\bar m>1$ by the assumption of \type{IA^\vee}
\ref{prop-imp-types}). 
In this case, by Corollary
\ref{cor-eqiv-imp} $s=4$ or $2$. If $s=4$, then $-K_{X'}\cdot
C'=s/\bar m=2$ (see Lemma \ref{lemma-KC}). Hence, $C'=f'^{-1}(o')$
and we have the case \ref{item-main-th-impr-barm=2-s=4} by
Proposition \ref{prop-imp-fiber}. But then $P$ is not of type
\type{IA^\vee}, a contradiction. Hence $s=2$ and the rest is easy.
\end{proof}

\begin{lemma}
\label{lemma-vspom-cover-DuVal}
Let $(S,Q)$ be a Du Val singularity of type $A_{n-1}$, $n\ge 3$ and
let $\pi \colon (S,Q)\to (\CC^2,0)$ be a double cover. Assume that
$\muu_d$ acts on $(S,Q)$ and $(\CC^2,0)$ freely in codimension one
and so that $\pi$ is $\muu_d$-equivariant. Then the quotient
$(S,Q)/\muu_d$ cannot be Du Val of type $A$.
\end{lemma}
\begin{proof}
Let $R\subset \CC^2$ be the branch divisor of $\pi$. Since $(S,Q)$ is
of type $A$, the equation of $R$ must contain a quadratic term.
Hence, in some $\muu_d$-semi-invariant coordinates $u,v$ in $\CC^2$,
the curve $R$ can be given by $u^2+v^n=0$. In this case, there is a
$\muu_d$-equivariant embedding $(S,Q)\hookrightarrow
(\CC^3_{u,v,w},0)$ such that $S$ is given by $w^2=u^2+v^n$ and $w$
is a semi-invariant. Assume that $S/\muu_d$ is Du Val. Since
$K_{S/\muu_d}$ is Cartier, we have $\wt (uvw)=\wt w^2=\wt u^2=\wt
v^n$. This implies $\wt w=\wt (uv)$ and $\wt v^2=0$. Since the
action of $\muu_d$ on $\CC^2$ is free in codimension one, $d=2$ and
$n$ is even. So, $n=2l$ for some $l\ge 2$. Hence, $S/\muu_d$ is a
quotient of $\{u^2+v^{2l}=w^2\}$ by $\muu_2(1,1,0)$. But if $l\ge
2$, this quotient is not of type $A$, a contradiction.
\end{proof}

\begin{propositionm}
\label{prop-imprim-d2-}
In notation of \xref{not-imprim} assume that $P$ is of type
\type{IA^\vee} \textup(resp. \type{II^\vee}\textup). Then $f\colon
(X,C)\to (Z,o)$ is as in \xref{item-main-th-impr-barm=2-s=2-cycl}
\textup(resp. \xref{item-main-th-impr-barm=2-s=2-cAx/4}\textup).
\end{propositionm}

\begin{proof}
In the case \type{II^\vee} $X$ has no other non-Gorenstein points by
\ref{prop-imp-IA-size}. Then applying \eqref{eq-diag--base-change}
we will see that $X/Z$ is the quotient of an index-two $\QQ$-conic
bundle $f'\colon (X',C')\to (Z',o')$ by $\muu_2$. The components of
the central curve $C'$ are permuted, so $C'$ has two components of
the same multiplicity. Hence $X'/Z'$ is in the case
\xref{cla-index-2-2-2}. The action on $X'$ is described in
\xref{lem-cla-index-2-2-2-mu2-action-2}.
\end{proof}

\begin{remark}
\label{rem-imp-ge}
We have treated all the types 
\type{IA^\vee}--\type{II^\vee}
of imprimitive points.
Finally we note that the existence of a good anicanonical divisor
(Proposition \ref{propo-main-ge}) in the imprimitive cases
\ref{item-main-th-impr-barm=2-s=4} -
\ref{item-main-th-impr-barm=2-s=2-cAx/4} can be shown exactly as in
\cite[7.3]{Mori-1988} (see also \cite{Prokhorov-1996b}).
\end{remark}

\begin{mparag}{Examples.}
Below we propose explicit examples of $\QQ$-conic bundles as in
\ref{item-main-th-impr-barm=2-s=4},
\ref{item-main-th-impr-barm=2-s=2-cycl} and
\ref{item-main-th-impr-barm=2-s=2-cAx/4}.
\begin{example}
\label{example-imp-4}
Under the notation of \xref{item-main-th-impr-barm=2-s=4} consider
the subvariety $X'$ defined by
\[
\left\{
\begin{array}{lll}
y_1^2-y_2^2&=&uy_4
\\[7pt]
y_1y_2-y_3^2&=&vy_4.
\end{array}
\right.
\]
The projection $f'\colon X'\to\CC^2$ is a $\QQ$-conic bundle of
index $2$ (see \xref{cla-index-2-4}). Then $X'/\muu_{4}\to
\CC^2/\muu_4$ is a $\QQ$-conic bundle
with an imprimitive point as in
\xref{item-main-th-impr-barm=2-s=4}.
The singular point is unique and is of type \type{IE^\vee}.
\end{example}

\begin{example}
\label{example-imp-IA}
Let $X'\subset\PP (1,1,1,2) \times\CC^2$ be the subvariety given by
the equations
\[
\left\{
\begin{array}{lll}
y_1^2-y_2^2&=&uy_4
\\[7pt]
y_3^2&=&vy_4+u^2y_2^2+\lambda uy_1y_2,\quad \lambda\in \CC.
\end{array}
\right.
\]
Consider the action of $\muu_{2}$ on $X'$:
\[
y_1\mapsto y_1,\quad y_2\mapsto - y_2,\quad y_3\mapsto y_3,\quad y_4\mapsto
-y_4,\quad u\mapsto -u,\quad v\mapsto -v.
\]
Then $X:=X'/\muu_{2}\to\CC^2/\muu_{2}$ is a
$\QQ$-conic bundle with an imprimitive point
as in \ref{item-main-th-impr-barm=2-s=2-cycl}. It
has a singularity of type \type{IA^\vee} which is the
cyclic quotient $\frac14(1,-1,1)$. If
$\lambda=0$, then $X$ also has a (Gorenstein) 
ordinary double point.
\end{example}

\begin{example}
\label{example-imp-II}
Let $X'\subset\PP(1,1,1,2)\times\CC^2$ be the subvariety given by
the equations
\[
\left\{
\begin{array}{lll}
y_1^2-y_2^2&=&u^3y_4+vy_4
\\[7pt]
y_3^2&=&vy_4+u^2y_2^2+\lambda uy_1y_2,\quad \lambda\in \CC.
\end{array}
\right.
\]
Define the action of $\muu_{2}$ on $X'$ as in \ref{example-imp-IA}.
Then $X:=X'/\muu_{2}\to\CC^2/\muu_{2}$ is a 
$\QQ$-conic bundle
with an imprimitive point  as in \ref{item-main-th-impr-barm=2-s=2-cAx/4}.
The non-Gorenstein point is of type \type{II^\vee}.
It is the only singular point 
if $\lambda\neq 0$. If $\lambda= 0$,
then $X$ has one more singular point which is of type \type{III}.

\end{example}
\end{mparag}

In case of a $\QQ$-conic bundle with an imprimitive point,
the proofs of Theorems \ref{th-main-1} 
and 
\ref{theo-main-DuVal}  are completed
here modulo the arguments in \ref{not-2-points-m} -- \ref{cor-2nonGor-IAIA}.

\section{The case where $X$ is locally primitive. Possible singularities.}
In this section we consider locally primitive $\QQ$-conic bundles.
The main result is summarized in Theorem \ref{th-prim-conf-summ}.

\begin{mparag}{Notation.}
\label{not-prim}
Let $f\colon (X,C\simeq \PP^1)\to (Z,o)$ be a locally primitive
$\QQ$-conic bundle germ. Let $P\in (X,C)$ be a (primitive)
non-Gorenstein singular point and let $m\ge 2$ be its index. We may
assume that $\gr_C^0\omega\simeq \OOO_C(-1)$ (see
\ref{cor-prop-grw=2-2-points-prim}).
\end{mparag}

\begin{lemma}
\label{lemma-3SingPts}
There are at most $3$ singular points of $X$ on $C$.
\end{lemma}
\begin{proof}
By $\deg \gr_C^0\omega<0$ and  \eqref{eq-grO-iP1-2}, we have
$\sum_Q i_Q(1) \le 3$,
and the proof of \cite[6.2(i)]{Mori-1988} works.
\end{proof}

\begin{lemma}
\label{lemma-IBICIIB}
If $P$ is a point of type \type{IB} or \type{IC}, then the base
$(Z,o)$ is smooth. In the case \type{IIB}, $(Z,o)$ is either smooth
or Du Val of type $A_1$.
\end{lemma}
\begin{proof}
Assume that $(Z,o)$ is singular and consider the base change as in
\eqref{eq-diag--base-change}. Let $P'\in g^{-1}(P)$ and let $m'$ be
the index of $(X',P')$. We note that $(X^\sharp,P^\sharp)$ is also
the index-one cover of $(X',P')$. 
Clearly, $m'$ divides $m$.

We claim that $m'<m$. Suppose $m'=m$. Then 
the Galois cover $X' \to X$ is
\'etale at $P'$ and $X'$ has at least two points of the
same index $m$ on $C'$. This means that $Z'$ is singular
by \ref{cor-cyclic-quo-new}, a contradiction. Thus $m'<m$ as claimed.

Since $C'$ is smooth, $m'\in \ord C^\sharp$. This is not
possible in cases \type{IB} and \type{IC} (because modulo
renumbering of $a_i$'s $a_4=m$, $a_1+a_2\ge m$, $\gcd (a_i,m')=1$,
and $a_i>1$ for $i=1,2,3$). In the case \type{IIB} the only
possibility is $m'=2$ and then the topological index of $f$ is $2$.
\end{proof}

\begin{propositionm}{}
\label{prop-IC}
Assume that $P\in (X,C)$ is a type \type{IC} point. Let $m$ be its
index. Then $(X,C)$ has no other singular points. Moreover,
$i_P(1)=2$, $w_P(0)=1-1/m$, $a_1=2$, and $a_4=m+1$.
\end{propositionm}

\begin{proof}
Assume that $(X,C)$ has one more singular point $Q$ of index $n\ge
1$. By Lemma \ref{lemma-IBICIIB} the base surface $Z$ is smooth.
Since $i_P(1)\ge 2$ \cite[Prop. 5.5]{Mori-1988}, we have $i_P(1)=2$
and $i_Q(1)=1$. We may assume that $Q$ is ordinary of type \type{IA}
or \type{III} by $L$-deformation. Further, by Corollary \ref{cor-wPstar}
$w_P^*(1)+w_Q^*(1)\le 1$.

If $w_P(0)\neq 1-1/m$, all the arguments of \cite[6.5.2]{Mori-1988}
can be applied and we derive a contradiction. Assume that $w_P(0)=
1-1/m$. We follow the arguments of \cite[6.5.3]{Mori-1988}. Since
$P$ is of type \type{IC}, $m\ge 5$, so
$w_Q(0)=1-w_P(0)-1/nm=1/m-1/nm<1/5$. Let $Q$ be of type \type{IA}
(resp. \type{III}). Then, for $1\le d\le 4$, 
by \cite[5.1]{Mori-1988} (resp.
\cite[4.9]{Mori-1988}), one has $w_Q^*(d)=d(d+1)/2$ (resp.
$w_Q^*(d)=\down{(d+1)^2/4}$) for $d\le 4$. On the other hand, by
\cite[5.5]{Mori-1988} $w_P^*(1)=1-\delta_{m,5}$. Therefore, $m=5$.
Further, by \cite[5.5 (v)]{Mori-1988} $w_P^*(2)=0$ and $w_P^*(3)=4$.
Thus,
\[ 
\begin{array}{l}
\sum_{d=1}^3 (1+d+\deg \gr_C^d\omega)= \sum_{d=1}^3\frac{d(d+1)}2
-4-\sum_{d=1}^3\frac{d(d+1)}2=-4<0
\\[5pt]
\qquad
\left(\mbox{resp.}\quad \sum_{d=1}^3\frac{d(d+1)}2 -4- \sum_{d=1}^3
\down{\frac{(d+1)^2}4}=-1<0\right).
\end{array}
\]
Therefore, $H^1(\OOO_X/I_C^{(4)}\totimes \omega_X)\neq 0$. Let
$V:=\Spec\OOO_X/I_C^{(4)}$. Then by Theorem \ref
{theo-sect-1-gen-fiber-m} $V\supset f^{-1}(o)$. Moreover, $V\neq
f^{-1}(o)$ because $V$ is not a local complete intersection inside
$X$. Hence,
\[
2=-K_X\cdot f^{-1}(o)< -K_X\cdot V\le -10K_X\cdot C =10/nm.
\]
This gives $m\le nm< 5$, a contradiction. Hence $P$ is the only
singular point.

By \ref{lemma-KC} and \eqref{eq-grw-w} $w_P(0)=1-1/m$. Hence,
$a_4\equiv 1\mod m$ (\cite[Th. 4.9]{Mori-1988}). Now assume that
$i_P(1)=3$. By \cite[Prop. 5.5]{Mori-1988} $w_P^*(1)\ge 2$. This
contradicts \ref{cor-wPstar}. Thus, $a_1=i_P(1)=2$ (\cite[Prop.
5.5]{Mori-1988}) and $m+1\in \ord C^\sharp$. By normalizedness
$a_4=m+1$ (see \ref{not-nazalo-loc}).
\end{proof}

\begin{propositionm}{}
\label{prop-IB}
$(X,C)$ has no type \type{IB} points.
\end{propositionm}
\begin{proof}
By \cite[Prop. 4.7]{Mori-1988} we can deform $(X,C)$ to
$(X_\lambda,C_\lambda\simeq\PP^1)$, where $X_\lambda$ has at least 
two
non-Gorenstein points of the same index $m$. If
$(X_\lambda,C_\lambda)$ is an extremal neighborhood, the assertion
follows as in the proof of \cite[Th. 6.3]{Mori-1988}. Otherwise
$(X_\lambda,C_\lambda)$ is a $\QQ$-conic bundle germ over a singular
base $(Z_\lambda,o_\lambda)$ by
\ref{cor-cyclic-quo-new} (iii). But $Z$ is smooth by Lemma
\xref{lemma-IBICIIB} and $Z_\lambda$ 
is a deformation
of $Z$, a contradiction. 
\end{proof}

\begin{propositionm}{}
\label{prop-IIB}
If $(X,C)$ has a point $P$ of type \type{IIB}, then $P$ is the only
singular point and the base surface is smooth.
\end{propositionm}
\begin{proof}
Assume that $(X,C)$ has a singular point $Q\neq P$ of index $m\ge
1$. By \cite[Prop. 4.7]{Mori-1988} we can deform $(X,C)$ to
$(X_\lambda,C_\lambda\simeq\PP^1)$, where $X_\lambda$ has three
singular points $P_{\lambda}$, $P'_{\lambda}$ and $Q_{\lambda}$ of
indices $2$, $4$ and $m$. If $(X_\lambda,C_\lambda)$ is an extremal
neighborhood, the assertion follows by \cite[Th. 6.2]{Mori-1988}.
Assume that $(X_\lambda,C_\lambda)$ is a $\QQ$-conic bundle germ
over $(Z_\lambda,o_\lambda)$. Since the indices of $P_{\lambda}$ and
$P'_{\lambda}$ are not coprime, the base surface
$(Z_\lambda,o_\lambda)$ is singular
by \ref{cor-cyclic-quo-new} (iii). So is the base surface $(Z,o)$
of $(X,C)$. By Lemma \ref{lemma-IBICIIB} $(Z,o)$ is of type $A_1$.
This implies that $m$ is even and $\Clsc X\simeq \ZZ\oplus \ZZ_2$
(see Corollary \ref{cor-cyclic-quo}). Then $(X_\lambda,C_\lambda)$
contains three non-Gorenstein points of even indices. But in this
case, the map $\varsigma$ in \eqref{eq-ex-seq-Pic_Cl} cannot be
surjective, a contradiction.
\end{proof}

\begin{mtparag}{Proposition (cf. {\cite[Th. 6.6]{Mori-1988}}).}
\label{prop-IA-size2}
Let $P\in X$ be a type \type{IA} point of index $m$. Then 
$\size_P=1$. If moreover $P$ is the only non-Gorenstein point on
$X$, then $(Z,o)$ is smooth, $w_P(0)=1-1/m$ and $a_2=1$.
\end{mtparag}

\begin{proof}
Assume that $\size_P\ge 2$ and let $m$ be the index of $P$.

\begin{parag}{}
\label{par-IA-a}
First we consider the case when $(Z,o)$ is smooth. We claim that $P$
is the only singular point of $X$. Let $Q\in X$ be a singular point
of index $n\ge 1$. To derive a contradiction we note that the proof
of \cite[Th. 6.6]{Mori-1988} works whenever
\[
H^1(\omega_X/F^2\omega_X)=H^1(\omega_X/F^3\omega_X)=0.
\]
Assume that one of the above vanishings does not hold.

Let $V:= \Spec \OOO_X/I_C^{(j)}$, \ $j=2$ or $3$. By Theorem
\xref{theo-sect-1-gen-fiber-m} (cf. \ref{imp-IC-fiber}) we have
$f^{-1}(o)\subset V$. Hence,
\[
2=-K_X\cdot f^{-1}(o) < -K_X \cdot V=-6K_X \cdot C.
\]
Taking account of $-K_{X}\cdot C=1/mn$ (see \ref{lemma-KC}) we
obtain $mn<3$. So, $m=2$ and $n=1$. On the other hand, a point of
type \type{IA} and index two has $\ord x=(1,1,1,2)$. Such a point is
of size $1$, a contradiction. Thus $P$ is the only singular point of
$X$. By \ref{lemma-KC} and \eqref{eq-grw-w} $w_P(0)=1-1/m$.
Therefore, $a_2\equiv 1\mod m$. On the other hand, $a_2<m$ by
definition of \type{IA} point. Hence, $a_2=1$, $m-a_1a_2=a_3\in \ord
C^\sharp$, and $\size_P=1$, a contradiction.
\end{parag}

\begin{parag}{}
\label{par-IA-b}
Now we consider the case when $(Z,o)$ is singular. So the
topological index of $f$ is $d>1$. Put $m':=m/d$. By definition of
size we have
\begin{equation}
\label{eq-pri_IA-1-a}
2\le U_P(a_1a_2):=\min \{ k \mid k m'd-a_1a_2\in \ord C^\sharp\}.
\end{equation}
Consider the base change \eqref{eq-diag--base-change}. Note that
$P'=g^{-1}(P)$ is also a point of type \type{IA}
of index $m'$  
having
the same index-one
cover as that of $P$.
At $P'$ we have
\[
U_{P'}(a_1a_2) =\min \{ l \mid l m'-a_1a_2\in \ord C^\sharp\}.
\]
Write $l:=U_{P'}(a_1a_2)$ as 
$l=qd-r$, where $0\le r< d$. Then
\[
qm'd-a_1a_2=lm'-a_1a_2+rm'\in \ord C^\sharp.
\]
(We used the fact that $m'\in \ZZ_{>0} a_1+\ZZ_{>0} a_2$, 
see \ref{prop-prim-types-def}.) It is easy to see now from
\eqref{eq-pri_IA-1-a} that $U_{P'}(a_1a_2)=l>qd-d\ge d
(U_P(a_1a_2)-1)\ge 2$. This contradicts the case \ref{par-IA-a}
above.
\end{parag}

\begin{parag}{}
Finally assume that $P$ is the only non-Gorenstein point on
$X$. Then $(Z,o)$ is smooth by Corollary \ref{cor-pr-top-1}. Hence
by \ref{lemma-KC} and \eqref{eq-grw-w} we have $-K_X\cdot C=1/m$ and
$w_P(0)=1-1/m$. Thus $a_2=1$, see \cite[Th. 4.9]{Mori-1988}.
\end{parag}
\end{proof}

Summarizing the results of this section we obtain

\begin{theoremm}
\label{th-prim-conf-summ}
Let $f\colon (X,C\simeq \PP^1)\to (Z,o)$ be a locally primitive
$\QQ$-conic bundle germ. Assume that $X$ is not Gorenstein. Then the
configuration of singular points is one of the following:
\begin{enumerate}
\item
type \type{IC} point $P$ of size $1$ and index $m$ with $i_P(1)=2$,
$w_P(0)=1-1/m$, $a_1=2$, and $a_4=m+1$;
\item
type \type{IIB} point;
\item
type \type{IA} point $P$ of size $1$ and index $m$ with
$w_P(0)=1-1/m$ and $a_2=1$, and possibly at most two more type
\type{III} points;
\item
type \type{IIA} point $P$, and possibly at most two more type
\type{III} points;
\item
two non-Gorenstein points which are of types \type{IA} or
\type{IIA}, and possibly at most one more type \type{III} point;
\item
three non-Gorenstein singular points and no other
singularities \textup(cf. \xref{lemma-3SingPts}\textup).
\end{enumerate}
\end{theoremm}

\begin{remark}
\label{rem-prim-ge}
The existence of a good member of $|-K_X|$ or $|-2K_X|$ in the cases
(i)-(iv) can be shown as in \cite[7.3]{Mori-1988}. The cases (v) and
(vi) will be studied in the following sections.
\end{remark}

\section{The case of three singular points}
In this section we consider $\QQ$-conic bundles with 
exactly three singular points \textup(cf. \xref{lemma-3SingPts}\textup). 
The main result is the following

\begin{mtparag}{Theorem (cf. {\cite[Th. 6.2]{Mori-1988}}).}
\label{th-main-3points}
Let $(X,C\simeq \PP^1)$ be a $\QQ$-conic bundle germ with three
singular points. Up to permutations the configuration of singular
points is one of the following:
\begin{enumerate}
\item
\type{IA}, \type{III}, \type{III} \textup(cf.
\xref{th-prim-conf-summ} {\rm (iii)}\textup);
\item
\type{IA}, \type{IA}, \type{III}. In this case, the indices are $2$,
odd $\ge 3$, and $1$.
\end{enumerate}
In both cases $(Z,o)$ is smooth.
\end{mtparag}

\begin{mparag}{Notation.}
\label{not-prim-3}
To the end of this section we assume that $f\colon (X,C\simeq
\PP^1)\to (Z,o)$ is a $\QQ$-conic bundle germ with three singular
points $P$, $Q$, and $R$. Let $k$, $m$, $n$ be the indices of $P$,
$Q$, and $R$, respectively.
\end{mparag}

\begin{parag}{}
\label{not-prim-3-a}
By \eqref{eq-grO-iP1-2} $i_P(1)=i_Q(1)=i_R(1)=1$ and by Lemma
\ref{lemma-IA-dual-IA-2} all these points are primitive. By
Propositions \ref{prop-IC}, \ref{prop-IB}, and \ref{prop-IIB} $P$,
$Q$, $R$ are of types \type{IA}, \type{IIA}, or \type{III}. We may
assume that $\gr_C^0\omega\simeq \OOO_C(-1)$ (see
\ref{cor-prop-grw=2-2-points-prim}).
\end{parag}

\begin{lemma}
\label{lemma-3-Gor}
A $\QQ$-conic bundle germ $(X,C\simeq\PP^1)$ cannot have three
Gorenstein singular points.
\end{lemma}
\begin{proof}
Indeed, assume that $(X,C)$ has three Gorenstein singular points
$P_1$, $P_2$, $P_3$. In this case, $(Z,o)$ is smooth and $(X,C)$ has
no other singular points. Applying $L$-deformation we may assume
that $P_i$ are ordinary. Then by \cite[Th. 4.9]{Mori-1988}
$w_{P_i}^*(1)=1$. This contradicts Corollary \ref{cor-wPstar}.
\end{proof}

\begin{tparag}{Lemma (cf. {\cite[0.4.13.3, Th. 6.2]{Mori-1988}}).}
\label{lem-prop-3points}
A $\QQ$-conic bundle germ $(X,C\simeq\PP^1)$ has at most two
non-Gorenstein points.
\end{tparag}
\begin{proof}[Proof \textup(following {\rm \cite[0.4.13.3]{Mori-1988}}\textup)]
Assume that $P$, $Q$, $R\in X$ are singular points of indices $k$,
$m$ and $n >1$. By \ref{zam-imp-Gor-fac} we may assume that the
subindex $>1$ for imprimitive points and hence that $(X,C)$ is
locally primitive
(cf. \ref{parag-new-imp-2points}). 
By $L$-deformation at $P$, $Q$ and $R$, and by
\cite[Th. (6.2)]{Mori-1988} we may assume that 
$P$, $Q$ and $R$ are cyclic quotient singularities. 
Using Van Kampen's theorem it is easy to compute
the fundamental group of $X\setminus\{ P,Q,R\}$:
\[
\pi_1(X\setminus\{ P,Q,R\})= \langle
\sigma_1,\sigma_2,\sigma_3\rangle/ \{ \sigma_1^{k}=\sigma_2^{m}=
\sigma_3^{n}=\sigma_1\sigma_2\sigma_3=1\}.
\]
The target group has a finite quotient group $G$ in which the images of
$\sigma_1$, $\sigma_2$, $\sigma_3$ are exactly of order $k$, $m$ and
$n$, respectively (see, e.g., \cite{Feuer-1971}). The above quotient
defines a finite Galois cover $g\colon (X',C')\to (X,C)$. By taking
Stein factorization we obtain a $\QQ$-conic bundle $f'\colon
(X',C')\to (Z',o')$ with irreducible central fiber $C'$. 
By Corollary~\xref{cor-cyclic-quo} $G$ is cyclic.
This
contradicts Corollary~\xref{cor-pr-top-1}.
\end{proof}

\begin{remark}
One can check also that the arguments of \cite[(6.2.2)]{Mori-1988}
work in this case without any changes.
\end{remark}

\begin{propositionm}
\label{prop-3-IIAIIIIII}
In notation \xref{not-prim-3} $(X,C)$ cannot have three singular
points of types \type{IIA}, \type{III}, and \type{III}.
\end{propositionm}

\begin{proof}
Assume that $X$ contains a type \type{IIA} point $P$ and two type
\type{III} points $Q$ and $R$. The base $(Z,o)$ is smooth
by Corollary~\xref{cor-pr-top-1}. Applying
$L$-deformations at $Q$ and $R$ and $L'$ deformation at $P$ (see
\cite[4.12.2]{Mori-1988}) we may assume that $Q$, $R$ are ordinary
and $(X,P)\simeq \{y_1y_2+y_3^2+y_4^3=0\} /\muu_4(1,1,3,2)$, where
$C^\sharp$ is the $y_1$-axis. This new $(X,C)$ is again a
$\QQ$-conic bundle germ by \cite[6.2]{Mori-1988}.

We claim that $H^1(\omega_X\totimes \OOO_X/I_C^{(3)})=0$. Indeed,
otherwise we can apply Theorem \xref{theo-sect-1-gen-fiber-m} 
to $V:=\Spec_X \OOO_X/I_C^{(3)}$ (cf.
\ref{imp-IC-fiber}):
\[
2=-K_{X}\cdot f^{-1}(o) < -K_{X} \cdot V =-6K_{X} \cdot C.
\]
Taking account of $-K_{X}\cdot C=1/4$ (see \ref{lemma-KC}) we obtain
a contradiction. Therefore, $H^1(\omega_X\totimes
\OOO_X/I_C^{(3)})=0$. This implies
\begin{equation}
\label{eq-grw-II-III-III}
\deg \gr_C^1\omega +2+\deg \gr_C^2\omega +3\ge 0
\end{equation}
(see \eqref{eq-grw-1-2}).

By \cite[4.9]{Mori-1988} $w_Q^*(1)=w_R^*(1)=1$ and
$w_Q^*(2)=w_R^*(2)=2$. By Lemma \ref{lemma-comput-w-cAx/4} below and
using \eqref{eq-grw-1} we obtain
\[
\deg \gr_C^1\omega=-2,\qquad \deg \gr_C^2\omega =-4.
\]
This contradicts \eqref{eq-grw-II-III-III}.
\end{proof}

\begin{lemma}
\label{lemma-comput-w-cAx/4}
Let $(X,P)$ be a $cAx/4$-singularity of the form
$\{y_1y_2+y_3^2+y_4^3=0\} /\muu_4(1,1,3,2)$ and let
$C=(\mbox{$y_1$-axis})/\muu_4$. Then
\[
i_P(1)=1,\ w_P(1)=2,\ w_P^*(1)=-1,\ i_P(2)=1,\ w_P(2)=3,\
w_P^*(2)=0.
\]
\end{lemma}
\begin{proof}
Let $I_{C^\sharp}$ be the ideal of ${C^\sharp}$ in $X^\sharp$. Since
$I_{C^\sharp}=(y_2,y_3,y_4)$, we have $I_{C^\sharp}^{(2)}=(y_2)+
I_{C^\sharp}^2$ and
$I_{C^\sharp}^{(3)}=y_2I_{C^\sharp}+I_{C^\sharp}^3$. Let $\bar
\omega$ be a semi-invariant generator of $\omega_{X^\sharp}$. For
example we can take
\[
\bar \omega= \frac{dy_1 \wedge dy_3 \wedge d y_4} {y_1}.
\]
Obviously, $\wt\bar \omega\equiv 1\mod 4$. By definitions of
$\gr_C^i\OOO$ and $\gr_C^i\omega$ we get
\begin{eqnarray*}
\gr_C^1\OOO&=& y_3y_1\cdot \OOO_C \oplus y_4y_1^2\cdot \OOO_C,
\\[7pt]
S^2\gr_C^1\OOO&=& y_3^2y_1^2\cdot \OOO_C \oplus y_3y_4 y_1^3\cdot
\OOO_C\oplus y_4^2y_1^4\cdot \OOO_C,
\\[7pt]
\gr_C^2\OOO&=& y_2y_1^3 \cdot \OOO_C \oplus y_3y_4y_1^3\cdot
\OOO_C\oplus y_4^2\cdot \OOO_C,
\\[7pt]
\gr_C^0\omega&=& y_1^3\bar \omega\cdot \OOO_C,
\\[7pt]
\gr_C^2\omega&=& y_2y_1^2\bar \omega\cdot \OOO_C \oplus
y_3y_4y_1^2\bar \omega\cdot \OOO_C\oplus y_4^2y_1^3\bar \omega\cdot
\OOO_C,
\end{eqnarray*}
where we note that $y_3^2=-y_1y_2$ in $\gr_C^2\OOO$. By definitions
of $i_P$ and $w_P$ (see \ref{not-iPwP}), $i_P(2)=\len \Coker
\alpha_2=1$ and $w_P(2)=\len \Coker \beta_2=3$. Hence, $w_P^*(2)=0$.
Computations for $i_P(1)$, $w_P(1)$ and $w_P^*(1)$ are similar
(for $i_P(1)$, see \cite[2.16]{Mori-1988}).
\end{proof}

\begin{propositionm}
\label{prop-3-points-2-odd}
Let $(X,C\simeq \PP^1)$ be a $\QQ$-conic bundle germ having two
non-Gorenstein points $P$ and $Q$ and one Gorenstein singular point
$R$. Then the indices of non-Gorenstein points are $2$ and odd $\ge
3$. In particular, both $P$ and $Q$ are of type \type{IA}
\textup(cf. \xref{not-prim-3-a}\textup).
\end{propositionm}

\begin{proof}
We use notation of \ref{not-prim-3}. Assume $k\ge 2$, $m\ge 2$,
$n=1$ by hypothesis. Up to permutation we also may assume that $k\le
m$. Apply $L$-deformation so that $P$, $Q$, $R$ become ordinary. In
particular, $P$ and $Q$ are of type \type{IA} points. If this new
$(X,C)$ is an extremal neighborhood, the fact follows by \cite[Th.
6.2]{Mori-1988}. Thus we may assume that $(X,C)$ is a $\QQ$-conic
bundle germ. By \ref{lemma-KC} and \ref{cor-cyclic-quo-new} (iii) we
have $-K_X\cdot C=d/km$ and $-K_{X'}\cdot C'=d^2/km$, where
$d=\gcd(k,m)$ and $(X',C')$ is as in \eqref{eq-diag--base-change}.
If $k=m$, then $X'$ is Gorenstein and $X$ cannot have three singular
points by \ref{prop-Gor-quot}. Thus, $d\le k<m$. If
$H^1(\omega_{X'}\totimes \OOO_{X'} /I_{C'}^{(2)})\neq 0$, then as in
\ref{imp-IC-fiber} by Theorem \xref{theo-sect-1-gen-fiber-m} we have
$f'^{-1}(o')\subset \Spec \OOO_{X'} /I_{C'}^{(2)}$, and so
\[
2=-K_{X'}\cdot f'^{-1}(o') < -3K_{X'} \cdot C'.
\]
We get $3d^2>2km$ and $d=k=m$, a contradiction. Therefore,
\begin{equation}
\label{eq-H1H1-3ppints-IAIAIII}
H^1(\omega_{X'}\totimes \OOO_{X'} /I_{C'}^{(2)})= 0, \qquad
H^1(\omega_X\totimes \OOO_X/ I_C^{(2)})= 0.
\end{equation}
If $H^1(\gr_C^2\omega)= 0$, the arguments of \cite[6.2.3]{Mori-1988}
apply and we are done. So we assume $H^1(\omega_X\totimes \OOO_X/
I_C^{(3)})\neq 0$ and
\[
H^1(\omega_{X'}\totimes \OOO_{X'} /I_{C'}^{(3)})\neq 0.
\]
Again apply Theorem \xref{theo-sect-1-gen-fiber-m} to $(X',C')$ with
$I_{C'}^{(3)}$. We obtain $4\le km<3d^2$. Thus $d>1$ and $X\neq X'$.
Note that in diagram \eqref{eq-diag--base-change} the preimage
$g^{-1}(R)$ consists of $d$ Gorenstein points. By Lemma
\ref{lemma-3-Gor} $d=2$. Hence, $k=2$ and $m=4$. Clearly,
$w_P(0)=1/2$. By Lemma \ref{lemma-KC} and \eqref{eq-grw-w}
$w_Q(0)=1-1/2-1/4=1/4$. Therefore, near $Q$ we have $\ord x_2=3$, so
$\ord x=(1,3,3)$. Further, by \cite[5.1 (ii), 4.9 (ii),
5.3]{Mori-1988}, $w_Q^*(2)=3$, $w_R^*(2)=2$, and $w_P^*(2)=0$. By
\eqref{eq-grw-1} $\deg \gr_C^2\omega=-5$. In particular,
$h^1(\gr_C^2\omega)\ge 2$ and $h^1(\gr_{C'}^2\omega')\ge 2$.

Now we claim $H^0(X',\gr_{C'}^1 \omega')=0$. Note that $(X',C')$ has
three singular points: $Q'$ of index $2$ and two \type{III} points
$R'$, $R''$. By \cite[4.9, 5.3]{Mori-1988} we have $w_{P'}^*(1)=-1$
and $w_{R'}^*(1)=w_{R''}^*(1)=1$. Thus, $\deg \gr_{C'}^1\omega' =
-2$ (see \eqref{eq-grw-1}). Since by \eqref{eq-H1H1-3ppints-IAIAIII}
$H^1(\gr_{C'}^1\omega')=0$, we have $\gr_{C'}^1\omega'\simeq
\OOO_{C'}(-1)\oplus \OOO_{C'}(-1)$ and $H^0(\gr_{C'}^1\omega')=0$.

Finally we apply Proposition \ref{prop-931} below for $(X',C')$ with
$a=2$ and derive a contradiction.
\end{proof}

\begin{proposition}
\label{prop-931}
Let $(X,C)$ be a $\QQ$-conic bundle germ with smooth $(Z,o)$
\textup($C$ may be reducible\textup). Assume that there exists a
positive integer $a$ such that $H^j(X,\gr_C^i\omega)=0$ for all
$j$ and all $i<a$ 
and such that $H^1(X,\gr_C^a\omega)\not=0$. Then
\begin{equation}
\label{eq-3-points-cas-owI}
H^j(X,\omega_X \totimes I_C^{(i)}) =
\left\{
\begin{array}{ll}
0 &(j=0)
\\[5pt]
\omega_Z &(j=1, \ i\le a),
\end{array}
\right.
\end{equation}
and $H^1(X, \gr_C^a \omega)\simeq\CC$.
\end{proposition}
\begin{proof}
We note that the first assertion follows from $H^0(X,\omega_X)=0$
when $j=0$ and from Lemma
\ref{lemma-omega-main} when $j = 1$ and $i=0$.
Consider the natural exact sequence
\[
0 \to \omega_X \totimes I_C^{(i+1)} \to \omega_X \totimes I_C^{(i)}
\to \gr_C^i \omega \to 0.
\]
If \eqref{eq-3-points-cas-owI} is proved for an $i$ ($< a$), we have
$H^1(X,\omega_X \totimes I_C^{(i+1)}) \simeq H^1(X,\omega_X \totimes
I_C^{(i)})$, which proves the assertion for $i+1$. If we set $i=a$,
we have a surjection $\OOO_Z \simeq \omega_Z \to H^1(X,\gr_C^a
\omega)$ which kills $\mathfrak{m}_{Z,o}$. Thus 
$H^1(X,\gr_C^a \omega)\simeq\CC$.
\end{proof}

\begin{proof}[Proof of Theorem {\rm\ref{th-main-3points}}]
By Lemmas \ref{lemma-3-Gor} and \ref{lem-prop-3points} $X$ has one or
two non-Gorenstein points and by \ref{not-prim-3-a} these points are
of types \type{IA}, \type{IIA}, or \type{III}. The case
\type{IIA}+\type{III}+\type{III} is disproved in Proposition
\ref{prop-3-IIAIIIIII} and the cases \type{IA}+\type{IIA}+\type{III}
and \type{IIA}+\type{IIA}+\type{III} are disproved in Proposition
\ref{prop-3-points-2-odd}.
\end{proof}
Finally we note that the existence of a good member of $|-K_X|$ or
$|-2K_X|$ in the cases (i) and (ii) of Theorem \ref{th-main-3points}
can be shown as in \cite[7.3]{Mori-1988}.

\section{Two non-Gorenstein points case:
general (bi)elephants} 
In this section we consider $\QQ$-conic
bundles with two non-Gorenstein points and no other singularities. The
main result of this section is Theorem \ref{th-semist-ge}.

\begin{mparag}{Notation.}
\label{not-2-points-m}
Let $(X,C\simeq \PP^1)$ be a $\QQ$-conic bundle germ having two
singular points $P$, $P'$ of indices $m,\, m'\ge 2$. We assume that
$(X,C)$ is not toroidal because in the toroidal case the existence
of a good divisor in $|-K_X|$ is an easy exercise (see (i) of
\ref{th-semist-ge}). Since $(X,C)$ has at most one imprimitive
point, we will assume that $P'$ is primitive. Let $s$ and $\bar m$
be the splitting degree and the subindex of $P$. Recall that
$m=s\bar m$, that
$m \ge 4$ by \ref{zam-imp-Gor-fac}
if $P$ is imprimitive,
and that $s=1$ if $P$ is primitive. Let $I_C$ be the
sheaf of ideals defining $C$ in $\OOO_X$. Let $\pi^{\sharp}\colon
(X^\sharp,P^\sharp)\to (X,P)$ (resp. $\pi^{\flat}\colon
(X^\flat,P'^\flat)\to (X,P')$) be the index-one cover and $C^\sharp
=\pi^{\sharp -1}(C)_{\red}$ (resp. $C^\flat =\pi^{\flat
-1}(C)_{\red}$). Let $I_C^\sharp$ (resp. $I_C^\flat$) be the
canonical lifting of $I_C$ at $P$ (resp. $P'$). Take normalized
$\ell$-coordinates $(x_1,\dots,x_4)$ (resp. $(x_1',\dots,x_4')$) at
$P$ (resp. $P'$) such that $a_i=\ord x_i$ (resp. $a_i'=\ord x_i'$).

By \ref{cor-prop-grw=2-2-points-prim} and \ref{cor-impr-grw-21cases} we
have $\gr_C^0\omega\simeq \OOO_C(-1)$.
We note that
if $m=2$ then $P$ is primitive
as seen above 
and we can reduce
$m=2$ to the case $m'=2$ by 
switching $P$ and  $P'$.

Thus we distinguish the following three cases:
\begin{parag}{}
\label{case-prim-2odd-two-points}
$m'=2$ and $m$ is odd,
\end{parag}

\begin{parag}{}
\label{case-prim-ge3-two-points}
$m,\, m'\ge 3$, and
\end{parag}

\begin{parag}{}
\label{case-prim-22n-two-points}
$m'=2$, $m=2n$, $n\in \ZZ_{>0}$.
\end{parag}
\par \medskip
\end{mparag}

\begin{mparag}{}
The case \ref{case-prim-2odd-two-points} is easy. Indeed, by
\ref{zam-imp-Gor-fac} and Corollary \ref{cor-imp-IA-size} $(X,C)$ is
primitive. Hence the base $(Z,o)$ is smooth by 
Corollary \ref{cor-pr-top-1} 
and both non-Gorenstein
points are of type \type{IA} by Theorem \ref{th-prim-conf-summ} (v).
We get the case \ref{item=main--th2-IAIA}. The existence of a good
member in $|-2K_X|$ (Proposition \ref{propo-main-ge}) can be shown
exactly as in \cite[7.3]{Mori-1988}. From now on we consider cases
\ref{case-prim-ge3-two-points} and \ref{case-prim-22n-two-points}.
\end{mparag}

\begin{mparag}{}
First, we treat 
the case \xref{case-prim-22n-two-points} till the end of \ref{prop-22n-MMMgrw}.
\begin{lemma}
In the case \xref{case-prim-22n-two-points}, $(X,C)$ is locally
primitive, $n$ is even, and we have
\[
a_1=1,\quad a_2=n+1,\quad a_3=2n-1,\quad a_4=2n.
\]
In particular, $(X,P)$ is of type \type{IA}.
\end{lemma}

\begin{proof}
By \ref{cor-imp-IA-size} $P$ is primitive. 
Hence by \ref{lemma-KC} and \ref{cor-cyclic-quo-new}, (iii)
we have $-K_X\cdot C=\frac 1{2n}$. Further, 
\ref{cor-prop-grw=2-2-points-prim} implies
$\gr_C^0\omega \simeq \OOO_C(-1)$.
Thus by \eqref{eq-grw-w} we have
\[
w_P(0)=1-w_{P'}(0)+K_X\cdot C=\frac12 -\frac{1}{2n}.
\]
Hence, by Proposition \ref{th-prim-conf-summ} and \cite[4.9
(i)]{Mori-1988} the point $P$ is of type \type{IA}. In this case,
$w_P(0)=1-a_2/2n$ (see \cite[4.9 (i)]{Mori-1988}). This gives us
$a_2=n+1$. Since $\gcd(2n,a_2)=1$, $n$ is even.
Finally, by Proposition \ref{prop-IA-size2} $\size_P=1$. Therefore,
$a_1a_2< 2n$ and $a_1=1$. The rest is obvious.
\end{proof}

\begin{tparag}{Lemma (cf. \cite[2.13.1]{Kollar-Mori-1992}).}
\label{lemma--prim-22n-two-points-PPprime}
In the case \xref{case-prim-22n-two-points}, we can write
\begin{eqnarray*}
(X,P)&=&(y_1,y_2,y_3,y_4;\phi)/\muu_m(1,a_2,-1,0;0)\supset
(C,P)=y_1\mbox{-axis}/\muu_m,
\\[5pt]
(X,P')&=&(y_1',y_2',y_3',y_4';\phi')/\muu_2(1,1,1,0;0)\supset
(C,P')=y_1'\mbox{-axis}/\muu_2,
\end{eqnarray*}
where $\phi\equiv y_1y_3 \mod (y_2,y_3)^2+(y_4)$.
\end{tparag}
\begin{proof}
We only need to prove the last equality, which follows from the fact
that $(X,P)$ is a point of type $cA/m$.
\end{proof}

Denote $\ell(P):=\len_P I^{\sharp (2)}/I^{\sharp 2}$, where
$I^\sharp$ is the ideal defining $C^\sharp$ in $X^\sharp$ near 
$P^\sharp$.

\begin{tparag}{Lemma (cf. \cite[2.13.2]{Kollar-Mori-1992}).}
In the case \xref{case-prim-22n-two-points} we have $\ell(P)\le 1$
and $i_P(1)=1$.
\end{tparag}
\begin{proof}
Follows by \ref{lemma--prim-22n-two-points-PPprime} and
\cite[2.16]{Mori-1988}.
\end{proof}

\begin{tparag}{Lemma (cf. \cite[2.13.3]{Kollar-Mori-1992}).}
In the case \xref{case-prim-22n-two-points} we have $\ell(P')\le 1$
and $i_{P'}(1)=1$.
\end{tparag}
\begin{proof}
Assume that $r:=\ell(P')\ge 2$. Then by \cite[2.16]{Mori-1988} the
equation of $X^\flat$ near $P'^\flat$ has the following form:
$\phi'\equiv y'^r_1y'_i \mod (y'_2,y'_3,y'_4)^2$, where $i=3$ (resp. $4$)
if $r$ is odd (resp. even). Consider the following $L$-deformation
$\phi_\lambda'=\phi'+\lambda y'^{r-2}_1y'_i$. Then $(X_\lambda,
C_\lambda)$ has three singular points of indices $m=2n$, $2$, and
$1$. This is impossible by \cite[6.2]{Mori-1988} and
\ref{th-main-3points}. Therefore, $r\le 1$. The last statement
follows by \cite[2.16 (ii)]{Mori-1988}.
\end{proof}
\begin{tparag}{Corollary.}
$\gr_C^1\OOO\simeq\OOO\oplus \OOO(-1)$ in the case
\xref{case-prim-22n-two-points}.
\end{tparag}

\begin{proof}
Follows by \eqref{eq-grO-iP1} because $H^1(\gr_C^1\OOO)=0$.
\end{proof}

\begin{parag}{}
Let $\LLL\subset \gr_C^1\OOO$ be a (unique) subsheaf such that
$\LLL\simeq \OOO$. Note that $\LLL$ is an $\ell$-invertible
$\OOO_C$-module. Let $u_1$ (resp. $u_1'$) be an $\ell$-free
$\ell$-basis at $P$ (resp. $P'$). By \cite[Cor. 9.1.7]{Mori-1988}
there is a subbundle $\MMM\simeq \OOO(-1)$ of $\gr_C^1\OOO$ such
that $\gr_C^1\OOO =\LLL\oplus \MMM$ is an $\ell$-splitting. Let
$u_2$ (resp. $u_2'$) be an $\ell$-free $\ell$-basis of $\MMM$ at $P$
(resp. $P'$).
\end{parag}

\begin{tparag}{Lemma (cf. \cite[2.13.8]{Kollar-Mori-1992}).}
\label{lemma-22n-2-points-qldegPprime}
$\qldeg (\MMM,P')=1$.
\end{tparag}
\begin{proof}
Since $\qldeg (\MMM,P')< m'=2$, it is sufficient only to disprove
the case $\qldeg (\MMM,P')=0$. Assume that $\qldeg (\MMM,P')=0$.
Since $y_2$, $y_3$, $y_4$ form an $\ell$-basis of $\gr_C^1\OOO$ at
$P$, we have $\MMM\simeq (-1+iP^\sharp)$, where $i=0$, $1$, or
$m-a_2\ (=n-1)$. Recall that $\gr_C^0\omega\simeq
(-1+(m-a_2)P^\sharp+P'^\flat)$. Therefore,
\[
\gr_C^1\omega\simeq \gr_C^1\OOO\totimes \gr_C^0\omega \simeq
\LLL\totimes \gr_C^0\omega \toplus (-2+(m-a_2+i)P^\sharp+P'^\flat).
\]
The last expression is normalized (because $m-a_2+i=n-1+i<m=2n$).
Hence, $H^1(\gr_C^1\omega)\neq 0$. Put $V:=\Spec \OOO_X/I_C^{(2)}$
and $V':=\Spec \OOO_{X'}/I_{C'}^{(2)}$ (notation of
\eqref{eq-diag--base-change}). As in the proof of Corollary
\ref{cor-lemma-sect-1-gen-omega-OOO-ne0} we get $H^1(\OOO_V\totimes
\omega_X)\neq 0$ and therefore $H^1(\OOO_{V'}\totimes
\omega_{X'})\neq 0$ (notation of \eqref{eq-diag--base-change}). By
Theorem \ref{theo-sect-1-gen-fiber-m} $V'\supset f'^{-1}(o')$. In
particular,
\[
2=-K_{X'}\cdot f'^{-1}(o')\le -K_{X'}\cdot V'=-3K_{X'}\cdot C'=3/n.
\]
This implies $n=1$, a contradiction.
\end{proof}

\begin{tparag}{Lemma (cf. \cite[2.13.10]{Kollar-Mori-1992}).}
$\qldeg (\MMM,P)=m-a_2$.
\end{tparag}

\begin{proof}
First we note that the arguments of 
\cite[2.13.10.1-2]{Kollar-Mori-1992}
apply to our case and
show in particular that if $\qldeg (\MMM,P)\neq m-a_2$ and if
$H^1(\omega_X/F^4(\omega,J))=0$, then
$m$ is odd while $m$ is even in our case.
Hence it is enough
to derive a contradiction assuming that
$\qldeg (\MMM,P)\neq m-a_2$ and
$H^1(\omega_X/F^4(\omega,J))\neq 0$.

Let $J$ be the $C$-laminal
ideal of width $2$ such that $J/I_C^{(2)}=\LLL$. Then
\[
0\neq H^1(\omega_X/F^4(\omega,J))=H^1(\omega_X/ J^2\omega_X)=
H^1(\omega_X\totimes \OOO_X/ J^{(2)}).
\]
As in the proof of Lemma 
\ref{lemma-22n-2-points-qldegPprime}
put $V:=\Spec_{X'} \OOO_{X'}/J'^{(2)}$, where
$J'$ is the pull-back of $J$ on $X'$
(we use notation of \eqref{eq-diag--base-change}).
Recall that $I_C\supset J\supset I_C^{2}$.
Thus $I_C'\supset J'\supset I_C'^{2}$, where $I_C'$
is the ideal sheaf of $C'$.
Since $H^1(\omega_X\totimes \OOO_X/ J^{(2)})\neq 0$,
we have $H^1(\omega_{X'}\totimes \OOO_{X'}/ J'^{(2)})\neq 0$.
By Theorem 
\xref{theo-sect-1-gen-fiber-m} 
$V\supset f'^{-1}(o')$. Let $Q'\in C'$ be a
general point. Then in a suitable coordinate system $(x,y,z)$ near
$Q'$ we may assume that $C'$ is the $z$-axis. So, $I_C'=(x,y)$ and
$I_C'^{(2)}=(x^2,xy,y^2)$. Since $J'/I_{C'}^{(2)}$ is of rank $1$,
by changing coordinates $x,y$ we may assume that $J'=(x,y^2)$
near $Q'$. Then $J'^{(2)}=(x^2, xy^2, y^4)$ and $V\equiv l C'$,
where $l=\len _0 \CC[x,y]/(x^2, xy^2, y^4)=6$.
Similarly 
to the proof of Lemma \ref{lemma-22n-2-points-qldegPprime} we have
\[
2=-K_{X'}\cdot f'^{-1}(o') \le -K_{X'} \cdot V =-6K_{X'} \cdot C'=6/n.
\]
Since $n$ is even $\ge 2$, we get only one possibility $n=2$.

As in the proof of Lemma \ref{lemma-22n-2-points-qldegPprime} we see
that $\qldeg(M,P)=0$ (because $m-a_2=1$).
Then by \ref{lemma-22n-2-points-qldegPprime} $\MMM\simeq (-1+P'^\flat)$.
Now we consider the base change \eqref{eq-diag--base-change}.
Here $g$ is a double cover and $X'$ is of index two.
Set $P^{\natural}:=g^{-1}(P)$,
the unique non-Gorenstein point of $X'$.
(Note that in our notation $P'\in X$, which is
different from the notation of \ref{base-change}--\ref{prop-cyclic-quo})).
Note that the index-one
covers of $(X,P)$
and $(X',P^{\natural})$ coincide.
Let $\MMM'$ be the $\ell$-invertible sheaf
on $X'$, the pull-back of $\MMM$
to $X'$. We have
$\MMM'\simeq (-1+0P^\sharp)$ and $\gr_{C'}^0
\omega_{X'} 
\simeq (-1+P^\sharp)$,
where $C':=g^{-1}(C)_{\red}$.
Hence, $H^1(X',\omega_{X'}\totimes \MMM') \neq 0$.
Let $J'$ be the ideal on $X'$ lifting $J$. Then
taking account of the exact sequence
\[
0 \longrightarrow
I_{C'}/J' (= \MMM')
\longrightarrow \OOO_{X'}/J'
\longrightarrow \OOO_{C'} \longrightarrow 0
\]
and isomorphisms
$\omega_{X'} \totimes\OOO_{C'}
\simeq \gr_{C'}^0\omega\simeq \OOO_{C'}(-1)$,
we get
$H^1(X', \omega_{X'} \totimes \OOO_{X'}/J')\simeq
H^1(X',\omega_{X'}\totimes \MMM') \neq 0$.
Hence by Theorem \xref{theo-sect-1-gen-fiber-m} we have
$f'^{-1}(o') \subset \Spec \OOO_{X'}/J'$
which means
$2 \le 2/2=1$, a contradiction.
\end{proof}

\begin{remark}
In the above notation the case $n=2$ can be disproved also by 
considering possible actions of involutions on index two 
$\QQ$-conic bundles $f'\colon X'\to Z'$ (see \ref{cla-index-4a} and
\ref{cla-index-4b}).
\end{remark}
Thus we have proved the following
\begin{mtparag}{Proposition.}
\label{prop-22n-MMMgrw}
In the case \xref{case-prim-22n-two-points} there is an
$\ell$-isomorphism $\MMM \simeq \gr_C^0\omega$.
\end{mtparag}
\end{mparag}

\begin{lemmam}
\label{lemma-mm13-types}
Up to permutations we may assume that $P'$ is of type \type{IA} and
$P$ is of type \type{IA^\vee}, \type{IA}, or \type{IIA}. Moreover,
$\size_P=\size_{P'}=1$.
\end{lemmam}
\begin{proof}
If $(X,C)$ is not locally primitive, the assertion follows by
\ref{lemma-IA-dual-IA-2} and \ref{prop-imp-IA-size}. We assume that
$(X,C)$ is locally primitive. By Propositions \ref{prop-IC},
\ref{prop-IB}, and \ref{prop-IIB} points $P$ and $P'$ are of types
\type{IA} or \type{IIA}. 
If both $P$ and $P'$ are of type
\type{IIA}, then $w_P(0)+w_{P'}(0)=3/2>1$ (see \cite[4.9
(i)]{Mori-1988}). This contradicts \eqref{eq-grw-w}. 
Thus we may
assume that $P'$ is of type \type{IA} modulo permutation of $P$ and $P'$. 
To prove the last statement
consider $L$-deformation $(X_\lambda, C_\lambda\ni P_\lambda,\,
P'_\lambda)$ of $(X,C\ni P,\, P')$ so that $P_\lambda,\, P'_\lambda$
are ordinary points. In particular, they are of type \type{IA}. By
\cite[4.7]{Mori-1988} $\size_{P_\lambda}=\size_{P}$ and
$\size_{P'_\lambda}=\size_{P'}$. If $(X_\lambda, C_\lambda)$ is a
$\QQ$-conic bundle germ, the assertion follows by Proposition
\ref{prop-IA-size2}. Otherwise we can apply \cite[Th.
6.6]{Mori-1988}.
\end{proof}

Temporarily we consider the following situation.

\begin{mparag}{Notation.}
\label{not-2-points-IA-ordin-u1u2}
Assume that $P$ and $P'$ are ordinary. Then $i_P(1)=\size_P=1$ and
$i_{P'}(1)=\size_{P'}=1$. Hence, $\gr_C^1\OOO\simeq \OOO\oplus
\OOO(-1)$ (because $H^1(\gr_C^1\OOO)=0$). Let $\LLL\subset
\gr_C^1\OOO$ be a (unique) subsheaf such that $\LLL\simeq \OOO$.
Note that $\LLL$ is an $\ell$-invertible $\OOO_C$-module. Let $u_1$
(resp. $u_1'$) be an $\ell$-free $\ell$-basis at $P$ (resp. $P'$).
\end{mparag}

\begin{tparag}{Theorem (cf. {\cite[9.3]{Mori-1988}}).}
\label{th-mm13-ordinary}
Notation as in \xref{not-2-points-IA-ordin-u1u2}. Then both
$(C^\sharp,P^\sharp)$ and $(C^\flat,P'^\flat)$ are smooth, $\wt
u_1\equiv -1\mod m$, $\wt u_1'\equiv -1\mod m'$ and furthermore we
may assume that $a_1=a_1'=1$.
\end{tparag}

The proof follows \cite[9.3]{Mori-1988}. 
We will treat \ref{th-mm13-ordinary}
in four cases
\ref{case-imprim-2points-imp}--\ref{case-prim-2points-93c} below.

\begin{parag}{Case: $P$ is of type \type{IA^\vee}.}
\label{case-imprim-2points-imp}

\end{parag}
If both $P$ and $P'$ are primitive, then in a suitable coordinate
system near $P^\sharp$ 
the ideal $I_C^\sharp$ is generated by $x_1^{a_2}-x_2^{a_1}$
and $x_3$ (because $P$ is ordinary). 
Hence, either $\wt u_1\equiv
a_1a_2$ or $\wt u_1\equiv -a_1 \mod m$
holds and the corresponding
assertion also holds for $P'$.
Modulo permutation of $P$ and $P'$ and if 
$a_2=1$
(resp. $a'_2=1$) modulo further
permutation of 
$a_1$ and $a_3$
(resp. $a'_1$ and $a'_3$) there are three cases:

\begin{parag}
{Case: $\wt u_1\equiv a_1a_2$, $a_2\not\equiv \pm 1\mod m$,\ $\wt
u_1'\equiv a_1'a_2'$, $a_2'\not\equiv \pm 1\mod m'$.}
\label{case-prim-2points-93a}
\end{parag}
\par
\begin{parag}
{Case: $\wt u_1\equiv a_1a_2$, $a_2\not\equiv \pm 1\mod m$,\ $\wt
u_1'\equiv a_3'\mod m'$.}
\label{case-prim-2points-93b}
\end{parag}
\par
\begin{parag}
{Case: $\wt u_1\equiv a_3\mod m$,\ $\wt u_1'\equiv a_3'\mod m'$.}
\label{case-prim-2points-93c}
\end{parag}
\par\medskip

We will show that only the case \ref{case-prim-2points-93c} is
possible and $a_1=a_1'=1$. By \cite[Cor. 9.1.7]{Mori-1988} there is
a subbundle $\MMM\simeq \OOO(-1)$ of $\gr_C^1\OOO$ such that
$\gr_C^1\OOO =\LLL\oplus \MMM$ is an $\ell$-splitting. Let $u_2$
(resp. $u_2'$) be an $\ell$-free $\ell$-basis of $\MMM$ at $P$
(resp. $P'$).

\begin{parag}{}
Let $J$ be the $C$-laminal ideal of width $2$ such that
$J/I^{(2)}_C=\LLL$,
and our symbols are compatible with those in
\cite[9.3.2]{Mori-1988}.
\end{parag}

Note that $w_P(0)=1-{a_2}/\bar m$ and $w_{P'}(0)=1-a_2'/m'$. 
Then by \eqref{eq-grw-w} and $(K_X \cdot C)<0$ we have
\begin{equation}
\label{eq-multl-pr-main-1-right}
1<\frac{a_2}{\bar m}+\frac{a_2'}{m'} \quad (\cite[9.3.4]{Mori-1988}).
\end{equation}
Using Lemma \ref{lemma-KC} that holds only for
$\QQ$-conic bundles, we have
\begin{equation}
\label{eq-multl-pr-main-1-left}
1+\frac{d}{m m'}=\frac{a_2}{\bar m}+\frac{a_2'}{m'},
\end{equation}
where $d=m m'/\lcm(\bar m, m')=s \gcd(\bar m, m')$ by \ref{cor-cyclic-quo-new}, (iii).

\begin{parag}{Disproof of the case
{\ref{case-imprim-2points-imp}}.}
\label{par-prop-imp-unique}
This case corresponds to \cite[9.3.ipr]{Mori-1988},
and is disproved by the same argument as
\cite[9.3.5]{Mori-1988}.

Hence we note that ${\bar m}=m$ below till the end
of the proof of \ref{th-mm13-ordinary}.
\end{parag}

\begin{parag}{Disproof of the case {\ref{case-prim-2points-93a}}.}
This case corresponds to \cite[9.3.a]{Mori-1988},
and is disproved by the same argument as \cite[9.3.6]{Mori-1988}.
\end{parag}

\begin{parag}{Disproof of the case {\ref{case-prim-2points-93b}}.}
This case corresponds to \cite[9.3.b]{Mori-1988},
and is disproved by the same argument as \cite[9.3.7]{Mori-1988}.
\end{parag}

\begin{parag}{Treatment of the case {\ref{case-prim-2points-93c}}.}
This case corresponds to \cite[9.3c]{Mori-1988}, and
the arguments of \cite[9.3.8]{Mori-1988} work
except for \cite[9.3.8.6]{Mori-1988}, which we prove below
using \eqref{eq-multl-pr-main-1-left}.

We have $\wt u_1\equiv a_3\mod m$, $\wt u_1'\equiv a_3'\mod m'$. We
will prove that $a_1=a_1'=1$. By symmetry we may assume that
$a_2'/m'>1/2$ (see \eqref{eq-multl-pr-main-1-right}). Since
$\size_{P'}=1$, $m'\ge a_1'a_2'$. This gives us $a_1'=1$. We will
prove $a_1=1$. Assume that $a_1\ge 2$. Then computations
\cite[9.3.8.3--9.3.8.4]{Mori-1988} apply and give us $a_1=2$,
$a_2=1$. In particular, $C^\sharp$ is smooth over $P$ and $P'$.
Further by \cite[9.3.8.5]{Mori-1988} 
we have
\begin{equation}
\label{eq-prim-maim-LM}
\begin{array}{ll}
ql_C(\LLL)&=2P^\sharp+P'^\sharp,
\\[5pt]
ql_C(\MMM)&=-1+(m-2)P^\sharp+(m'-a_2')P'^\sharp,
\end{array}
\end{equation}
and $ql_C(\gr_C^0\omega)=-1+(m-1)P^\sharp+(m'-a_2')P'^\sharp$.
\end{parag}

\begin{tparag}{Claim (cf. {\cite[9.3.8.6]{Mori-1988}}).}
$m\ge 5$.
\label{10-1-2-m-ge-5}
\end{tparag}
\begin{proof}
Assume that $m<5$. Since $a_1=2$ and $\gcd (m,a_1)=1$, we have
$m=3$. Then by \eqref{eq-multl-pr-main-1-left}
\[
2m'+d=3a_2',\quad d=\gcd(3,m').
\]
Now one can see that computations of \cite[9.3.8.6]{Mori-1988} apply
and give us $\chi(F^1(\omega,J)/F^4(\omega,J))<0$. From the exact
sequence
\[
0 \longrightarrow F^1(\omega,J)/F^4(\omega,J) \longrightarrow
\omega_X/F^4(\omega,J) \longrightarrow \gr_C^0\omega \longrightarrow
0
\]
we get
\[
\chi(\omega_X/F^4(\omega,J))=
\chi(F^1(\omega,J)/F^4(\omega,J))+\chi(\omega_X/F^1(\omega,J))<0.
\]
(we note that $F^1(\omega,J)= \mt{Sat}_{\omega_X}(I_C\omega_X)$ and
$\omega_X/F^1(\omega,J)=\gr_C^0\omega$). In particular, we have
$H^1(\omega_X/F^4(\omega,J))\neq 0$. Recall that $I_C\supset
J\supset I_C^{(2)}$ and
$F^4(\omega,J)=\mt{Sat}_{\omega_X}(J^2\omega_X)$. Assume that
$(Z,o)$ is smooth, i.e., $3\nmid m'$. Then by Theorem
\ref{theo-sect-1-gen-fiber-m} $f^{-1}(o)\subset \Spec \OOO_X/J^{(2)}
\subset \Spec \OOO_X/I_C^{(4)}$. We get a contradiction (cf.
\ref{imp-IC-fiber}):
\[
2=-K_X\cdot f^{-1}(o) < -10 K_X \cdot C =10/(3m'), \quad m'= 1.
\]
Now assume that $m'=3m''$, $m''\ge 2$. Take a Weil divisor $\xi$
such that $ql_C \xi=P^\sharp -m''P'^\sharp$. Then $\xi$ is a
$3$-torsion in $\Clsc X$. Taking \eqref{eq-prim-maim-LM} into
account we obtain
\[ 
\begin{array}{l}
ql_C(\MMM\totimes \xi)=-1+(m-1)P^\sharp+(2m''-a_2')P'^\sharp =
\\
\hspace{60pt}
=-2+(m-1)P^\sharp+(5m''-a_2')P'^\sharp.
\end{array}
\]
Since $a_2'=2m''+1$, $5m''-a_2'=3m''-1$. So, the last expression is
normalized. By \cite[8.9.1 (iii)]{Mori-1988} $\deg_C \MMM\totimes
\xi=-2$. Note that $\gr_C^1\OOO \totimes \xi= (\LLL \totimes
\xi)\toplus(\MMM\totimes \xi)$. Hence, $H^1(\gr_C^1\OOO \totimes
\xi)\neq 0$. 
This is a contradiction, and $m \ge 5$ as proved.
\end{proof}

The remainder of the proof is the same as 
\cite[9.3.8.7]{Mori-1988}.
Thus Theorem \ref{th-mm13-ordinary} is proved.

Now we treat the case \ref{case-prim-ge3-two-points}
from here till the end  of 
\ref{prop-mm13-MMMgrw}.

\begin{propositionm}
\label{prop-mm13-IIA}
In notation and assumptions of \xref{case-prim-ge3-two-points}
$(X,C)$ has no \type{IIA} type points.
\end{propositionm}
\begin{proof}
Assume that $P$ 
is of type \type{IIA} (and $P'$ is of type \type{IA}
by Lemma 
\ref{lemma-mm13-types}).
Let $d:=\gcd (4,m')$. By \cite[Th. 4.9, (i)]{Mori-1988} 
$w_P(0)=3/4$. Then by Lemma \ref{lemma-KC} and
\eqref{eq-grw-w} we have
\[
w_{P'}(0)=1-w_{P}(0)+K_X\cdot C=1/4-d/(4m')<1/2.
\]
Hence, $a_1'=1$ (see \cite[Prop. 5.1]{Mori-1988}). Moreover,
\begin{equation}
\label{eq-main-CAXCA}
d+3m'=4a_2'
\end{equation}
(because $w_P(0)=1-a_2'/m'$ by 
\cite[Th. 4.9, (i)]{Mori-1988}). Applying $L$-deformations at $P'$ and
$L'$ deformation at $P$ (see \cite[4.12.2]{Mori-1988}) we may assume
that $P'$ is ordinary and $(X,P)\simeq \{y_1y_2+y_3^2+y_4^3=0\}
/\muu_4(1,1,3,2)$, where $C^\sharp$ is the $y_1$-axis. This new
$(X,C)$ is again a $\QQ$-conic bundle germ by \cite[9.4]{Mori-1988}.
Applying \cite[9.4.3-9.4.5]{Mori-1988} we get an $\ell$-splitting
$\gr_C^1\OOO=\LLL \toplus \MMM$, where $\LLL\simeq\OOO$ and
$\MMM\simeq \OOO(-1)$. Moreover, $ql_C(\LLL)=P^\sharp+P'^\sharp$ and
$ql_C(\MMM)=-1+2P^\sharp +(m'-a_2')P'^\sharp$ and we may assume that
$y_3$ is an $\ell$-free $\ell$-basis of $\LLL$ at $P$.

Now one can see that computations of \cite[9.4.6]{Mori-1988} apply
and give us $\chi(F^1(\omega,J)/F^4(\omega,J))<0$.
If $2\not|m'$, then we can apply the second half of the proof of
\ref{10-1-2-m-ge-5} to get a contradiction:
\[
2=-K_X\cdot f^{-1}(o) < -10 K_X \cdot C =10/(4m'), \quad m'=1.
\]
Thus $2 | m'$. 
If $m'=2m''$ and $2 \not|m''$, then considering diagram
\eqref{eq-diag--base-change} one has
$H^1(\omega_{X'}/\mt{Sat}_{\omega_{X'}}(J^2\omega_{X'}))\neq 0$ 
and similarly gets:
\[
2=-K_{X'}\cdot f'^{-1}(o') < -10 K_{X'} \cdot C' = -20 K_{X} \cdot
C= 20 /2m'=5/m''.
\]
Thus $m'' = 1$, and one sees $d=2$ and $a_2'=2$ by \eqref{eq-main-CAXCA},
which contradicts $\gcd (m',a_2')=1$,
a condition on \type{IA} points. 
Hence $4 | m'$ and we set $m'=4m''$.

Take a Weil divisor $\xi$ such that 
$ql_C (\xi)=P^\sharp -m''P'^\sharp$. 
Then $\xi$ is a $4$-torsion in $\Clsc X$. By
\cite[9.4.5]{Mori-1988} we have
\[
ql_C(\MMM)=-1+2P^\sharp+(4m''-a_2')P'^\sharp
\]
Recall that $a_2'=3m''+1$. Taking this into account we obtain
\[
ql_C(\MMM\totimes \xi)=-1+3P^\sharp-P'^\sharp
=-2+3P^\sharp+(4m''-1)P'^\sharp.
\]
The last expression is normalized. By \cite[8.9.1 (iii)]{Mori-1988}
$\deg_C \MMM\totimes \xi=-2$. Note that $\gr_C^1\OOO \totimes \xi=
(\LLL \totimes \xi)\toplus(\MMM\totimes \xi)$. Hence,
$H^1(\gr_C^1\OOO \totimes \xi)\neq 0$. This is a contradiction.
\end{proof}

\begin{tparag}{Corollary (\cite[9.4.7]{Mori-1988}).}
\label{cor-2nonGor-IAIA}
In notation and assumptions of \xref{case-prim-ge3-two-points}
points $P$ and $P'$ are type \type{IA} points such that
$a_1=a_1'=1$, and moreover $\ell(P)<m$ and $\ell(P')<m'$.
\end{tparag}

\begin{proof}
By \ref{lemma-mm13-types} and \ref{prop-mm13-IIA} $P$ is of type
\type{IA} and $P'$ is of type \type{IA} or \type{IA^\vee}. Replacing
$(X,C)$ with $L$-deformation we may assume that both $P$ and $P'$
are ordinary. Then by \ref{par-prop-imp-unique}, \ref{th-mm13-ordinary},
and
\cite[9.3, 9.4]{Mori-1988} $P'$ is of type \type{IA} and $a_1=a_1'=1$
($L$-deformation does not change $a_i$'s because $P$ and $P'$ are of
type \type{IA} or \type{IA^\vee}). If $\ell(P)\ge m$, then an
$L'$-deformation $(X_\lambda,C_\lambda)$ (see
\cite[4.12.2]{Mori-1988}) has at least one Gorenstein singular point
besides $P_\lambda$ and $P_\lambda'$. This contradicts
\ref{th-main-3points} and \cite[6.2]{Mori-1988}. Thus $\ell(P)<m$.
By symmetry we also have $\ell(P')<m'$.
\end{proof}

\begin{tparag}{Corollary (\cite[9.4.8]{Mori-1988}).}
In notation and assumptions of \xref{case-prim-ge3-two-points} we
have $i_P(1)=i_{P'}(1)=1$ and an isomorphism
$\gr_C^1\OOO\simeq\OOO\oplus\OOO(-1)$.
\end{tparag}

\begin{proof}
Since $\ell(P)<m$ and $\ell(P')<m'$, by \cite[2.16 (ii)]{Mori-1988}
we have $i_P(1)=i_{P'}(1)=1$. Hence $\deg \gr_C^1\OOO=-1$, see
\eqref{eq-grO-iP1}. Taking account of $H^1(\gr_C^1\OOO)=0$ we obtain
the last statement.
\end{proof}

\begin{mtparag}{Proposition (\cite[9.8]{Mori-1988}).}
\label{prop-th-prim-2-pts-ell-qldeg}
We have
\[
\ell(P)+ \qldeg (\LLL,P)=\ell(P')+ \qldeg (\LLL,P')=1.
\]
\end{mtparag}

\begin{proof}
Similar to the proof of \cite[Theorem 9.8]{Mori-1988}.
\end{proof}

\begin{mtparag}{Proposition (\cite[9.9.1]{Mori-1988}).}
\label{prop-mm13-MMMgrw}
In the case \xref{case-prim-ge3-two-points},
there is an $\ell$-isomorphism $\MMM \simeq \gr_C^0\omega$.
\end{mtparag}

\begin{proof}
Since $\gr_C^0\omega\simeq\OOO(-1)$, it is sufficient to show that
$\qldeg (\MMM, P) =\qldeg (\gr_C^0\omega, P)=R(a_2)$ and $\qldeg
(\MMM, P') =\qldeg (\gr_C^0\omega, P')=R'(a_2')$. By symmetry it is
sufficient to prove for example the first equality. According to
\ref{prop-th-prim-2-pts-ell-qldeg} there are two cases.

\begin{parag}{Case: $\ell(P)=0$, $\qldeg (\LLL,P)=1$.}
Then $(X,P)$ is a cyclic quotient singularity of type $\frac 1m
(1,a_2,-1)$. If $u_1$ is an $\ell$-free $\ell$-basis of $\LLL$ at
$P$, then $R(\wt u_1)=1$, so $\wt u_1 \equiv -1\mod m$. An
$\ell$-free $\ell$-basis of $\gr_C^1\OOO$ at $P$ is $x_2$, $x_3$.
Hence we can put $u_1=x_3$ and $x_2$ is an $\ell$-free $\ell$-basis
of $\MMM$. Therefore, $\qldeg (\MMM, P)=R(\wt x_2)=R(a_2)$.
\end{parag}

\begin{parag}{Case: $\ell(P)=1$, $\qldeg (\LLL,P)=0$.}
Then we can choose a some coordinate system so that
$(X^\sharp,P^\sharp)$ is given by $\phi=0$ with $\phi\equiv
x_1x_3\mod (x_2,x_3,x_4)^2$ and $C^\sharp$ is the $x_1$-axis (see
\cite[2.16]{Mori-1988}). If $u_1$ is an $\ell$-free $\ell$-basis of
$\LLL$ at $P$, then $R(\wt u_1)=0$, so $\wt u_1 \equiv 0\mod m$.
Again an $\ell$-free $\ell$-basis of $\gr_C^1\OOO$ at $P$ is $x_2$,
$x_4$. Hence we can put $u_1=x_4$ and $x_2$ is an $\ell$-free
$\ell$-basis of $\MMM$. Therefore, $\qldeg (\MMM, P)=R(\wt
x_2)=R(a_2)$.
\end{parag}
\end{proof}

Taking Propositions \ref{prop-22n-MMMgrw} and \ref{prop-mm13-MMMgrw}
into account one can see that all the arguments and computations
from \cite[9.9.2-9.9.10]{Mori-1988} apply in our case. This proves
(iii) of the following theorem (cf. \cite[2.2.4]{Kollar-Mori-1992}).

\begin{mtparag}{Theorem (cf. \cite[Th. 9.10]{Mori-1988},
\cite[2.2.4]{Kollar-Mori-1992}).}
\label{th-semist-ge}
Let $(X,C\simeq \PP^1)$ be a $\QQ$-conic bundle germ having two
non-Gorenstein points $P$, $P'$ of indices $m,\, m'\ge 2$ and no
other singularities. Then $P$, $P'$ are of type \type{IA} by
\xref{lemma-mm13-types} and \xref{prop-mm13-IIA}, and the following
assertions 
hold.
\begin{enumerate}
\item
If $m'=2$ and $m$ is odd, then the general member of $|-2K_X|$ does
not contain $C$ and has only log terminal singularities.
\item
If $(X,C)$ is toroidal, then a general member $F\in |-K_X|$ does not
contain $C$. It has two connected components. Each of them is a Du
Val singularity of $A$-type.
\item
If $(X,C)$ is not toroidal, we further assume either $m,\, m'\ge 3$
or $m'=2$ and $m$ is even. Then a general member 
$F \in |-K_X|$ is a
normal surface containing $C$, smooth outside of $\{P,\, P'\}$, with
Du Val points of $A$-type at $P$, $P'$.
Furthermore $C$ on $F$ is contractible to a Du Val point of
$A$-type.
\end{enumerate}
\end{mtparag}

Note that (i) and (ii) of \ref{th-semist-ge} are easy (cf. \cite[Th.
7.3]{Mori-1988}). In (ii) one can also take as $F$ the sum of two
horizontal toric divisors.

\section{Two non-Gorenstein points case: the classification}
The following is the main result of this section.

\begin{theoremm}
\label{th-semist-ge-conseq}
Let $(X,C\simeq\PP^1)$ be a $\QQ$-conic bundle germ having two
points of indices $m,\, m'\ge 2$
and no other singularities. Assume either $m,\, m'\ge
3$ or $m'=2$ and $m$ is even. Then $(X,C)$ is either toroidal or as
in \xref{item=main--th-pr-ex3}.
\end{theoremm}

The above theorem is an easy consequence of Theorem
\ref{th-semist-ge} and Proposition \ref{prop-prim-2points-ge->1}
below.

\begin{mtparag}{Proposition (cf. \cite[\S 4]{Prokhorov-1997_e}).}
\label{prop-prim-2points-ge->1}
Let $f\colon (X,C)\to (Z,o)$ be a non-Gorenstein $\QQ$-conic bundle
germ with $C\simeq \PP^1$. Assume that the general element $F\in
|-K_X|$ containing $C$ has only Du Val singularities. Let
$F\stackrel{f_1}{\longrightarrow} \bar F\to Z$ be the Stein
factorization and let $\bar P=f_1(C)$. Assume that $(\bar F,\bar P)$
is a singularity of type $A$. Then one of the following holds:
\begin{enumerate}
\item
$f$ is as in \xref{item=main--th-pr-ex3}, or
\item
$X$ is of index $2$ and $(Z,o)$ is smooth.
\end{enumerate}
\end{mtparag}

\begin{proof}[Proof of Proposition \xref{prop-prim-2points-ge->1}]
By the inversion of adjunction \cite[17.6]{Utah} the log divisor
$K_X+F$ is plt. Consider diagram \eqref{eq-diag--base-change} and
put $F':=g^*F$. We may assume that $Z'\simeq \CC^2$ and $Z\simeq
\CC^2/\muu_d(1,q)$, where $\gcd (d,q)=1$. By \cite[20.3]{Utah}
$K_{X'}+F'=g^*(K_X+F)\sim 0$ is plt. In particular, $F'$ is normal
and irreducible. Further, diagram \eqref{eq-diag--base-change}
induces the following diagram
\begin{equation}
\label{eq-diag--base-change-F}
\begin{CD}
(F',C')@>{g_{F'}}>> (F,C)
\\
@V{f'_1}VV @V{f_1}VV
\\
(\bar F',\bar P')@>{\bar g}>> (\bar F,\bar P)
\\
@V{f'_2}VV @V{f_2}VV
\\
(Z',o')@>{h}>> (Z,o)
\end{CD}
\end{equation}
where the vertical arrows are Stein factorizations of restrictions
$f'|_{F'}$ and $f|_F$. It is clear that $f_2'$ and $f_2$ are double
covers. 
By adjunction
$K_{F'}\sim 0$ and $f_1'$ is a
crepant morphism contracting $C'$. Since $\bar g$ is \'etale in
codimension one
and $(F,P)$ is a singularity of type $A$,
$(\bar{F}',\bar{P}')$ is also of type $A$.
Note that $(\bar F',\bar P')$ cannot be smooth (because $f_1'$ is
non-trivial).

Consider the case $d\ge 2$. Then by Lemma
\ref{lemma-vspom-cover-DuVal} $(\bar F',\bar P')$ is of type $A_1$.
In this case, $F'$ is smooth and so is $X'$ (see, e.g., \cite[Lemma
1.4]{Prokhorov-1997_e}). Therefore, $f$ is the quotient of a smooth
conic bundle by $\muu_d$. By Proposition \ref{prop-Gor-quot} we get
the case \ref{item=main--th-pr-ex3}.

Thus we may assume that $d=1$ (and $X'=X$). Let $R\subset Z$ be the
ramification divisor of $f_2$. Since $(\bar F,\bar P)$ is of type
$A$, in some coordinate system on $Z=\CC^2$, $R$ is given by the
equation $x^k+y^2=0$. Let $\Gamma\subset Z$ is given by $x=0$ and
let $S:=f^*\Gamma$. By the inversion of adjunction the log divisor
$K_Z+\Gamma+\frac 12R$ is log canonical (lc). So are $K_{\bar
F}+f_2^*\Gamma= f_2^*(K_Z+\Gamma+\frac 12R)$ and
$K_F+f^*\Gamma=K_F+S|_F$. Again by the inversion of adjunction
$K_X+F+S$ is lc near $F$. Shrinking $X$ we may assume that $K_X+F+S$
is lc everywhere. Replacing $\Gamma$ with a general hyperplane
section through $o$, we may assume that $S$ is smooth outside of
$C$. Then $K_X+S$ is plt. In particular, $S$ is normal and has only
log terminal singularities of type $T$ \cite{Kollar-ShB-1988}. Let
$D:=F|_S$. Then $K_S+D\sim 0$ is lc and $D\supset C$. By the
classification of two-dimensional log canonical singularities
\cite{Kawamata-1988-crep}, \cite[Ch. 3]{Utah} $K_S+C$ is plt.

The restriction $f_S\colon S\to \Gamma$ is a rational curve
fibration such that $-K_S$ is $f_S$-ample. If $C$ is a Cartier
divisor on $S$, then $S$ is smooth and so is $X$. Take the minimal
positive $n$ such that $nC$ is Cartier. Then $nC\sim 0$. This
induces an \'etale in codimension one $\muu_n$-cover $\pi \colon
S^{\natural} \to S$ such that $C^{\natural}:=\pi^*C\sim 0$. The
divisor $K_{S^{\natural}}+C^{\natural}=\pi^*(K_S+C)$ is plt (see,
e.g., \cite[20.3]{Utah}). Hence, $C^{\natural}$ is smooth and so is
$S^{\natural}$. Thus $S$ is a quotient of $S^{\natural}\simeq
\CC\times \PP^1$ by $\muu_n$. It is easy to see that $S$ has
singular points of types $\frac1n (1,q)$ and $\frac1n (-1,q)$, where
$\gcd (n,q)=1$. These points are of type $T$ if and only if
\[
(q+1)^2\equiv (q-1)^2 \equiv 0\mod n
\]
(see \cite{Kollar-ShB-1988}). This implies $n=2$ or $4$.
If $n=2$, then $S$ is Gorenstein and
so is $X$, a contradiction. Hence $n= 4$, so the
singularities of $S$ are of types $\frac14(1,1)$ and $A_3$. By
\cite{Kollar-ShB-1988} $(X,C)$ has exactly one non-Gorenstein point
which is of index $2$.
\end{proof}

\section{Index two $\QQ$-conic bundles}
Index two $\QQ$-conic bundles were classified in \cite[\S
3]{Prokhorov-1997_e}. Under the condition that the base $(Z,o)$ is
smooth, these are quotients of some elliptic fibrations by an
involution. Here we propose an alternative description and sketch a
different proof. (Note that a $\QQ$-conic bundle of index two over a
singular base is either of type \ref{item-main-th-impr-barm=1} or toroidal 
\cite[\S 3]{Prokhorov-1997_e}).

\begin{theoremm}
\label{th-index=2}
Let $f\colon (X,C)\to (Z,o)$ be a $\QQ$-conic bundle germ of index
two. Assume that $(Z,o)$ is smooth. Fix an isomorphism $(Z,o)\simeq
(\CC^2,0)$. Then there is an embedding
\begin{equation}
\label{eq-diag-last-2}
\xymatrix{X \ar@{^{(}->}[r] \ar[rd]_{f}& \PP(1,1,1,2)\times \CC^2
\ar[d]^{p}
\\
&\CC^2}
\end{equation}
such that $X$ is given by two equations
\begin{equation}
\label{eq-eq-index2}
\begin{array}{l}
q_1(y_1,y_2,y_3)-\psi_1(y_1,\dots,y_4;u,v)=0,
\\[7pt]
q_2(y_1,y_2,y_3)-\psi_2(y_1,\dots,y_4;u,v)=0,
\end{array}
\end{equation}
where $\psi_i$ and $q_i$ are weighted quadratic in $y_1,\dots,y_4$
with respect to $\wt (y_1,\dots,y_4)=(1,1,1,2)$ and
$\psi_i(y_1,\dots,y_4;0,0)=0$. The only non-Gorenstein point of $X$
is $(0,0,0,1; 0,0)$. Up to projective transformations, the
following are the possibilities for $q_1$ and $q_2$:

\begin{emptytheorem}
\label{cla-index-2-4}
$q_1 = y_1^2-y_2^2$ and $q_2 = y_1y_2-y_3^2$; then $f^{-1}(o)$ is
reduced and has exactly four irreducible components;
\end{emptytheorem}

\begin{emptytheorem}
$q_1 = y_1y_2$ and $q_2 = (y_1 +y_2)y_3$; then $f^{-1}(o)$ has three
irreducible components, one of them has multiplicity $2$;
\end{emptytheorem}
\begin{emptytheorem}
\label{cla-index-2-1-3}
$q_1 = y_1y_2 - y_3^2$ and $q_2 = y_1y_3$; then $f^{-1}(o)$ has two
irreducible components, one of them has multiplicity $3$;
\end{emptytheorem}

\begin{emptytheorem}
\label{cla-index-2-2-2}
$q_1 = y_1^2-y_2^2$ and $q_2 = y_3^2$; then $f^{-1}(o)$ has two
irreducible components, both of multiplicity $2$;
\end{emptytheorem}
\begin{emptytheorem}
\label{cla-index-4a}
$q_1 = y_1y_2 - y_3^2$ and $q_2 = y_1^2$; then $f^{-1}(o)$ is
irreducible of multiplicity $4$;
\end{emptytheorem}
\begin{emptytheorem}
\label{cla-index-4b}
$q_1 = y_1^2$ and $q_2 = y_2^2$; then $f^{-1}(o)$ is also
irreducible of multiplicity $4$.
\end{emptytheorem}
\par\noindent
Conversely, if $X\subset\PP(1,1,1,2)\times \CC^2$ is given by
equations of the form \eqref{eq-eq-index2} and singularities of $X$
are terminal, then the projection $f\colon (X,f^{-1}(0)_{\red}) \to
(\CC^2,0)$ is a $\QQ$-conic bundle of index $2$.
\end{theoremm}

\begin{proof}[Sketch of the proof]
First we prove the last statement. Assume that $X$ has only terminal
singularities. Then $X$ does not contain the surface
$\{y_1=y_2=y_3=0\} = \Sing \PP\times \CC^2$ (otherwise both $\psi_1$
and $\psi_2$ do not depend on $y_4$). By the adjunction formula,
$K_X=-L|_X$, where $L$ is a Weil divisor on $\PP\times \CC^2$ such
that the restriction $L|_{\PP}$ is $\OOO_{\PP}(1)$. Therefore, $X\to
\CC^2$ is a $\QQ$-conic bundle. It is easy to see that the only
non-Gorenstein point of $X$ is $(0,0,0,1;0,0)$ and it is of index
$2$.

Now let $f\colon (X,C)\to (Z,o)\simeq (\CC^2,0)$ be a $\QQ$-conic
bundle germ of index two. Let $P\in X$ be a point of index $2$. We
claim that $P$ is the only non-Gorenstein point. Indeed, if $C$ is
irreducible, the assertion follows by Corollary 
\ref{cor-cyclic-quo-new}, (iii).
If $C=\cup C_i$ is reducible, the same holds by Lemma
\ref{lemma-int-non-Gor} and \cite[Th. 4.2, Prop. 4.6]{Kollar-Mori-1992}. Thus
$P$ is the only non-Gorenstein point on $X$. By
\ref{prop-prim-types-def} each $(X,C_i)$ is of type \type{IA} at
$P$. 
(The case \type{IB} is excluded
by Proposition \ref{prop-IB} and
\cite[Th. 6.3]{Mori-1988}).
Hence the general member $F\in |-K_X|$ satisfies $F\cap
C=\{P\}$ and has only Du Val singularity at $P$ (see 
\cite[Th. 7.3]{Mori-1988}).

Let $\pi\colon (X^\sharp,P^\sharp)\to (X,P)$ be the 
index-one cover and let $F^\sharp=\pi^{-1} (F)_{\red}$ 
be the pull-back of $F$.
Let $\Gamma:=f^{-1}(o)$ be the scheme fiber and let 
$\Gamma^\sharp=\pi^{-1} (\Gamma)$.
\begin{lemma}
\label{lemma-12-fiber-new}
\[
\OOO_{F^\sharp \cap \Gamma^\sharp} \simeq \CC[x,y]/(xy,\, x^2+y^2).
\]
Furthermore $\muu_2$-action is given by $\wt(x,y) \equiv (1,1) \mod 2$.
\end{lemma}

\begin{proof}
Since $(F^\sharp,P^\sharp)$ is a Du Val singularity,
we may assume that $(F^\sharp,P^\sharp)\subset (\CC^3_{x,y,z},0)$.
The scheme $F^\sharp \cap \Gamma^\sharp$ is defined in 
$\CC^3_{x,y,z}$ by three equations 
$\alpha =\beta =\gamma =0$, where 
two of them are coordinates on $Z=\CC^2$, and the rest is 
the defining equation of $F^\sharp\subset \CC^3$. 
Since the morphism $F^\sharp \to Z$ is flat and of degree 4,
we have
\[
\OOO_{F^\sharp \cap \Gamma^\sharp} \simeq \CC\{x,y,z\}/
(\alpha, \beta,\gamma)
\]
is of length $4$. Furthermore $\muu_2$ acts on the ring so
that $\wt(x,y,z,\alpha, \beta,\gamma) \equiv (1,1,0,0,0,0) \mod 2$
because the quotient $(F^\sharp,P^\sharp)/\muu_2$
is Du Val, and, in particular, Gorenstein.
If $\alpha, \beta, \gamma \in (x,y,z)^2$, then 
$\len \CC\{x,y,z\}/(\alpha, \beta,\gamma) \ge 8$, which is a contradiction.
Hence, in view of the weights, we may assume that 
$\alpha =(\mbox{unit}) \cdot z + \alpha_1\cdot \beta + \alpha_2 \cdot \gamma $
modulo permutation of $\alpha, \beta, \gamma$. Thus we have
\[
\OOO_{F^\sharp \cap \Gamma^\sharp}\simeq \CC\{x,y\}/(\beta,\, \gamma).
\]
Since $\wt (x,y,\beta,\gamma) \equiv (1,1,0,0) \mod 2$, we see 
$\beta, \gamma \in (x,y)^2$.
Hence we may assume that $\beta \equiv xy \mod (x,y)^3$ modulo coordinate
change of $x,y$ and change of $\beta $, $\gamma $. 
Modulo analytic change of
coordinates $x$, $y$, we may assume 
$(\beta,\, \gamma)=(xy,\, x^a+y^b)$
for some $a,\, b \ge 2$. Since the ring 
$\OOO_{F^\sharp 
\cap \Gamma^\sharp}$ is of length $4$, we have
$4=a+b$ and hence $a=b=2$.
\end{proof}

Using this lemma one can apply 
arguments of \cite[pp. 631--633]{Mori-1975}
to get the desired embedding $X\subset \PP(1,1,1,2)\times Z$
considering the graded anti-canonical $\OOO_Z$-algebra
\[
\mathcal{R}:=\bigoplus_{i\ge 0}\mathcal{R}_i,\quad\mbox{where}
\quad \mathcal{R}_i:= H^0(\OOO_X(-iK_X)).
\]
We sketch the main idea.

Let $w$ be a local generator of $\OOO_{X^\sharp} (-K_{X})$ at
$P^\sharp$, let $u$, $v$ 
be coordinates on $Z=\CC^2$,
and let $z=0$ be the local equation of $F^\sharp$ in 
$(X^\sharp, P^\sharp)$.
Using the vanishing of $H^1(\OOO_X(-K_X))$ for $i>0$
and the exact sequence
\[
0 \to \OOO_X(-(i-1)K_X) \to \OOO_X(-iK_X)\to
\OOO_F(-iK_X)\to 0
\]
one can see 
\[
\mathcal{R}_i/(zw)\mathcal{R}_{i-1}\simeq H^0(\OOO_F(-iK_X)), \quad i>0.
\]
Therefore,
\[
\mathcal{R}_i/(zw)\mathcal{R}_{i-1}+(u,\, v)\mathcal{R}_i
=\bigl(\OOO_{F^\sharp\cap \Gamma^\sharp}(-iK_X)\bigr)^{\muu_2}.
\]
By Lemma \ref{lemma-12-fiber-new} we have an embedding
\[
\mathcal{R}/(zw,\, u,\, v)\mathcal{R}
\hookrightarrow
\left(\CC[x,\, y,\, w]/(xy,\, x^2+y^2)\right)^{\muu_2}.\]
Using $R_0/(u,v)R_0=\CC$, one can easily see that
\[
\mathcal{R}/(zw,\, u,\, v)\mathcal{R}=
\CC[y_1,y_2,y_4]/(y_1y_2,\, y_1^2+y_2^2),
\]
where $y_1=xw$, $y_2=yw$, $y_4=w^2$.
Put $y_3:=zw$. Then similarly to \cite[pp. 631--633]{Mori-1975}
we obtain
\[
\mathcal{R} \simeq 
\OOO_Z[y_1,y_2,y_3,y_4]/\mathcal I,
\]
where $\mathcal I$ is generated by the following regular
sequence
\[
\begin{array}{ll}
y_1y_2+y_3\ell_1(y_1,\dots,y_3)&+\psi_1(y_1,\dots,y_4;u,v),
\\[5pt]
y_1^2+y_2^2+y_3\ell_2(y_1,\dots,y_3)&+\psi_2(y_1,\dots,y_4;u,v)
\end{array}
\]
with $\psi_i(y_1,\dots,y_4;0,0)=0$.
\end{proof}

As is seen in Theorem \ref{th-main-1}, a $\QQ$-conic
bundle is often 
constructed as a
quotient of one of index two by a cyclic
group. 
Theorem \ref{th-index=2} is useful in  such
a context.
Finally we provide facts which 
are used in the study of 
$\QQ$-conic bundles 
with imprimitive points
(cf. Proposition \ref{prop-imprim-d2-}).

\begin{proposition}
\label{prop-cla-index-2-2-2-mu2-action-0}
Assume that in the notation of Theorem \xref{th-index=2} a cyclic
group $\muu_d$ acts on $X$ and $Z$ so that $f$ is
$\muu_d$-equivariant. Then the diagram \eqref{eq-diag-last-2} can be
chosen to be $\muu_d$-equivariant.
\end{proposition}

\begin{proof}
The sheaf $\OOO_X(-K_X)$ has a natural $\muu_d$-linearization. Hence
the embedding $X= \mt{Proj}\mathcal R\hookrightarrow \PP(1,1,1,2)$
is $\muu_d$-equivariant.
\end{proof}

The following is obvious.

\begin{lemma}
\label{lem-cla-index-2-2-2-mu2-action-l}
In notation of Theorem \xref{th-index=2} assume that 
$f$ has an equivariant $\muu_d$-action and that
$f^{-1}(o)$ has
two irreducible components, both of multiplicity $2$
\textup(i.e., we are in case
\xref{cla-index-2-2-2}\textup), which are permuted by
some element of $\muu_d$. 
Then the coordinates $y_1,\dots,y_4,u,v$ can
be chosen so that they and the equations \eqref{eq-eq-index2} are
semi-invariant.
\end{lemma}

\begin{proof}
Indeed, by \ref{prop-cla-index-2-2-2-mu2-action-0} the action of $\muu_d$
preserves the pencil $\lambda_1 (q_1-\psi_1)+\lambda_2(q_2-\psi_2)$.
It remains to note that in \ref{cla-index-2-2-2} $q_1$ and $q_2$ are
the only degenerate quadratic forms in this pencil and they cannot
be interchanged.
\end{proof}

\begin{lemma}
\label{lem-cla-index-2-2-2-mu2-action-2}
In notation and
assumptions of Theorem \xref{th-index=2} and Lemma
\xref{lem-cla-index-2-2-2-mu2-action-l} assume additionally that $d=2$.
Furthermore assume that $\muu_2$ acts on $X$ and $Z$ so that the
action is free in codimension one, has a unique fixed point $P=
(0,0,0,1;0,0)$, and the quotient $(X,P)/\muu_2$ is a terminal
singularity. Then modulo change of coordinates,
we are in case
\xref{cla-index-2-2-2}
with the action written as follows:
\[
y_1\mapsto y_1,\quad y_2\mapsto - y_2,\quad y_3\mapsto y_3,\quad y_4\mapsto
-y_4,\quad u\mapsto -u,\quad v\mapsto -v.
\]
\end{lemma}

\begin{proof}
By Lemma \ref{lem-cla-index-2-2-2-mu2-action-l} we can choose
the coordinates $y_1,\dots,y_4,u,v$ 
so that they and the equations \eqref{eq-eq-index2} are
semi-invariant.
Since the action of $\muu_2$ on $Z\simeq \CC^2$ is free outside of
$o$, this action is
given by $u\mapsto -u$, $v\mapsto -v$.
Modulo multiplication of $\pm1$ on the $\muu_2$-linearization
of $\OOO(-K_X)$, we may assume also that
$y_3\mapsto y_3$. Then $y_i\mapsto \pm y_i$ for all $i$.

Recall that we are in the case
\ref{cla-index-2-2-2} with $X \subset \PP(1,1,1,2)_{y_1,y_2,y_3,y_4}\times
\CC^2_{u,v}$.
Since $(1,1,0,0;0,0) \in X$, the point
is not $\muu_2$-fixed by the assumption. Hence
$y_1y_2 \mapsto -y_1y_2$. Modulo permutation of $y_1, y_2$,
we have $y_1 \mapsto y_1$ and $y_2 \mapsto -y_2$.
It remains to show only
that $y_4\mapsto - y_4$. Assume 
to the contrary that $y_4\mapsto y_4$.

Let $U\subset \PP(1,1,1,2)$ be the chart $y_4\neq 0$.
Then $U\simeq \CC^3_{z_1,z_2,z_3}/\muu_2(1,1,1)$.
Let $X^\sharp$ be the pull-back of $X\cap (U\times \CC^2_{u,v})$
on $\CC^3_{z_1,z_2,z_3}\times \CC^2_{u,v}$ and let $P^\sharp\in X^\sharp$
be the preimage of $P$.

Since the induced map $X^\sharp\to X$ is \'etale in codimension one,
$(X^\sharp, P^\sharp)\to (X,P)$ is the index-one cover.
Hence $(X^\sharp, P^\sharp)\to (X,P)/\muu_2$ is also the index-one cover
of
the terminal point
$(X,P)/\muu_2$ of index $4$
(the last is true because
the action of $\muu_2$ is free in codimension one).
Hence the morphism is a $\muu_4$-covering by the
structure of terminal  singularities.
However $(X,P)/\muu_2$ is the quotient of $(X^\sharp,P^\sharp)$ by
commuting $\muu_2$-actions:
\[
(z_1,z_2,z_3,u,v)\mapsto(-z_1,-z_2,-z_3,u,v),
(z_1,-z_2,z_3,-u,-v)
\]
This is a contradiction, and we have $y_4 \mapsto -y_4$ as
claimed.
\end{proof}

\subsection*{Acknowledgments}
The work was carried out at Research Institute for Mathematical
Sciences (RIMS), Kyoto University. The second author would like to
thank RIMS for invitations to work there in 2005-2006, for
hospitality and wonderful conditions of work.

The research of the first author was supported in part by JSPS Grant-in-Aid
for Scientific Research (B)(2), No. 16340004. The second author was
partially supported by Grant CRDF-RUM1-2692-MO-05.

Finally we are extremely
grateful to the referee,
who has pointed out so
many mistakes and 
helped us improve the
presentation.


\begin{thebibliography}{Mum66}

\bibitem[Cat87]{Catanese-1987}
Catanese, F.,
Automorphisms of rational double points and moduli spaces of surfaces
of general type,
\textit{Compositio Math.}, \textbf{61} (1987), 81--102.

\bibitem[Cut88]{Cutkosky-1988}
Cutkosky, S.,
Elementary contractions of {G}orenstein threefolds,
\textit{Math. Ann.}, \textbf{280} (1988), 521--525.

\bibitem[Feu71]{Feuer-1971}
Feuer, R.~D.,
Torsion-free subgroups of triangle groups,
\textit{Proc. Amer. Math. Soc.}, \textbf{30} (1971), 235--240.

\bibitem[Isk96]{Iskovskikh-1996-conic-re}
Iskovskikh, V.~A.,
On a rationality criterion for conic bundles,
\textit{Mat. Sb.}, \textbf{187} (1996),
75--92.

\bibitem[Kaw88]{Kawamata-1988-crep}
Kawamata, Y.,
Crepant blowing-up of {$3$}-dimensional 
canonical singularities and
its application to degenerations of surfaces,
\textit{Ann. of Math.}, \textbf{127} (1988),
93--163.

\bibitem[KM92]{Kollar-Mori-1992}
Koll{\'a}r, J. and Mori, S.,
Classification of three-dimensional flips,
\textit{J. Amer. Math. Soc.}, \textbf{5} (1992), 533--703.

\bibitem[KM98]{KM-1998}
Koll{\'a}r, J. and Mori, S.,
\textit{Birational geometry of algebraic varieties}, 
Cambridge University Press, 1998.

\bibitem[Kol86]{Kollar-1986-I}
Koll{\'a}r, J.,
Higher direct images of dualizing sheaves, {I},
\textit{Ann. of Math.}, \textbf{123} (1986), 11--42.

\bibitem[Kol92]{Utah}
Koll{\'a}r, J. et al., 
\textit{Flips and abundance for algebraic threefolds},
Soci\'et\'e Math\'ematique de France, Paris, 1992.
Papers from the Second Summer Seminar on 
Algebraic Geometry held at
the University of Utah, Salt Lake City, Utah, August 1991, Ast\'erisque \textbf{211}, (1992).

\bibitem[Kol99]{Kollar-1999-R}
Koll{\'a}r, J.,
Real algebraic threefolds, {III}, {C}onic bundles,
\textit{J. Math. Sci. (New York)}, \textbf{94}
(1999), 996--1020.


\bibitem[KSB88]{Kollar-ShB-1988}
Koll{\'a}r, J. and Shepherd-Barron, N.~I., 
Threefolds and deformations of surface singularities,
\textit{Invent. Math.}, \textbf{91} (1988), 299--338.

\bibitem[Mor75]{Mori-1975}
Mori, S.,
On a generalization of complete intersections,
\textit{J. Math. Kyoto Univ.}, \textbf{15} (1975), 619--646.

\bibitem[Mor88]{Mori-1988}
Mori, S.,
Flip theorem and the existence of minimal models for {$3$}-folds,
\textit{J. Amer. Math. Soc.}, \textbf{1} (1988), 117--253.

\bibitem[Mum66]{Mumford-Lectures-on-curves}
Mumford, D.,
\textit{Lectures on curves on an algebraic surface},
Princeton University Press, 1966.

\bibitem[Pro97a]{Prokhorov-1997_e}
Prokhorov, Yu.,
On the complementability of the canonical divisor for {M}ori
fibrations on conics,
\textit{Sbornik. Math.}, \textbf{188} (1997), 1665--1685.

\bibitem[Pro97b]{Prokhorov-1996b}
Prokhorov, Yu.,
On extremal contractions from threefolds to surfaces: the case of one
non-{G}orenstein point,
\textit{Birational algebraic geometry (Baltimore, MD, 1996)}, \textit{Contemp. Math.}, \textbf{207} (1997), 119--141. 

\bibitem[Spr77]{Springer-1977}
Springer, T.~A.,
\textit{Invariant theory},
\textit{Lecture Notes in Mathematics}, \textbf{585}
(Springer-Verlag, Berlin, 1977).

\bibitem[Ste88]{Stevens-1988}
Stevens, J.,
On canonical singularities as total spaces of deformations,
\textit{Abh. Math. Sem. Univ. Hamburg}, \textbf{58} (1988), 275--283.

\end{thebibliography}
\end{document}